\documentclass[12pt, leqno]{amsart}
\usepackage{
amsmath,
amssymb,
amsthm,
mathrsfs,
amsfonts,
ytableau,
enumitem,
comment,
upgreek,
soul,
yhmath,
mathtools,
kotex,
stackengine,
caption,
tabularx,
stackrel,
amscd,
bbm}
\usepackage[mathscr]{eucal}
\usepackage{tikz,tikz-cd}
\usetikzlibrary{decorations.pathreplacing,decorations.pathmorphing}
\usetikzlibrary{arrows,arrows.meta}
\usetikzlibrary{shapes,shapes.geometric,shapes.symbols}
\usetikzlibrary{matrix,calc}
\usepackage[poly,all]{xy}
\usepackage[colorinlistoftodos, textwidth = 2.3cm]{todonotes}
\usepackage{graphicx}
\usepackage{xcolor}
\usepackage{hyperref}
\hypersetup{
    colorlinks=true,
    linkcolor=blue,
    citecolor=magenta,
    urlcolor=cyan
}
\usepackage[capitalise, nameinlink]{cleveref}
\ytableausetup{mathmode, boxsize=1.1em,aligntableaux=center}

\crefname{figure}{{\sc Figure}}{{\sc Figure}}
\crefname{section}{Section}{Sections}
\crefname{example}{Example}{Examples}
\crefname{theorem}{Theorem}{Theorems}
\crefname{lemma}{Lemma}{Lemmas}
\crefname{proposition}{Proposition}{Propositions}
\crefname{thm}{Theorem}{Theorems}
\crefname{lem}{Lemma}{Lemmas}
\crefname{prop}{Proposition}{Propositions}
\crefname{figure}{Figure}{Figures}
\crefname{fig}{Figure}{Figures}
\crefname{remark}{Remark}{Remarks}
\crefname{rem}{Remark}{Remarks}
\crefname{cor}{Corollary}{Corollaries}
\crefname{corollary}{Corollary}{Corollaries}
\crefname{conjecture}{Conjecture}{Conjectures}
\crefname{conj}{Conjecture}{Conjectures}
\crefname{ex}{Example}{Examples}

\newtheorem{theorem}{Theorem}[section]
\newtheorem{proposition}[theorem]{Proposition}
\newtheorem{lemma}[theorem]{Lemma}
\newtheorem{corollary}[theorem]{Corollary}

\newtheorem*{claim*}{Claim}

\theoremstyle{definition}
\newtheorem{algorithm}[theorem]{Algorithm}
\newtheorem{example}[theorem]{Example}
\newtheorem{definition}[theorem]{Definition}
\newtheorem{summary}[theorem]{Summary}
\newtheorem{remark}[theorem]{Remark}

\numberwithin{equation}{section} \numberwithin{figure}{section}
\numberwithin{table}{section}

\def \hp {0.35}
\def \vp {0.45}

\def \Z {\mathbb Z}
\def \C {\mathbb C}

\newcommand{\nc}{\newcommand}
\nc{\tcd}{\mathtt{cd}}
\nc{\QSQ}[1]{{\widetilde{Q}}_{#1}}
\nc{\PYQS}[1]{\widetilde{\mathscr{S}}_{#1}}
\nc{\YQS}[1]{\ensuremath{\hat{\mathscr{S}}_{#1}}}
\nc{\QS}[1]{\mathscr{S}_{#1}}
\nc{\DIF}[1]{\ensuremath{\mathfrak{S}_{#1}^*}}
\nc{\calT}{\mathscr{T}}
\nc{\comp}{\ensuremath{\mathrm{comp}}}
\nc{\Sym}{\ensuremath{\mathsf{Sym}}}
\nc{\QSym}{\ensuremath{\mathsf{QSym}}}
\nc{\sort}{\ensuremath{\mathtt{sort}}}
\nc{\rect}{\ensuremath{\mathtt{rect}}}
\nc{\rw}{\ensuremath{\mathrm{rw}}}
\nc{\rwRC}[2]{\ensuremath{\mathrm{w}_{\rm #1}(#2)}}
\nc{\stan}{\ensuremath{\mathsf{stan}}}
\nc{\Peak}{\ensuremath{\mathrm{Peak}}}
\nc{\Des}{\ensuremath{\mathrm{Des}}}
\nc{\BDes}{\ensuremath{\mathbf{Des}}}
\nc{\Asc}{\ensuremath{\mathrm{Asc}}}
\nc{\SYT}{\ensuremath{\mathsf{SYT}}}
\nc{\SRYT}{\ensuremath{\mathsf{SRYT}}}
\nc{\SIT}{\ensuremath{\mathsf{SIT}}}
\nc{\DIRT}{\ensuremath{\mathsf{DIRT}}}
\nc{\SPIT}{\ensuremath{\mathsf{SPIT}}}
\nc{\SYCT}{\ensuremath{\mathsf{SYCT}}}
\nc{\SPYCT}{\ensuremath{\mathsf{SPYCT}}}
\nc{\SShT}{\ensuremath{\mathsf{SShT}}}
\nc{\SPYT}{\ensuremath{\mathsf{SPYT}}}
\nc{\SRCT}{\ensuremath{\mathsf{SRCT}}}
\nc{\SVT}{\ensuremath{\mathsf{SVT}}}
\nc{\SG}{\ensuremath{\mathfrak{S}}}
\nc{\calV}{\ensuremath{\mathcal{V}}}
\nc{\PCp}{\mathsf{Pcomp}}
\nc{\Cp}{\mathsf{Comp}}
\nc{\Par}{\mathsf{Par}}
\nc{\set}{\mathrm{set}}
\nc{\rmr}{\mathrm{r}}
\nc{\rmc}{\mathrm{c}}
\nc{\Qsym}{\mathrm{QSym}}
\nc{\inv}[1]{\mathrm{Inv}(#1)}
\nc{\coinv}[1]{\mathrm{Coinv}(#1)}
\nc{\calB}{\ensuremath{\mathcal{B}}}
\nc{\calR}{\ensuremath{\mathcal{R}}}
\nc{\setA}[1]{\mathsf{A}(#1)}
\nc{\setB}[1]{\mathsf{B}(#1)}
\nc{\HnZmodC}{\ensuremath{\mathbf{Q}^{\rm c}}}
\nc{\HnZmodR}{\ensuremath{\mathbf{Q}^{\rm r}}}
\nc{\ch}{\mathrm{ch}}
\nc{\calG}{\ensuremath{\mathcal{G}}}
\nc{\bfF}{\mathbf{F}}
\nc{\bfd}{\mathbf{d}}
\nc{\bfu}{\mathbf{u}}
\nc{\opi}{\overline{\pi}}
\nc{\RWC}{\ensuremath{{\mathsf {RW}}_{\rm c}}}
\nc{\RWR}{\ensuremath{{\mathsf {RW}}_{\rm r}}}
\nc{\rev}{\ensuremath{{\textrm{r}}}}
\nc{\DPsecond}[1]{\ensuremath{\mathsf{Dyck}^{1/2}(#1)}}
\nc{\readSYCT}{\ensuremath{\mathrm{w}_{\rm Y}}}
\nc{\SPRCT}{\ensuremath{\mathsf{SPRCT}}}
\nc{\alphamax}{\ensuremath{\alpha_{\rm max}}}
\nc{\subtableauT}[2]{\ensuremath{{#1}^{\rm sub}_{#2}}}
\nc{\sh}{\ensuremath{\textrm{sh}}}
\nc{\bfs}{\ensuremath{\mathbf{s}}}
\nc{\SPITtoYT}[1]{\ensuremath{{\rm sort}(#1)}}
\nc{\Tsup}[1]{T^{\rm sup}_{#1}}
\nc{\ourmapS}{\Xi_\alpha}
\nc{\ourmapSS}{\Theta_\alpha}
\nc{\ourc}{\mathfrak{c}_\alpha}
\nc{\id}{\mathrm{id}}
\nc{\mergef}{\ensuremath{\mathfrak{m}}}
\nc{\splitf}{\ensuremath{\mathfrak{s}}}
\nc{\reflectf}{\ensuremath{\mathfrak{r}}}
\nc{\AHM}{\underline{\textrm{AHM}}}
\nc{\setWord}[1]{\ensuremath{{\sf W}_{#1}}}
\nc{\Tsupalpha}{T^{\rm sup}_\alpha}
\nc{\Tsourcealpha}{T^{\rm source}_\alpha}
\nc{\Tsinkalpha}{T^{\rm sink}_\alpha}
\nc{\Tspalpha}{T^{\rm sp}_\alpha}
\nc{\modL}{H_n(0)\text{-\bf mod}}
\nc{\modLb}{H_\bullet(0)\text{-\bf mod}}
\nc{\osfB}{\overline{\mathsf{B}}}
\nc{\sfB}{\mathsf{B}}
\nc{\Downseq}[2]{{\bf d}_{#1;#2}}
\nc{\Upseq}[2]{{\bf u}_{#1;#2}}
\nc{\actionS}{\scalebox{0.6}{\ $\bullet$ }}
\newcommand{\oeis}[1]{\href{http://oeis.org/#1}{#1}}

\newlength\cellsize 

\savebox2{%
\begin{picture}(12,12)
\put(0,0){\line(1,0){12}}
\put(0,0){\line(0,1){12}}
\put(12,0){\line(0,1){12}}
\put(0,12){\line(1,0){12}}
\end{picture}}

\nc\cellify[1]{\def\thearg{#1}\def\nothing{}%
\ifx\thearg\nothing\vrule width0pt height\cellsize depth0pt%
  \else\hbox to 0pt{\usebox2\hss}\fi%
  \vbox to 12\unitlength{\vss\hbox to 12\unitlength{\hss$#1$\hss}\vss}}
  
\nc\tableau[1]{\vtop{\let\\=\cr
\setlength\baselineskip{-12000pt}
\setlength\lineskiplimit{12000pt}
\setlength\lineskip{0pt}
\halign{&\cellify{##}\cr#1\crcr}}}

\nc\ctab[2]{\ensuremath{{#1} = \begin{array}{c}
\tableau{#2}
\end{array}
}}
\nc\ctabT[1]{\ensuremath{\begin{array}{c}
\tableau{#1}
\end{array}
}}
\nc{\DYCKPsix}[8]{
\begin{tikzpicture}[baseline=5mm]
\def \hp {5mm}
\def \vp {5mm}
\foreach \c in {0,1,2,3,4}{
    \draw[-,black!50] (\hp*\c,\vp*0) -- (\hp*\c,\vp*2);
}
\foreach \c in {0,1,2}{
    \draw[-,black!50] (\hp*0,\vp*\c) -- (\hp*4,\vp*\c);
}
\draw[dotted] (0,0) -- (\hp*4,\vp*2);
\draw[-,thick] (\hp*0,\vp*0) -- (\hp*1,\vp*0) -- (\hp*#1,\vp*#2) -- (\hp*#3,\vp*#4) -- (\hp*#5,\vp*#6) -- (\hp*#7,\vp*#8) -- (\hp*4,\vp*2);
\end{tikzpicture}
}

\nc{\nt}[1]{\todo[size=\tiny,color=red!40]{#1 \\ \hfill -- Note}}

\title[Quasisymmetric Schur $Q$-functions]
{Quasisymmetric Schur $Q$-functions and peak Young quasisymmetric Schur functions}

\author[S.-I. Choi]{Seung-Il Choi}
\address{Center for Quantum structures in Modules and Spaces, Seoul National University, Seoul 08826, Republic of Korea}
\email{ignatioschoi@snu.ac.kr}

\author[S.-Y. Nam]{Sun-Young Nam}
\address{Department of Mathematics, Sogang University, Seoul 04107, Republic of Korea}
\email{synam.math@gmail.com}

\author[Y.-T. Oh]{Young-Tak Oh}
\address{Department of Mathematics / Institute for Mathematical and Data Sciences, Sogang University, Seoul 04107, Republic of Korea}
\email{ytoh@sogang.ac.kr}

\keywords{Quasisymmetric Schur $Q$-function, Peak Young quasisymmetric function, $0$-Hecke algebra, Weak Bruhat order, Fuss--Catalan number}

\date{\today}
\subjclass[2020]{20C08, 05E05, 05E10}

\begin{document}
\begin{abstract}
In this paper, we explore the relationship between quasisymmetric Schur $Q$-functions and peak Young quasisymmetric Schur functions.
We introduce a bijection on $\mathsf{SPIT}(\alpha)$ such that $\{\mathrm{w}_{\rm c}(T) \mid T \in \mathsf{SPIT}(\alpha)\}$ and $\{\mathrm{w}_{\rm r}(T) \mid T \in \mathsf{SPIT}(\alpha)\}$ share identical descent distributions. Here, $\mathsf{SPIT}(\alpha)$ is the set of standard peak immaculate tableaux of shape $\alpha$, and $\mathrm{w}_{\rm c}$ and $\mathrm{w}_{\rm r}$ denote column reading and row reading, respectively.
By combining this equidistribution with the algorithm developed by Allen, Hallam, and Mason, we demonstrate that the transition matrix from the basis of quasisymmetric Schur $Q$-functions to the basis of peak Young quasisymmetric Schur functions is upper triangular, with entries being non-negative integers.
Furthermore, we provide explicit descriptions of the expansion of peak Young quasisymmetric Schur functions in specific cases, in terms of quasisymmetric Schur $Q$-functions.
We also investigate the combinatorial properties of standard peak immaculate tableaux, standard Young composition tableaux, and standard peak Young composition tableaux.
We provide a hook length formula for $\mathsf{SPIT}(\alpha)$ and show that standard Young composition tableaux and standard peak Young composition tableaux can be each bijectively mapped to words satisfying suitable conditions.
Especially, cases of compositions with rectangular shape are examined in detail.
\end{abstract}

\maketitle
  
\tableofcontents

\section{Introduction}
Quasisymmetric Schur functions and dual immaculate quasisymmetric functions emerged in the early 2010s as significant quasi-analogues of Schur functions. 
Specifically, quasisymmetric Schur functions were introduced by Haglund, Luoto, Mason, and van Willigenburg \cite{11HLMW} as specializations of integral form nonsymmetric Macdonald polynomials obtained by setting $q=t=0$ in the combinatorial formula described in \cite{09Mason} and summing the resulting Demazure atoms over all weak compositions that collapse to the same strong composition. 
On the other hand, dual immaculate quasisymmetric functions were introduced by Berg, Bergeron, Saliola, Serrano, and Zabrocki \cite{14BBSSZ} as the quasisymmetric duals of the corresponding immaculate noncommutative symmetric functions.
These functions respectively form a basis for the ring $\QSym$ of quasisymmetric functions and possess many properties similar to those of Schur functions (for instance, see \cite{11BLW,11HLMW2,13LMvW} and \cite{17BBSSZ, 16BSZ, 16BTvW}).
They have also been extensively investigated within the framework of $0$-Hecke algebras (see \cite{15TW,15BBSSZ, 19TW,21CKNO,22CKNO1,24CKO}).

At the end of the 2010s, it was demonstrated that dual immaculate quasisymmetric functions have a closer relationship with Young quasisymmetric Schur functions than with quasisymmetric Schur functions. 
Young quasisymmetric Schur functions were introduced by Luoto, Mykytiuk, and van Willigenburg in \cite{13LMvW} as a variant of quasisymmetric Schur functions (for the precise definition, see \cref{The ring of quasisymmetric functions and its bases}).
For a composition $\alpha$, let $\DIF{\alpha}$ and $\YQS{\alpha}$ denote the dual immaculate quasisymmetric function and the Young quasisymmetric Schur function associated with $\alpha$.
Allen, Hallam, and Mason \cite{18AHM} showed that 
\begin{equation}\label{expansion of DIMF into YQS}
\DIF{\alpha} = \sum_{\beta} c_{\alpha,\beta} \YQS{\beta},
\end{equation}
where $c_{\alpha, \beta}$ counts the number of dual immaculate recording tableaux of shape $\beta$ with row strip shape $\alpha^\rmr$, and the transition matrix from the basis $\{\DIF{\alpha}\}$ to $\{\YQS{\alpha}\}$ is unitriangular. For more details, refer to \cref{A generating function-based interpretation}.

Let us recall the $F$-expansions of $\DIF{\alpha}$ and $\YQS{\alpha}$, which play a fundamental role throughout this paper.
It was shown in \cite[Proposition 3.48]{14BBSSZ} that  
$\DIF{\alpha}=\sum_{T}F_{\comp(\Des(\rwRC{r}{T}))}$,
where the summation is over the set of standard immaculate tableaux of shape $\alpha$, $\rwRC{r}{T}$ is the word obtained by reading the entries from left to right, starting from the uppermost row of $T$, and 
$F_{\gamma}$ is the {\em fundamental quasisymmetric function} associated with $\gamma$.
Similarly, by combining \cite[Proposition 5.2.2]{13LMvW} with \cite[Lemma 5.3]{22Searles}, one derives that  
$\YQS{\alpha} = \sum_{T}F_{\comp(\Des(\readSYCT(T)))}$,
where the summation is over the set of standard Young composition tableaux of shape $\alpha$
and $\readSYCT(T)$ is the word obtained from $T$ by reading its entries from top to bottom in the first column, and then bottom to top in subsequent columns, reading the columns in order from left to right.

Now, we introduce the main subjects of this paper.
Jing and Li \cite{15JL} introduced {\em quasisymmetric Schur $Q$-functions} by dualizing noncommutative lifts of Schur $Q$-functions. 
Meanwhile, Searles \cite{22Searles} introduced {\em peak Young quasisymmetric Schur functions} as the images of certain $0$-Hecke--Clifford modules under the peak quasisymmetric characteristic.
For a peak composition $\alpha$, 
let $\QSQ{\alpha}$ denote the quasisymmetric Schur $Q$-function 
and $\PYQS{\alpha}$ denote the peak Young quasisymmetric Schur function $\PYQS{\alpha}$ associated with $\alpha$.  
The $K$-expansions of $\QSQ{\alpha}$ and $\PYQS{\alpha}$ are given in \cite[Theorem 3.9]{19Oguz} and \cite[Theorem 5.9]{22Searles}, respectively. To be more precise,  
\begin{equation}\label{K-expansions of quasi Schur Q and Young quasiSchur}
\QSQ{\alpha} = \sum_{T \in \SPIT(\alpha)} K_{\comp(\Peak(\Des(\rwRC{c}{T})))} \quad \text{and} \quad \PYQS{\alpha} = \sum_{T \in \SPYCT(\alpha)} K_{\comp(\Peak(\Des(\readSYCT(T))))},
\end{equation}
where $\SPIT(\alpha)$ is the set  of standard peak immaculate tableaux of shape $\alpha$, $\SPYCT(\alpha)$ is the set of standard peak Young composition tableaux of shape $\alpha$, 
$\rwRC{c}{T}$ is the word obtained by reading the entries from bottom to top, starting from the leftmost column of $T$, and $K_{\gamma}$ is the {\em peak quasisymmetric function} associated with the peak composition $\gamma$. 
By comparing these expansions with the $F$-expansions of $\DIF{\alpha}$ and $\YQS{\alpha}$, we can regard $\QSQ{\alpha}$ and $\PYQS{\alpha}$ as natural peak analogues of $\DIF{\alpha}$ and $\YQS{\alpha}$.
However, it should be noted that in the right-hand side of the first equality in \eqref{K-expansions of quasi Schur Q and Young quasiSchur}, the column reading word $\rwRC{c}{T}$ is used instead of the row reading word $\rwRC{r}{T}$.

In this paper, we examine the relationship between quasisymmetric Schur $Q$-functions and peak Young quasisymmetric Schur functions by studying their combinatorial models.

Given a positive integer $n$, the classification of subsets of $\SG_n$ up to descent-preserving isomorphism is known to be a highly challenging problem
(for instance, see \cite{72Stanley,14MW,23KLO}).
In \cref{A combinatorial bijection on SPIT preserving its descent}, we introduce a function $\ourmapS: \SPIT(\alpha) \to \SPIT(\alpha)$, which is constructed using an algorithm operating on $\SPIT(\alpha)$. 
We confirm the well-definedness of this function and demonstrate that it satisfies the equality:
\[
\Des(\rwRC{c}{T}) = \Des(\rwRC{r}{\ourmapS(T)}) \quad \text{for all $T \in \SPIT(\alpha)$} 
\]
(\cref{Thm3-1: bijection1 on SPITs}).
Furthermore, we establish the bijectivity of this function (\cref{Thm3-2: our map is bij.}).
Consequently, we conclude that $\{\rwRC{c}{T} \mid T \in \SPIT(\alpha)\}$ and $\{\rwRC{r}{T} \mid T \in \SPIT(\alpha)\}$ share the same descent distribution (\cref{Coro: new definition of QSQ}).
From this equidistribution, we obtain the equality:
\begin{equation}\label{Eq: new def. of QSQ}
\QSQ{\alpha} = \sum_{T \in \SPIT(\alpha)} K_{\comp(\Peak(\Des(\rwRC{r}{T})))}.
\end{equation}
It is worth noting that the algorithm by Allen, Hallam, and Mason \cite{18AHM} is designed for row-reading words. Thanks to this equality, in \cref{A generating function-based interpretation}, we can utilize their algorithm.
In \cref{A representation-theoretical interpretation}, 
we investigate two $0$-Hecke modules 
$\HnZmodC_{\alpha}$ and $\HnZmodR_{\alpha}$ satisfying that  
\[
\ch(\HnZmodC_{\alpha}) = \sum_{\sigma \in \{\rwRC{c}{T} \mid T \in \SPIT(\alpha)\}} F_{\comp(\Des(\sigma))} \quad \text{and} \quad
\ch(\HnZmodR_{\alpha}) = \sum_{\sigma \in \{\rwRC{r}{T} \mid T \in \SPIT(\alpha)\}} F_{\comp(\Des(\sigma))},
\]
where $\ch$ is the quasisymmetric characteristic (see \eqref{quasi characteristic}).
The equidistribution mentioned above indicates that these two modules share the same image under the quasisymmetric characteristic.
We demonstrate that if $\alpha$ has at most one part greater than $2$, then $\HnZmodC_{\alpha}=\HnZmodR_{\alpha}$ as $0$-Hecke modules, and if $\alpha$ has more than one part greater than $2$, then $\HnZmodC_{\alpha}\not\cong\HnZmodR_{\alpha}$ as $0$-Hecke modules. 
Furthermore, we show that when $\alpha$ has at most one part greater than $2$, both $\{\rwRC{c}{T} \mid T \in \SPIT(\alpha)\}$ and $\{\rwRC{c}{T} \mid T \in \SPIT(\alpha)\}$ are left weak Bruhat intervals, and there exists a descent-preserving colored digraph between them (\cref{Thm: CRWC is interval if alpha one part2}).

In \cref{A generating function-based interpretation}, we resolve the problem proposed by Searles \cite[Section 8]{22Searles} by combining \eqref{Eq: new def. of QSQ} with the algorithm developed by Allen, Hallam, and Mason.
To be more precise, we prove that for a peak composition $\alpha$ of $n$,
\[
\QSQ{\alpha} = \sum_{\beta} c_{\alpha,\beta} \PYQS{\beta} 
\]
where the coefficients $c_{\alpha,\beta}$ are given in \eqref{expansion of DIMF into YQS}
(\cref{Thm: QSQ into PYQS}).
From this expansion, we can derive that the transition matrix from $\{\QSQ{\alpha}\}$ to $\{\PYQS{\alpha}\}$ is upper triangular, composed of nonnegative integer entries.
This result can be viewed as a peak version of \cite[Theorem 1.1]{18AHM}.
On the other hand, the expansion of a peak Young quasisymmetric Schur function in terms of quasisymmetric Schur $Q$-functions has not been well understood until now. 
In \cref{The expansion of PYQS in the basis QSQ in special cases}, we successfully describe the expansions of peak Young quasisymmetric Schur functions associated with strict partitions, peak compositions with entries $\leq 3$, and peak compositions containing at most one part greater than 2 (\cref{peak analogue:strict partition}, \cref{Thm: peak Young Quasi Schur to quasi Schur Q}, and \cref{Prop: peak hook shape1}). Additionally, we discover a characterization of peak compositions $\alpha$ for which $\PYQS{\alpha} = \QSQ{\alpha}$ (\cref{characterization:PYQS=QSD}).
These results are of significant importance in partially resolving the peak version conjectures of \cite[Conjecture 6.1 and Conjecture 6.2]{18AHM}.
    
The final section is devoted to exploring the combinatorial properties of $\SPIT(\alpha)$ and $\SPYCT(\alpha)$. 
In \cref{A hook length formula}, we present a hook length formula for $\SPIT(\alpha)$ that applies to every peak composition (\cref{hook length formula for SPIT}). 
However, closed-form enumeration formulas for standard Young composition tableaux and standard peak Young composition tableaux have yet to be discovered.
In \cref{SPYCT for rectangular shaped peak compositions}, we demonstrate that standard Young composition tableaux and standard peak Young composition tableaux can be each bijectively mapped to words
satisfying suitable conditions
(\cref{words for SYCT} and \cref{prop: word model for spyct}). 
The case where $\alpha$ is of rectangular shape is investigated in more detail.
Particularly intriguing is the fact that $|\SPYCT((3^k))|$ emerges as the well-known Fuss--Catalan number $\frac{1}{2k+1}\binom{3k}{k}$ (\cref{Example: 3^k and 4^k}). 
Therefore, $|\SPYCT((m^k))|$ can be regarded as a new generalization of Catalan numbers.
In \cref{Other bijections}, we introduce two sets $\SPYT((m^k))$ and $\SVT((k^{m-1}), \rho_m)$, which are in bijection with $\SPYCT((m^k))$ for rectangular-shaped peak compositions $(m^k)$ with $m\geq 2$ and $k\geq 1$ (\cref{Prop: bijections among SPYCT SPYT and SVT }).

\bigskip
\noindent
{\bf Convention.}
\ 
In this paper, we use many notations
related to the fillings of composition diagrams 
(see \cref{Compositions and peak compositions}). 
To help readers' comprehension, we collect them below.  

Let $\alpha$ be a composition and $T$ be a filling of the composition diagram of $\alpha$. 
\begin{itemize}
\item $T_{i, j}$ : the entry at the box $(i, j) \in \tcd(\alpha)$, i.e., the entry at column $i$ and row $j$ 

\item $T_{(i)}$ : the $i$th row from bottom to top

\item $T_{(\leq i)}$ : the subfilling consisting of the first $i$ rows

\item $T^{(j)}$ : the $j$th column from left to right

\item $T^{(\leq j)}$ : the subfilling consisting of the first $j$ columns

\item $T[\leq k]$ : the subfilling of $T$ consisting of the boxes with entries $\le k$

\item $|T|$ : the number of boxes in $T$

\item $T_{[1:j]}$ : the subfilling of $T$ whose composition diagram is composed of the first $j$ boxes from bottom to top and left to right
\end{itemize}
We also employ the notation $x \in T$ to denote that $x$ appears as an entry of $T$.

\section{Preliminaries}
\label{Sec: Prelim.}

For nonnegative integers $a, b$ with $a \leq b$, we define $[a, b]$ as the set $\{a, a+1, \dots, b\}$ and define $[a]$ as $[1, a]$ if $a \ge 1$, and as the empty set otherwise. 
Throughout this paper, $n$ denotes a nonnegative integer unless explicitly stated otherwise.

\subsection{Compositions and peak compositions}
\label{Compositions and peak compositions}

A \emph{composition} $\alpha$ of $n$ is a finite ordered list of positive integers $(\alpha_1, \alpha_2, \ldots, \alpha_k)$ such that $|\alpha|:=\sum_{i=1}^k \alpha_i = n$.
We call $k$ the \emph{length} of $\alpha$ and denote it by $\ell(\alpha)$.
For convenience, we regard the empty composition $\emptyset$ as a unique composition of size and length $0$.
We represent each composition visually using its composition diagram.
The \emph{composition diagram} $\tcd(\alpha)$ of $\alpha$ is the left-justified array of boxes with $\alpha_i$ boxes in row $i$, where the rows are numbered from bottom to top. 
In this paper, when there is no danger of confusion, we frequently identify each box in $\tcd(\alpha)$ by its $xy$-plane coordinate.
For example, when $\alpha = (3,2,5)$, the coordinate $(5,3)$ denotes the shaded box in the diagram
\[
\tcd(\alpha) =
\begin{ytableau}
~ & ~ & ~ & ~ & *(black!40) \\
~ & ~ \\
~ & ~ & ~
\end{ytableau} \ .
\]
For later use, we denote the largest part of $\alpha$ as $\alphamax$, and represent the number of boxes in the $i$th column of $\tcd(\alpha)$ as $c_i(\alpha)$.
Let $\Cp_n$ be the set of compositions of $n$,
and let $\Cp:=\bigcup\limits_{n\geq 0}\Cp_n$.

For $\alpha \in \Cp_n$ and $I = \{i_1 < i_2 < \cdots < i_p\} \subset [n-1]$,
define
\begin{align*}
\set(\alpha)&:= \{\alpha_1,\alpha_1+\alpha_2,\ldots, \alpha_1 + \alpha_2 + \cdots + \alpha_{\ell(\alpha)-1}\} \\
\comp(I) &:= (i_1,i_2 - i_1,\ldots,n-i_p).
\end{align*}
Then $\Cp_n$ is in bijection with the set of subsets of $[n-1]$ under the correspondence $\alpha \mapsto \set(\alpha)$ (or $I \mapsto \comp(I)$).
For $\alpha, \beta \in \Cp_n$, we say that $\alpha$ refines $\beta$ if one can obtain $\beta$ from $\alpha$ by combining some of its adjacent parts. Alternatively, this means that $\set(\beta) \subseteq \set(\alpha)$. Denote this by $\alpha \preceq \beta$.
For a composition $\alpha$, let $\alpha^\rmr$ denote the composition $(\alpha_{\ell(\alpha)}, \ldots, \alpha_1)$.
The partition obtained by sorting the parts of $\alpha$ in a weakly decreasing order is denoted by $\lambda(\alpha)$.

Now, let us introduce peak compositions, which are central objects in this paper.
A {\em peak composition} of $n$ is a composition of $n$ whose parts (except possibly the last part) are all greater than $1$. 
On the other hand, a subset $P$ of $[n]$ is called a {\em peak set} in $[n]$ if $P \subseteq [2, n-1]$ and for every $i$ in $P$, $i-1$ is not in $P$.
One can see that under the correspondence $\alpha \mapsto \set(\alpha)$, the peak compositions of $n$  precisely correspond to the peak sets of $[n-1]$.
Let $\PCp_n$ denote the set of peak compositions of $n$ and let  
$\PCp:=\bigcup\limits_{n\geq 0}\PCp_n$.
For a finite set $X$ with positive integer entries, define
\[
\Peak(X):=\{i\in X \mid i > 1 \mbox{ and } i-1 \notin X\},
\]
which is referred to as the {\em peak set} of $X$.

\subsection{The ring of quasisymmetric functions and its bases}
\label{The ring of quasisymmetric functions and its bases}
Quasisymmetric functions are formal power series of bounded degree in variables $x_{1},x_{2},x_{3},\ldots$  with coefficients in $\Z$, which are shift invariant in the sense that the coefficient of  $x_{1}^{\alpha_{1}} x_{2}^{\alpha_{2}} \cdots x_{k}^{\alpha_{k}}$ is equal to the coefficient of $x_{i_{1}}^{\alpha_{1}} x_{i_{2}}^{\alpha _{2}} \cdots x_{i_{k}}^{\alpha _{k}}$ for any strictly increasing sequence of positive integers $i_{1}<i_{2}<\cdots <i_{k}$ indexing the variables and any positive integer sequence $(\alpha _{1},\alpha _{2},\ldots,\alpha _{k})$ of exponents.
Let $\QSym$ be the ring of quasisymmetric functions. It 
is a graded ring, decomposing as
\[
\QSym =\bigoplus _{n\geq 0} \QSym_{n},
\]
where $\QSym_n$ is the $\Z$-module consisting of all quasisymmetric functions that are homogeneous of degree $n$. 
For $\alpha\in \Cp_n$, we define the \emph{monomial quasisymmetric function} $M_\alpha$ by 
\[
M_\alpha=\sum_{i_1 < i_2 < \cdots < i_k} x_{i_1}^{\alpha_1} x_{i_2}^{\alpha_2}
\cdots
x_{i_k}^{\alpha_k}.
\]
Clearly $\{M_\alpha \mid \alpha \in \Cp_n\}$ is a basis for $\QSym_n$.
Similarly, we define the \emph{fundamental quasisymmetric function} $F_\alpha$ by 
\[
F_\alpha=\sum_{\substack{1 \le i_1 \le \cdots \le i_n \\ i_j < i_{j+1} \text{ if } j \in \set(\alpha)}} x_{i_1} x_{i_2} \cdots x_{i_n}.
\]
Since 
$F_{\alpha}=\sum_{\beta \preceq \alpha} M_{\beta}$,
it follows that $\{F_{\alpha} \mid \alpha \in \Cp_n \}$ 
is also a basis for $\QSym_n$.

Besides these bases, we introduce the basis of dual immaculate quasisymmetric functions and the basis of Young quasisymmetric Schur functions.
To do this, we recall that for a permutation $\sigma \in \SG_n$, the {\em $($left$)$ descent set} of $\sigma$ is defined by 
\[
\Des(\sigma):=\{i \in [n-1] \mid \text{$i$ is to the right of $i+1$ in $\sigma = \sigma(1) \, \sigma(2) \, \cdots \, \sigma(n)$}\}.
\]
It is important to note that this set is equal to the descent set of $\sigma^{-1}$ for the right weak Bruhat order on $\SG_n$, 
in other words,
\[
\Des(\sigma) = \{ i \in [n-1] \mid \sigma^{-1}(i) > \sigma^{-1}(i+1) \}.
\]

For $\alpha\in \Cp$, the {\em dual immaculate quasisymmetric function} $\DIF{\alpha}$ was initially introduced as 
the quasisymmetric dual of the corresponding immaculate noncommutative symmetric function in~\cite[Section 3.7]{14BBSSZ}. 
A \emph{standard immaculate tableau} (SIT) of \emph{shape} $\alpha$ is a filling of $\tcd(\alpha)$ with distinct entries in $[\,|\alpha|\,]$ such that the entries in each row are increasing from left to right and that in the first column are increasing from bottom to top.
Let $\SIT(\alpha)$ denote the set of all standard immaculate tableaux of shape $\alpha$.
Rather than utilizing the original definition, we adopt the following identity presented in \cite[Proposition 3.48]{14BBSSZ} 
\begin{equation}\label{def of dual immaculate ftn}
\DIF{\alpha}:=\sum_{T \in \SIT(\alpha)}F_{\comp(\Des(\rwRC{r}{T}))}
\end{equation}
as the formal definition for $\DIF{\alpha}$. 
Here, $\rwRC{r}{T}$ represents the word obtained by reading the entries from left to right, starting from the uppermost row of $T$.
It follows from \cite[Proposition 3.37]{14BBSSZ} that $\{\DIF{\alpha} \mid \alpha \in \Cp\}$ is a basis for $\QSym$.

For $\alpha \in \Cp$,  
the {\it Young quasisymmetric Schur function} $\YQS{\alpha}$ is defined by 
$\uprho(\QS{\alpha^\rmr})$,
where $\QS{\alpha^\rmr}$ is the quasisymmetric Schur function associated with $\alpha^\rmr$ and 
$\uprho: \QSym \rightarrow \QSym$ is the automorphism given by $F_\alpha \mapsto F_{\alpha^\rmr}$.
It follows from \cite[Proposition 5.5]{11HLMW} that  $\{\QS{\alpha} \mid \alpha \in \Cp\}$ is a basis for $\QSym$. 
Therefore, these functions form a basis for $\QSym$ and exhibit many analogous properties to quasisymmetric Schur functions.
For detailed information, refer to \cite[Section 4 and Section 5]{13LMvW}.
A \emph{Young composition tableau} $T$ (YCT) of shape $\alpha$ is a filling of $\tcd(\alpha)$ with entries in $\Z_{>0}$ such that
\begin{enumerate}[label = {\rm (\arabic*)}]
\item the entries in each row are weakly increasing from left to right,

\item the entries in the first column are strictly increasing from bottom to top, and 

\item if $i<j$ and $T_{k,j} \leq T_{k+1,i}$, then $(k+1,j) \in \tcd(\alpha)$ and $T_{k+1,j} < T_{k+1,i}$.
\end{enumerate}
The third condition is called the \emph{Young triple condition}.
Pictorially, this condition states that if 
\[
\begin{tikzpicture}
\def \hp {2.7em}
\def \vp {2em}
\node at (\hp*0.5,\vp*0.5) {{\tiny $(k,j)$}};
\node at (\hp*1.5,\vp*0.5) {{\tiny $(k+1,j)$}};
\node at (\hp*1.5,\vp*-0.4) {$\vdots$};
\node at (\hp*1.5,\vp*-1.5) {{\tiny $(k+1,i)$}};

\draw (0,0) rectangle (\hp,\vp);
\draw (\hp,0) rectangle (\hp*2,\vp);
\draw (\hp,\vp*-2) rectangle (\hp*2,\vp*-1);

\node at (\hp*4,\vp*-0.5) {\small with $T_{k,j} \leq T_{k+1,i}$,};
\end{tikzpicture}
\]
then $T_{k+1,j} < T_{k+1,i}$.
A {\em standard Young composition tableau} (SYCT) of shape $\alpha$ is a Young composition tableau in which each of entries in $[\,|\alpha|\,]$ appears exactly once. 
Let $\SYCT(\alpha)$ denote the set of all standard Young composition tableaux of shape $\alpha$.
Then it holds that  
\[
\YQS{\alpha} = \sum_{T\in \SYCT(\alpha)}F_{\comp(\BDes(T))},
\]
where $\BDes(T)$ is the set of $i$'s such that $i$ is weakly right of $i+1$ in $T$
(see \cite[Proposition 5.2.2]{13LMvW}).
In \cite[Lemma 5.3]{22Searles}, Searles showed that $\BDes(T) = \Des(\readSYCT(T))$ for every SYCT $T$,
where $\readSYCT(T)$ denotes the word obtained from $T$ by reading its entries from top to bottom in the first column, and then bottom to top in subsequent columns, reading the columns in order from left to right.
It follows that 
\[
\YQS{\alpha} = \sum_{T\in \SYCT(\alpha)}F_{\comp(\Des(\readSYCT(T)))}.
\]

\subsection{Stembridge's peak ring and its bases}\label{Sec: Stembridge peak ring}

For each peak set $I$ in $[n]$, the peak quasisymmetric function $K_I$ was introduced by Stembridge \cite[Proposition 2.2]{97Stem} as the weight enumerator for enriched $(P,\gamma)$-partitions, where $(P,\gamma)$ is a specific labeled poset. 
Furthermore, he showed in \cite[Proposition 3.5]{97Stem} that 
\[
K_I=2^{|I|+1}\sum_{\substack{\alpha \in \Cp_n \\ I\subseteq \set(\alpha)\triangle(\set(\alpha)+1)}} F_\alpha,
\]
where $\triangle$ represents the symmetric difference of sets, defined as $A \triangle B := (A-B) \cup (B-A)$, 
and $\set(\alpha)+1$ denotes the set $\{i+1 \mid i\in \set(\alpha)\}$.
In this paper, we adopt this identity as the definition of $K_I$. 
Additionally, to maintain consistency with other bases, we parameterize peak quasisymmetric functions using peak compositions rather than peak sets.
To be precise, for each $\alpha\in \PCp$, we use $K_\alpha$ to denote $K_{\set(\alpha)}$.
One can see that $\{K_\alpha \mid \alpha\in \PCp\}$ is linearly independent over $\mathbb Z$.  
For each nonnegative integer $n$, let $\calB_n:=\bigoplus_{\alpha\in \PCp_n}\mathbb Z K_\alpha$.
It is shown in \cite[Proposition 3.5]{97Stem} that $\calB:=\bigoplus_{n\ge 0}\calB_n$
is a graded subring of $\QSym$, called {\it Stembridge's peak ring}. 

In the following, we introduce the basis of quasisymmetric Schur $Q$-functions and the basis of peak Young quasisymmetric Schur functions. These serve as peak analogues to $\{\DIF{\alpha} \mid \alpha \in \Cp \}$ and $\{\YQS{\alpha} \mid \alpha \in \Cp \}$, respectively.

In \cite{15JL}, Jing and Li introduced quasisymmetric Schur $Q$-functions by dualizing noncommutative lifts of Schur $Q$-functions. 
They conjectured that these functions expand positively in the basis $\{K_\alpha \mid \alpha \in \PCp\}$. 
Subsequently, O\u{g}uz confirmed this conjecture in \cite{19Oguz} by providing a detailed description of the expansion. 

\begin{definition}\label{def:SPIT}
Let $\alpha$ be a peak composition of $n$. 
A \emph{standard peak immaculate tableau} (SPIT)\footnote{
In \cite[Definition 4.1]{15JL}, this tableau is referred to as a \emph{standard peak composition tableau}. However, we use this terminology to emphasize that it is a peak analogue of a standard immaculate tableau.} 
of \emph{shape} $\alpha$ is a standard immaculate tableau of shape $\alpha$ such that for every $1\le k \le n$, the subdiagram of $\tcd(\alpha)$ formed by the boxes filled with entries $\le k$ represents the diagram of a peak composition.
\end{definition}

For example, consider the two SITs  
\[
\ctab{S}{  
5 \\  
3 & 4 & 6 \\  
1 & 2 & 7  
}  
\quad \text{and} \quad  
\ctab{T}{  
4 \\  
3 & 5 & 6 \\  
1 & 2 & 7  
} \,.  
\]
It can be observed that \( S \) is an SPIT. In contrast, \( T \) is not an SPIT because the subdiagram containing the entries \( 1, 2, 3, 4 \) is not the diagram of a peak composition.

In this paper, we refer to the last condition in \cref{def:SPIT} as the ``\emph{peak-tableau condition}".
Given a filling $S$ of a composition diagram, we say that $S$ is {\em canonical} if its entries are strictly increasing in order from left to right and from bottom to top. 
With this terminology, standard peak immaculate tableaux can be characterized in the following simple form.

\begin{lemma}\label{lem: peak condition iff canonical}
Let $\alpha$ be a peak composition and $T \in \SIT(\alpha)$.
Then $T$ satisfies the peak-tableau condition if and only if $T^{(\leq 2)}$ is canonical.
\end{lemma}
\begin{proof}
For the `if' part we prove its contraposition.
Suppose that $T$ does not satisfy the peak-tableau condition.
Then there exist entries $c<d$ of $T$ such that 
\begin{itemize}
\item $c$ appears in column $1$, 
\item $d$ is immediately above $c$, and 
\item 
the entries $c+1, c+2, \ldots, d-1$ are located in rows distinct from the one containing $c$.
\end{itemize}
This says that the entry immediately to the right of $c$ is greater than $d$.
Therefore $T^{(\leq 2)}$ is not canonical.

For the `only if' part, we observe that the peak-tableau condition implies the inequalities $T_{1,i+1} > T_{2,i}$ ($1 \leq i \leq \ell(\alpha)-1$).
Since $T_{1,i} < T_{2,i}$ ($1 \leq i \leq \ell(\alpha)$), $T^{(\leq 2)}$ is canonical.
\end{proof}

Let $\SPIT(\alpha)$ denote the set of all standard peak immaculate tableaux of shape $\alpha$.
Rather than relying on the original definition due to Jing and Li, we adopt the expansion obtained by Oğuz as the definition for quasisymmetric Schur $Q$-functions.

\begin{definition}{\rm (\cite[Theorem 3.9]{19Oguz})}\label{def:SchurQformula}
Let $\alpha$ be a peak composition of $n$. The \emph{quasisymmetric Schur $Q$-function} $\QSQ{\alpha}$ is defined by
\[
\QSQ{\alpha}=\sum_{T \in \SPIT(\alpha)} K_{\comp(\Peak(\Des(\rwRC{c}{T})))},
\]
where $\rwRC{c}{T}$ is the word obtained by reading the entries from bottom to top, starting from the leftmost column of $T$.
\end{definition}

It was also stated in \cite[Theorem 3.9]{19Oguz} that $\{\QSQ{\alpha} \mid \alpha \in \PCp\}$ is a basis for $\calB$.

\begin{example}\label{ex:SPIT}
The tableaux in $\SPIT((3,3,1))$ are
\[
\begin{array}{c @{\hskip2\cellsize}c@{\hskip2\cellsize}c@{\hskip2\cellsize}c@{\hskip2\cellsize}c}
 \tableau{ 7 \\ 4 & 5 & {\color{red} 6} \\ 1 & 2 & {\color{red} 3}} & 
\tableau{ 6 \\ 4 & {\color{red} 5} & 7 \\ 1 & 2 & {\color{red} 3}}  &  
\tableau{ 7 \\ 3 & 5 & {\color{red} 6} \\ 1 & {\color{red} 2} & {\color{red} 4}} & 
\tableau{ 6 \\ 3 & {\color{red} 5} & 7 \\ 1 & {\color{red} 2} & {\color{red} 4}} & 
\tableau{ 7 \\ 3 & 4 & {\color{red} 6} \\ 1 & {\color{red} 2} & 5 } \\ \\ 
\tableau{ 6 \\ 3 & 4 & 7 \\ 1 & {\color{red} 2} & {\color{red} 5}}  &  
\tableau{ 7 \\ 3 & 4 & {\color{red} 5} \\ 1 & {\color{red} 2} & {\color{red} 6}} &  
\tableau{ 5 \\ 3 & {\color{red} 4} & 7 \\ 1 & {\color{red} 2} & 6 } & 
\tableau{ 6 \\ 3 & 4 & {\color{red} 5} \\ 1 & {\color{red} 2} & 7 } & 
\tableau{ 5 \\ 3 & {\color{red} 4} & {\color{red} 6} \\ 1 & {\color{red} 2} & 7 } 
\end{array} \ .
\]
Here, the entries in red indicate the descents of $\rwRC{c}{T}$ for $T \in \SPIT((3,3,1))$.
Therefore we have 
\[
\QSQ{(3,3,1)} = K_{(3,3,1)} +  K_{(3,2,2)} +  2K_{(2,2,2,1)} + 2K_{(2,2,3)} + K_{(2,4,1)} +  3K_{(2,3,2)}.   
\]
\end{example}

To introduce peak Young quasisymmetric Schur functions, we need a peak analogue of standard Young composition tableaux.

\begin{definition}{\rm (\cite[Definition 5.6]{22Searles})}\label{def: SPYCT}
Given a peak composition $\alpha$ of $n$, a \emph{standard peak Young composition tableau} $T$ (SPYCT) of shape $\alpha$ is a standard Young composition tableau of shape $\alpha$ satisfying the peak-tableau condition.
\end{definition}

Similar to \cref{lem: peak condition iff canonical}, standard peak Young tableaux can be characterized in the following simple form.

\begin{lemma}{\rm (\cite[Lemma 5.7]{22Searles})}\label{lem: SPYCT canonical iff}
Let $\alpha$ be a peak composition of $n$ 
and $T \in \SYCT(\alpha)$.
Then $T \in \SPYCT(\alpha)$ if and only if the entries in the second column of $T$ increase from bottom to top.    
\end{lemma}

Let $\SPYCT(\alpha)$ denote the set of all standard peak Young composition tableaux of shape $\alpha$.
While peak Young quasisymmetric Schur functions were initially defined as the images of certain $0$-Hecke--Clifford modules under the {\it peak quasisymmetric characteristic}, we choose to adopt the following expansion as the definition. 
For the definition of the peak quasisymmetric characteristic, see \cite[Theorem 5.7]{04BHT}.

\begin{definition}{\rm (\cite[Theorem 5.9]{22Searles})}
Let $\alpha$ be a peak composition of $n$. 
Then the \emph{peak Young quasisymmetric Schur function} $\PYQS{\alpha}$ is given by
\[
\PYQS{\alpha}=\sum_{T \in \SPYCT(\alpha)} K_{\comp(\Peak(\Des(\readSYCT(T))))}.
\]
\end{definition}

It follows from \cite[Theorem 5.11]{22Searles} that $\{\PYQS{\alpha} \mid \alpha \in \PCp\}$ forms a basis for $\calB$.

\begin{example}\label{ex:SPYCT}
Continuing \cref{ex:SPIT}, the elements of $\SPYCT((3,3,1))$ are
\[
\begin{array}{c @{\hskip2\cellsize}c@{\hskip2\cellsize}c@{\hskip2\cellsize}c@{\hskip2\cellsize}c@{\hskip2\cellsize}c@{\hskip2\cellsize}c}
\tableau{ 7 \\ 4 & 5 & {\color{red} 6} \\ 1 & 2 & {\color{red} 3}}  &
\tableau{ 7 \\ 3 & 5 & {\color{red} 6} \\ 1 & {\color{red} 2} & {\color{red} 4}} &
\tableau{ 7 \\ 3 & 4 & {\color{red} 5} \\ 1 & {\color{red} 2} & {\color{red} 6}} &
\tableau{ 6 \\ 4 & {\color{red} 5} & 7 \\ 1 & 2 & {\color{red} 3}} &
\tableau{ 6 \\ 3 & {\color{red} 5} & 7 \\ 1 & {\color{red} 2} & {\color{red} 4}} &
\tableau{ 6 \\ 3 & 4 & {\color{red} 5} \\ 1 & {\color{red} 2} & 7}  &
\tableau{ 5 \\ 3 & {\color{red} 4} & {\color{red} 6} \\ 1 & {\color{red} 2} & 7}
\end{array} \ .
\]
Here the entries in red indicate the descents of $\text{w}_{\text{Y}}(T)$ for $T\in \SPYCT((3,3,1))$.
Therefore we have
\[
\PYQS{(3,3,1)} = K_{(3,3,1)} +  2 K_{(2,2,2,1)} +  2K_{(2,3,2)} + K_{(3,2,2)} + K_{(2,2,3)}.
\]
\end{example}

\section{An equidistribution arsing from \texorpdfstring{$\SPIT(\alpha)$}
{Lg}}\label{A descent preserving bijection}

In this section, we demonstrate that the subsets of $\SG_n$ obtained by applying the column reading $\mathrm{w}_{\rm c}$ and the row reading $\mathrm{w}_{\rm r}$ to $\SPIT(\alpha)$ exhibit the same descent distribution for every $\alpha \in \PCp_n$. 
Using this equidistribution, we develop a combinatorial formula for the quasisymmetric Schur $Q$-functions. 
Moreover, we examine this equidistribution from the viewpoint of the representation theory of $0$-Hecke algebras.

\subsection{A combinatorial bijection on \texorpdfstring{$\SPIT(\alpha)$}{Lg} preserving its descent}
\label{A combinatorial bijection on SPIT preserving its descent}
Let $n \in \Z_{>0}$.
For $\alpha \in \PCp_n$, let 
\begin{align*}
\RWC(\alpha) & := \{\rwRC{c}{T} \mid T \in \SPIT(\alpha)\}\subseteq\SG_n \quad \text{and}\\
\RWR(\alpha) & := \{\rwRC{r}{T} \mid T \in \SPIT(\alpha)\}\subseteq\SG_n.
\end{align*}
We will provide a bijection
\[
\ourmapS: \SPIT(\alpha) \to \SPIT(\alpha)
\]
satisfying that 
$\Des(\rwRC{c}{T}) = \Des(\rwRC{r}{\ourmapS(T)})$ for all $T \in \SPIT(\alpha)$. 
Consequently we deduce that $\RWC(\alpha)$ and $\RWR(\alpha)$ have the same descent distribution, that is, 
\[
\left\{\Des(\sigma) \mid \sigma \in \RWC(\alpha)\right\}
=\left\{\Des(\sigma) \mid \sigma \in \RWR(\alpha)\right\} \
\text{ (as multisets).}
\]
Indeed, our bijection is defined through an appropriate algorithm. Initially, we will introduce the essential definitions and notation needed for the construction of this algorithm.

For $\alpha \in \PCp$, a \emph{generalized standard peak immaculate tableau} (generalized SPIT) $T$ of shape $\alpha$ is a filling of $\tcd(\alpha)$ with distinct entries $a_1 < a_2 < \cdots < a_{|\alpha|}$ in $\Z_{>0}$ such that 
the tableau obtained from $T$ by the replacement $a_i \mapsto i$ is an SPIT of shape $\alpha$.
For a generalized SPIT $T$ and a positive integer $i$,
we define $s_i \cdot T$ as the tableau obtained from $T$ by swapping $i$ and $i+1$ if $i,i+1 \in T$, and the empty tableau otherwise.

For a word $w$ with distinct entries in $\Z_{>0}$, let
\begin{align*}
&\Des(w) := \{i \in \Z_{>0} \mid i,i+1 \in w \text{ and $i$ is to the right of $i+1$ in $w$}\} \quad \text{and}\\ 
&\Asc(w) := \{i \in \Z_{>0} \mid i,i+1 \in w \text{ and $i$ is to the left of $i+1$ in $w$}\}.
\end{align*}

In our algorithm, we will consider pairs $(T,U)$ of generalized SPITs satisfying the following conditions:
\begin{itemize} 
\item (size relation) \ $|U|=|T|-1$.

\item (shape condition) \ The shape of $T$ has at least one part greater than or equal to $3$, and the shape of $U$ equals that of $T_{[1:|U|]}$. 

\item (entry correspondence) \ The set of entries of $U$ equals that of $T_{[1:|U|]}$. 

\item (descent correspondence) \  $\Des(\rwRC{r}{U}) = \Des(\rwRC{c}{T_{[1:|U|]}})$.
\end{itemize}
In this case, we say that $(T,U)$ satisfies the {\it shape-descent alignment condition}.

\begin{algorithm}\label{Subalgorithm on SPIT}
Let $(T,U)$ be a pair of generalized SPITs satisfying 
the shape-descent alignment condition.
Let $h$ be the entry in $T$, but not in $U$.
And, let $(x,y)$ be the position of $h$ in $T$.

\begin{enumerate}[label = {\it Step \arabic*.}, leftmargin=18mm]
\item Set $j:=1$.
Let $U^h_j$ be the tableau obtained from $U$ by adding a box at position $(x, y)$ and filling that box with $h$.

\item 
Let 
\begin{align*}
\setA{j} &:= \Des(\rwRC{r}{U^h_j}) \cap  \Asc(\rwRC{c}{T}) \quad \text{and} \\
\setB{j} &:= \Asc(\rwRC{r}{U^h_j})\cap \Des(\rwRC{c}{T}).
\end{align*}
If $(\setA{j},\setB{j}) = (\emptyset,\emptyset)$, then return $U^h_j$.
Otherwise, go to {\it Step 3}.

\item 
If $\setA{j} \neq \emptyset$, then let $a$ denote the unique element in $\setA{j}$.
Similarly, if $\setB{j} \neq \emptyset$, then let $b$ denote the unique element in $\setB{j}$.
(It will be shown in \cref{Lem1: case j=1} and \cref{Lem2: case j geq 2} that if $\setA{j},\setB{j}$ are nonempty, then they are singletons, respectively). 
Then, we define $z_j$ by
\begin{align*}
\begin{cases}
a & \parbox[t]{12cm}{if $\setA{j},\setB{j} \ne\emptyset$, $b=a+1$, and $a$ is strictly below $a+2$ in $U^h_j$,}\\
a+1 & \parbox[t]{12cm}{if $\setA{j},\setB{j} \ne\emptyset$, $b=a+1$, and $a$ is weakly above $a+2$ in $U^h_j$,}\\
a &  \parbox[t]{12cm}{if $\setA{j}\ne \emptyset,\setB{j} =\emptyset$,}\\
b &  \parbox[t]{12cm}{if $\setA{j}=\emptyset,\setB{j} \ne \emptyset$.}
\end{cases}
\end{align*}
For the $z_j$ defined in this way, set $U^h_{j+1} := s_{z_j} \cdot U^h_j$, increase $j$ by $1$, and then proceed to {\it Step 2}.
\end{enumerate}
\end{algorithm}

\begin{example}
Let 
\[
\ctab{T}{12 \\ 
5 & 6 & 8 \\
3 & 4 & 7 \\
1 & 2 & 9 
}
\quad \text{and} \quad 
\ctab{U}{12 \\
5 & 6  \\ 
3 & 4 & 7 \\
1 & 2 & 9}.
\]
The pair $(T,U)$ satisfies  
the shape-descent alignment condition and $h=8$.
Thus the filling $U^h_1$ in {\it Step 1} is given by {\footnotesize $\ctab{U^h_1}{12 \\ 5 & 6 & 8 \\ 3 & 4 & 7 \\ 1 & 2 & 9}$}.
Note that $\setA{1} = \{7\}$ and $\setB{1} = \{8\}$, hence we go to {\it Step 3}.
Since $7$ is weakly above $9$ in $U^h_1$, $z_1 = 8$, thus  $U^h_2:=s_8 \cdot U^h_1$, that is, {\footnotesize $\ctab{U^h_2}{12 \\ 5 & 6 & 9 \\ 3 & 4 & 7 \\ 1 & 2 & 8}$}. 
Increase $j$ by $1$ and proceed to {\it Step 2}.
Since $\setA{2} = \setB{2} = \emptyset$, our algorithm returns $U^h_2$ as the final value.
\end{example}

Recall that $(T,U)$ is a pair of generalized SPITs satisfying the shape-descent alignment condition.
By the construction of $U^h_1$ in \cref{Subalgorithm on SPIT}, it is evident that the pair $(\setA{1},\setB{1})$ takes one of the following forms:
\[
(\emptyset,\emptyset), \quad (\{h-1\},\emptyset), \quad (\emptyset,\{h\}), \quad (\{h-1\},\{h\}).
\]
The next two lemmas show how $(\setA{j},\setB{j})$ appears when $j\geq 2$.

\begin{lemma}\label{Lem1: case j=1}
Let $(T,U)$ be a pair of generalized SPITs satisfying the shape-descent alignment condition.
If $(\setA{1},\setB{1}) = (\{h-1\},\{h\})$, then either $(\setA{2},\setB{2}) = (\emptyset,\emptyset)$ or 
\[
(\setA{2},\setB{2}) = 
\begin{cases}
(\{h-2\},\emptyset) & \text{if $h-1$ is strictly below $h+1$ in $U^h_1$}, \\ 
(\emptyset,\{h+1\}) & \text{if $h-1$ is weakly above $h+1$ in $U^h_1$}.
\end{cases}
\]
\end{lemma}
\begin{proof}
First, we deal with the case where $h-1$ is strictly below $h+1$ in $U^h_1$. 
Following \cref{Subalgorithm on SPIT}, we see that 
$z_1 = h-1$ and $U^h_2 = s_{h-1} \cdot U^h_1$.
To prove that $(\setA{2},\setB{2})$ is either $(\emptyset,\emptyset)$ or $({h-2},\emptyset)$, it is necessary to examine how the positions of $h$, $h-1$, and $h-2$ change when transitioning from $U^h_1$ to $U^h_2$.

(i) ({\it the change in the position of $h$})
From $\setB{1} = \{h\}$
it follows that $h \in \Asc(\rwRC{r}{U^h_1})$ and $h \in \Des(\rwRC{c}{T})$.
Also, from $U^h_2 = s_{h-1} \cdot U^h_1$, 
it follows that $h$ is strictly below $h+1$ in $U^h_2$.
Therefore, $h \in \Des(\rwRC{r}{U^h_2}) \cap \Des(\rwRC{c}{T})$. 
This implies that $h$ belongs to neither $\setA{2}$ nor $\setB{2}$.

(ii) ({\it the change in the position of $h-1$})
From $\setA{1} = \{h-1\}$
it follows that $h-1 \in \Des(\rwRC{r}{U^h_1})$ and $h-1 \in \Asc(\rwRC{c}{T})$.
Therefore $h-1 \in \Asc(\rwRC{r}{U^h_2})$. This implies that 
$h-1$ belongs to neither $\setA{2}$ nor $\setB{2}$.

(iii) ({\it the change in the position of $h-2$})
If $h-2 \notin T$, then $h-2 \notin U^h_1$. 
Therefore, $(\setA{2},\setB{2}) = (\emptyset,\emptyset)$.
If $h-2 \in T$, then $h-2$ is neither in $\setA{1}$ nor $\setB{1}$.
This implies that the entry $h-2$ appears in either $\Des(\rwRC{r}{U^h_1}) \cap \Des(\rwRC{c}{T})$ or $\Asc(\rwRC{r}{U^h_1}) \cap \Asc(\rwRC{c}{T})$. 

\begin{itemize}
\item If $h-2 \in \Des(\rwRC{r}{U^h_1})$, then $h-2 \in \Des(\rwRC{r}{U^h_2})$.

\item If $h-2 \in \Asc(\rwRC{r}{U^h_1})$ and $h-2$ is weakly above $h$ in $U^h_1$, then $h-2 \in \Asc(\rwRC{r}{U^h_2})$.

\item If $h-2 \in \Asc(\rwRC{r}{U^h_1})$ and $h-2$ is strictly below $h$ in $U^h_1$, then $h-2 \in \Des(\rwRC{r}{U^h_2})$.
\end{itemize}
In the first two cases, $h-2$ does not belong to either $\setA{2}$ or $\setB{2}$, hence $(\setA{2},\setB{2}) = (\emptyset,\emptyset)$.
In the third case, $h-2 \in \setA{2}$, hence $(\setA{2},\setB{2}) = (\{h-2 \},\emptyset)$.

Next, we deal with the case where $h-1$ is weakly above $h+1$ in $U^h_1$. 
Following \cref{Subalgorithm on SPIT}, we see that 
$z_1 = h$ and $U^h_2 = s_{h} \cdot U^h_1$.
Similarly to the previous paragraph, we analyze how the positions of entries $h$, $h-1$, and $h-2$ change during the transition from $U^h_1$ to $U^h_2$.

(i) ({\it the change in the position of $h-1$})
From $\setA{1} = \{h-1\}$,
it follows that $h-1 \in \Des(\rwRC{r}{U^h_1})$ and $h-1 \in \Asc(\rwRC{c}{T})$.
Also, from $U^h_2 = s_{h} \cdot U^h_1$,
it follows that $h-1$ is weakly above $h$ in $U^h_2$.
Therefore, 
$h-1 \in \Asc(\rwRC{r}{U^h_2}) \cap \Asc(\rwRC{c}{T})$. 
This implies that $h-1$ belongs to neither $\setA{2}$ nor $\setB{2}$.

(ii) ({\it the change in the position of $h$})
From $\setB{1} = \{h \}$,
it follows that $h \in \Asc(\rwRC{r}{U^h_1})$ and $h \in \Des(\rwRC{c}{T})$.
Therefore $h \in \Des(\rwRC{r}{U^h_2})$. 
This implies that 
$h$ belongs to neither $\setA{2}$ nor $\setB{2}$.

(iii) ({\it the change in the position of $h+1$})
If $h+2 \notin T$, then $h+2 \notin U^h_1$. 
Therefore, $(\setA{2},\setB{2}) = (\emptyset,\emptyset)$.
If $h+2 \in T$, then $h+1$ is neither in $\setA{1}$ nor $\setB{1}$.
This implies that the entry $h+1$ appears in either $\Des(\rwRC{r}{U^h_1}) \cap \Des(\rwRC{c}{T})$ or $\Asc(\rwRC{r}{U^h_1}) \cap \Asc(\rwRC{c}{T})$. 

\begin{itemize}
\item
If $h+1 \in \Des(\rwRC{r}{U^h_1})$ and $h+2$ is strictly above $h$ in $U^h_1$, then $h+1 \in \Des(\rwRC{r}{U^h_2})$.

\item If $h+1 \in \Des(\rwRC{r}{U^h_1})$ and $h+2$ is weakly below $h$ in $U^h_1$, then $h+1 \in \Asc(\rwRC{r}{U^h_2})$.

\item If $h+1 \in \Asc(\rwRC{r}{U^h_1})$, then $h+1 \in \Asc(\rwRC{r}{U^h_2})$.
\end{itemize}
In the first and third case, $h+1$ does not belong to either $\setA{2}$ or $\setB{2}$, hence $(\setA{2},\setB{2}) = (\emptyset,\emptyset)$.
In the second case, $h+1 \in \setB{2}$, hence $(\setA{2},\setB{2}) = (\emptyset, \{h+1 \})$. 
\end{proof}

\begin{lemma}\label{Lem2: case j geq 2}
Let $(T,U)$ be a pair of generalized SPITs satisfying the shape-descent alignment condition, and $j\geq 1$.
\begin{enumerate}[label = {\rm (\alph*)}]
\item 
If $(\setA{j},\setB{j}) = (\{h-j\},\emptyset)$, then either $(\setA{j+1},\setB{j+1}) = (\emptyset,\emptyset)$ or $(\setA{j+1},\setB{j+1}) = (\{h-j-1\},\emptyset)$.

\item 
If $(\setA{j},\setB{j}) = (\emptyset,\{h+j-1\})$, then either $(\setA{j+1},\setB{j+1}) = (\emptyset,\emptyset)$ or  $(\setA{j+1},\setB{j+1}) = (\emptyset,\{h+j\})$.
\end{enumerate}  
\end{lemma}
\begin{proof}
(a) We use induction on $j$.
For the base case where $j =1$, 
following \cref{Subalgorithm on SPIT}, we see that $z_1 = h-1$ and $U^h_2 = s_{h-1} \cdot U^h_1$.
To prove that $(\setA{2},\setB{2})$ is either $(\emptyset,\emptyset)$ or $(\{h-2\},\emptyset)$,
we analyze how the positions of $h$, $h-1$, and $h-2$ change during the transition from $U^h_1$ to $U^h_2$.

(i) ({\it the change in the position of $h$})
If $h+1 \notin T$, then $h+1 \notin U^h_1$, 
so $h$ does not belong to either $\setA{2}$ or $\setB{2}$.
If $h+1 \in T$, then $h$ is neither in $\setA{1}$ nor $\setB{1}$.
This implies that the entry $h$  appears in either $\Des(\rwRC{r}{U^h_1}) \cap \Des(\rwRC{c}{T})$ or $\Asc(\rwRC{r}{U^h_1}) \cap \Asc(\rwRC{c}{T})$. 
\begin{itemize}
\item 
If $h \in \Des(\rwRC{r}{U^h_1})$, 
then $h \in \Des(\rwRC{r}{U^h_2})$ since $h$ is strictly below $h-1$ in $U^h_1$·

\item 
If $h \in \Asc(\rwRC{r}{U^h_1})$ and $h-1$ is weakly above $h+1$ in $U^h_1$,
then $h \in \Asc(\rwRC{r}{U^h_2})$.

\item 
If $h \in \Asc(\rwRC{r}{U^h_1})$ and $h-1$ is strictly below $h+1$ in $U^h_1$,
then $h+1$ must be strictly to the right of $h$ in $T$ since $h \in \Asc(\rwRC{c}{T})$.  
\end{itemize}
However, the third case contradicts the fact that $h$ is the final letter to appear in $\rwRC{c}{T}$.
Thus, we deduce that if $h+1$ is in $U^h_1$, 
then $h$ does not belong to either $\setA{2}$ or $\setB{2}$.

(ii) ({\it the change in the position of $h-1$})
From $\setA{1} = \{h-1\}$,
it follows that $h-1 \in \Des(\rwRC{r}{U^h_1})$ and $h-1 \in \Asc(\rwRC{c}{T})$.
So $h-1 \in \Asc(\rwRC{r}{U^h_2})$,
and thus we see that $h-1$ does not belong to either $\setA{2}$ or $\setB{2}$.

(iii) ({\it the change in the position of $h-2$})
If $h-2 \notin T$, then $h-2 \notin U^h_1$. 
Therefore, $(\setA{2},\setB{2}) = (\emptyset,\emptyset)$.
If $h-2 \in T$, then $h-2$ is neither in $\setA{1}$ nor $\setB{1}$.
This implies that the entry $h-2$ appears in either $\Des(\rwRC{r}{U^h_1}) \cap \Des(\rwRC{c}{T})$ or $\Asc(\rwRC{r}{U^h_1}) \cap \Asc(\rwRC{c}{T})$.  
\begin{itemize}
\item If $h-2 \in \Des(\rwRC{r}{U^h_1})$, then $h-2 \in \Des(\rwRC{r}{U^h_2})$.

\item If $h-2 \in \Asc(\rwRC{r}{U^h_1})$ and $h-2$ is weakly above $h$ in $U^h_1$, then $h-2 \in \Asc(\rwRC{r}{U^h_2})$.

\item If $h-2 \in \Asc(\rwRC{r}{U^h_1})$ and $h-2$ is strictly below $h$ in $U^h_1$, then $h-2 \in \Des(\rwRC{r}{U^h_2})$.
\end{itemize}
In the first two cases, $h-2$ does not belong to either $\setA{2}$ or $\setB{2}$, hence $(\setA{2},\setB{2}) = (\emptyset,\emptyset)$.
In the third case, $h-2 \in \setA{2}$, hence $(\setA{2},\setB{2}) = (\{h-2 \},\emptyset)$.

\smallskip
Now, assume that the assertion holds for all $j \leq m$.
Suppose that $(\setA{m+1}, \setB{m+1}) = (\{h-m-1\}, \emptyset)$.
Following \cref{Subalgorithm on SPIT}, we see that $z_{m+1} = h-m-1$ and $U^h_{m+2} = s_{h-m-1} \cdot U^h_{m+1}$.
To prove that $(\setA{m+2},\setB{m+2})$ is either $(\emptyset,\emptyset)$ or $(\{h-m-2\},\emptyset)$,
we analyze how the positions of $h-m$, $h-m-1$ and $h-m-2$ change during the transition from $U^h_{m+1}$ to $U^h_{m+2}$. 

(i) ({\it the change in the position of $h-m$})
If $h-m+1 \notin T$, then $h-m+1 \notin U^h_{m+1}$, 
so $h-m$ does not belong to either $\setA{m+2}$ or $\setB{m+2}$.
If $h-m+1 \in T$, then $h-m$ is neither in $\setA{m+1}$ or $\setB{m+1}$.
This implies that the entry $h-m$  appears in either
$\Des(\rwRC{r}{U^h_{m+1}}) \cap \Des(\rwRC{c}{T})$ or $\Asc(\rwRC{r}{U^h_{m+1}}) \cap \Asc(\rwRC{c}{T})$.

\begin{itemize}
\item 
If $h-m \in \Des(\rwRC{r}{U^h_{m+1}})$, 
then $h-m \in \Des(\rwRC{r}{U^h_{m+1}})$.

\item 
If $h-m \in \Asc(\rwRC{r}{U^h_{m+1}})$ and $h-m-1$ is weakly above $h-m+1$ in $U^h_{m+1}$,
then $h-m \in \Asc(\rwRC{r}{U^h_{m+2}})$.

\item 
If $h-m \in \Asc(\rwRC{r}{U^h_{m+1}})$ and $h-m-1$ is strictly below $h-m+1$ in $U^h_{m+1}$,
then $h-m-1 \in \setA{m}$.
This is because that
$(\setA{m}, \setB{m}) = (\{h-m\}, \emptyset)$ and $U^h_{m+1} = s_{h-m} \cdot U^h_{m}$, 
so the placement of $h-m-1$ in both $U^h_{m}$ and $U^h_{m+1}$ remains unchanged.
\end{itemize}
However, the third case contradicts to the induction hypothesis.
Therefore, we deduce that if $h-m+1$ is in $U^h_{m+1}$, then $h-m$ does not belong to either $\setA{m+2}$ or $\setB{m+2}$.

(ii) ({\it the change in the position of $h-m-1$ and $h-m-2$})
Using the same approach as for the letters $h-1$ and $h-2$ in case of $j=1$, it becomes clear that $h-m-1$ does not belong to either $\setA{m+2}$ or $\setB{m+2}$.
Furthermore, $h-m-2 \in \setA{m+2}$ if and only if $h-m-2 \in \Asc(\rwRC{r}{U^h_{m+1}}) \cap \Asc(\rwRC{c}{T})$ and $h-m-2$ is strictly below $h-m$ in $U^h_{m+1}$.
Combining all the observations mentioned,
we conclude that $(\setA{m+2},\setB{m+2})$ is either $(\emptyset,\emptyset)$ or $(\{h-m-2\},\emptyset)$, as desired.

\smallskip
(b) We prove the assertion by using induction on $j$ with a similar approach in (a).
In the case where $j=1$, it can be done by analyzing how the positions of $h+1$, $h$, and $h-1$ change during the transition from $U^h_1$ to $U^h_2$.
Under the assumption that the assertion holds for all $1 \leq j \leq m$,
and by analyzing the positions of $h+j$, $h+j-1$, and $h+j-1$.
\end{proof}

Let $j\ge 2$. 
Due to \cref{Lem1: case j=1} and \cref{Lem2: case j geq 2}, 
if $(\setA{j'},\setB{j'}) \neq (\emptyset,\emptyset)$ for all $1\le j' < j$, then $U^h_j$ is given by
\begin{equation}\label{Eq: decreasing or increasing}
s_{h-j} s_{h-j+1} \cdots s_{h-1} \cdot U^h_1 \quad \text{or} \quad s_{h+j-1} \cdots s_{h+1} s_{h} \cdot U^h_1.
\end{equation}
This implies that there exists a positive integer $l$ such that $\setA{l} = \setB{l} = \emptyset$.
Therefore, \cref{Subalgorithm on SPIT} always terminates after a finite number of steps.

We will provide a clear summary of the changes in entry positions from $U^h_j$ to $U^h_{j+1}$ in the proofs of \cref{Lem1: case j=1} and \cref{Lem2: case j geq 2}. 
This summary will help demonstrate the reversibility of this transition (see \cref{Thm3-2: our map is bij.}).

\begin{summary}\label{Summary: observation h-2 h h-1}
Let $(T,U)$ be a pair satisfying the shape-descent alignment condition.
Let $l$ be the positive integer such that $(\setA{j},\setB{j}) \neq (\emptyset,\emptyset)$ for $1 \leq j < l$, but $(\setA{l},\setB{l}) =  (\emptyset,\emptyset)$.
Suppose that $l \geq 2$, equivalently, $(\setA{1},\setB{1}) \neq (\emptyset,\emptyset)$.
For the sequence $\{z_j\}_{j=1,2,\ldots,l-1}$ in \cref{Subalgorithm on SPIT}, we have the two cases.

{\it Case 1. $z_j=h-j$:} \   
For $1 \leq j < l-1$, the entry $h-j+1$ is strictly below $h-j$ and $h-j-1$ is strictly below $h-j+1$ in $U^h_{j+1}$.
And, $h-l+2$ is strictly below $h-l+1$ in $U^h_{l}$, and one of the following conditions holds:
\begin{itemize}
\item $h-l \notin U^h_l$.
\item $h-l \in U^h_l$ and $h-l$ is strictly below $h-l+2$.
\item $h-l \in U^h_l$ and $h-l$ is weakly above $h-l+1$.
\end{itemize}

{\it Case 2. $z_j=h+j-1$:} \ 
For $1 \leq j < l-1$, the entry $h+j-1$ is strictly below $h+j$ and $h+j+1$ is weakly below $h+j-1$ in $U^h_{j+1}$.
And, $h+l-1$ is strictly below $h+l$ in $U^h_l$, and one of the following conditions holds:
\begin{itemize}
\item $h+l+1 \notin U^h_l$.
\item $h+l+1 \in U^h_l$ and $h+l+1$ is weakly below $h+l-1$.
\item $h+l+1 \in U^h_l$ and $h+l+1$ is strictly above $h+l$.
\end{itemize}
\end{summary}

The next lemma tells us that the final tableau returned by \cref{Subalgorithm on SPIT} is a generalized SPIT.

\begin{lemma}\label{Lem: The resulting tableau via first algorithm is a PIT}
Let $(T,U)$ be a pair of generalized SPITs satisfying the shape-descent alignment condition.
Then \cref{Subalgorithm on SPIT} returns a generalized SPIT.
\end{lemma}
\begin{proof}
The construction of $U^h_1$ in \cref{Subalgorithm on SPIT} ensures that $U^h_1$ is a generalized SPIT. 
Due to the previous discussion, the assertion can be verified by showing that if $U^h_{j}$ is a generalized SPIT and $(\setA{j}, \setB{j}) \neq (\emptyset, \emptyset)$, then $U^h_{j+1}$ is a generalized SPIT.

Suppose that $U^h_{j}$ is a generalized SPIT with $(\setA{j}, \setB{j}) \neq (\emptyset, \emptyset)$.
Let $(x,y)$ be the position of $h$ in $T$.
Since $(T,U)$ is a pair of generalized SPITs satisfying the shape-descent alignment condition, $(x,y)$ cannot be in the first two columns of $T$.
Moreover,  
referring to \cref{Lem1: case j=1} and \cref{Lem2: case j geq 2},
if $z_j$ is in $\setA{j}$, then $z_{j}+1$ is placed at the position $(x,y)$ in $U^h_j$, with $z_j$ positioned strictly below $z_{j}+1$. 
On the other hand, if $z_j \in \setB{j}$, 
then $z_{j}$ is placed at the position $(x,y)$ in $U^h_j$, with $z_{j}+1$ with $z_j$ positioned weakly above $z_{j}+1$.
According to \cref{lem: peak condition iff canonical}, neither $z_j$ nor $z_j+1$ can be in the first two columns of $U^h_j$.
Therefore, $U^h_{j+1} := s_{z_j} \cdot U^h_j$ is a generalized SPIT.
\end{proof}

Let $\alpha \in \PCp_n$ and $\ourc:=c_1(\alpha)+c_2(\alpha)$.
For each $T \in \SPIT(\alpha)$, we recursively define a sequence $\{U_i(T)\}_{i=0,1,\ldots,n-\ourc}$.
Let 
\[
U_0(T) := T_{[1:\ourc]}.
\]
For $1 \leq i \leq n-\ourc$, assume that $U_j(T)$'s have been determined for all $0 \leq j \leq i-1$. 
Then, define $U_{i}(T)$ to be the filling obtained by applying \cref{Subalgorithm on SPIT} to $(T_{[1:\ourc+i]},U_{i-1}(T))$.
We consider the function 
\[
\ourmapS: \SPIT(\alpha) \to \SPIT(\alpha), \quad T \mapsto U_{n-\ourc}(T).
\]

\begin{example}\label{Ex: C to R case}
Let $\alpha=(4,3,4,1)$ and 
\[
\ctab{T}{9 \\ 5 & 6 & 8 & 10 \\ 3 & 4 & 12 \\ 1 & 2 & 7 & 11} \in \SPIT(\alpha). 
\]
For $1 \leq i \leq 5$,
applying \cref{Subalgorithm on SPIT} to each pair $(T_{[1:\ourc+i]},U_{i-1}(T))$, we can observe the following sequence of steps:
\begin{align*}
&\ctab{U_0(T)}{9 \\ 5 & 6 \\ 3 & 4 \\ 1 & 2} 
\stackrel[]{}{\longrightarrow}
\ctab{U_1(T)}{9 \\ 5 & 6 \\ 3 & 4  \\ 1 & 2 & 7} 
\stackrel[]{}{\longrightarrow}
\ctab{U_2(T)}{9 \\ 5 & 6 \\ 3 & 4 & 12 \\ 1 & 2 & 7} \\ 
& \stackrel[]{}{\longrightarrow}
\ctab{U_3(T)}{9 \\ 5 & 6 & 7\\ 3 & 4 & 12 \\ 1 & 2 & 8} 
\quad \left(\ctab{(U_3(T))^{8}_1}{9 \\ 5 & 6 & 8 \\ 3 & 4 & 12 \\ 1 & 2 & 7} \stackrel[]{s_{7}}{\longrightarrow}  
\ctab{(U_3(T))^{8}_2}{9 \\ 5 & 6 & 7\\ 3 & 4 & 12 \\ 1 & 2 & 8} 
  \right) \\
& \stackrel[]{}{\longrightarrow}
\ctab{U_4(T)}{9 \\ 5 & 6 & 7\\ 3 & 4 & 12 \\ 1 & 2 & 8 & 11} 
\stackrel[]{}{\longrightarrow} 
\ctab{U_5(T)}{9 \\ 5 & 6 & 7 & 12\\ 3 & 4 & 11 \\ 1 & 2 & 8 & 10} \\
& \left(\ctab{(U_5(T))^{10}_1}{9 \\ 5 & 6 & 7 & 10 \\ 3 & 4 & 12 \\ 1 & 2 & 8 & 11} 
\stackrel[]{s_{10}}{\longrightarrow}  \ctab{(U_5(T))^{10}_2}{9 \\ 5 & 6 & 7 & 11 \\ 3 & 4 & 12 \\ 1 & 2 & 8 & 10} 
\stackrel[]{s_{11}}{\longrightarrow}  \ctab{(U_5(T))^{10}_3}{9 \\ 5 & 6 & 7 & 12 \\ 3 & 4 & 11 \\ 1 & 2 & 8 & 10} 
\right).
\end{align*}
Then $\ourmapS(T) = U_5(T)$.
Note that $\Des(\rwRC{c}{T}) = \{2,4,8,10,11\} = \Des(\rwRC{r}{\ourmapS(T)})$.
\end{example}

\begin{theorem}\label{Thm3-1: bijection1 on SPITs}
For every peak composition $\alpha$, 
the map $\ourmapS$ is well-defined, and for all $T\in \SPIT(\alpha)$, it holds that $\Des(\rwRC{c}{T})=\Des(\rwRC{r}{\ourmapS(T)})$.
\end{theorem}
\begin{proof}
Consider the first pair $(T_{[1:\ourc+1]},U_0(T))$.
Since $U_0(T) = T_{[1:\ourc]}$, it satisfies the shape-descent alignment condition.
By \cref{Lem: The resulting tableau via first algorithm is a PIT}, the resulting filling $U_1(T)$ is a generalized SPIT. 
And, by the construction of \cref{Subalgorithm on SPIT}, $\Des(\rwRC{r}{U_1(T)}) = \Des(T_{[1:\ourc+1]})$.
In the same manner as above, the resulting tableaux $U_i(T)$'s ($2 \leq i \leq n-\ourc$) are generalized SPITs and $\Des(\rwRC{r}{U_i(T)}) = \Des(T_{[1:\ourc+i]})$ for $1 \leq i \leq n-\ourc$.
Hence, we conclude that $\ourmapS(T) \in \SPIT(\alpha)$ and $\Des(\rwRC{c}{T})=\Des(\rwRC{r}{\ourmapS(T)})$.
\end{proof}

Consider the situation where each element in $\{\Des(\sigma) \mid \sigma \in \RWC(\alpha)\}$ occurs with multiplicity 1.
In this instance, \cref{Thm3-1: bijection1 on SPITs}
implies that $\ourmapS$ is a bijection.
In what follows, we will show that this function is a bijection in general.
This can be accomplished by showing the reversibility of each step in \cref{Subalgorithm on SPIT}.

\begin{algorithm}\label{subalgorithm: down}
Let $V$ be a generalized SPIT and $g \in V$.

\begin{enumerate}[label = {\it Step \arabic*.}, leftmargin=18mm]
\item Set $j:=1$.
Let $a_j := g-1$ and $V_{j} := V$ and $\bfd_{(j)}$ be the sequence $\{a_j\}$.
\item 
Increase $j$ by $1$ and then 
let $V_{j} := s_{a_{j-1}} \cdot V_{j-1}$.
Examine the conditions listed below:
\begin{itemize}
\item $g-j \in V_{j}$, 
     
\item $g-j$ is weakly above $g-j+2$ in $V_{j}$, and

\item $g-j$ is strictly below $g-j+1$ in $V_{j}$.
\end{itemize}
If any of the conditions fail, then return $\bfd_{(j-1)}$.
Otherwise, go to {\it Step 3}.
\item 
Let $a_{j} := g-j$ and define $\bfd_{(j)}$ to be the sequence obtained from $\bfd_{(j-1)}$ by appending $a_{j}$ at the end.
Then, proceed to {\it Step 2}.
\end{enumerate}
\end{algorithm}

\begin{algorithm}\label{subalgorithm: up}
Let $V$ be a generalized SPIT and $g \in V$.
\begin{enumerate}[label = {\it Step \arabic*.}, leftmargin=18mm]
\item Set $j:=1$.
Let $b_j = g$ and $V_{(j)} = V$ and $\bfu_{(j)}$ be the sequence $\{b_j\}$.

\item 
Increase $j$ by $1$ and then 
let $V_{j} = s_{b_{j-1}} \cdot V_{j-1}$.
Examine the conditions listed below:
\begin{itemize}
\item $g+j \in V_{j}$, 

\item $g+j$ is strictly below $g+j-1$ in $V_{j}$, and

\item $g+j$ is weakly below $g+j-2$ in $V_{j}$.
\end{itemize}
If any of the conditions fail, then return $\bfu_{(j-1)}$.
Otherwise, go to {\it Step 3}.

\item 
Let $b_{j} = g+j-1$ and define $\bfu_{(j)}$ to be the sequence obtained from $\bfu_{(j-1)}$ by appending $b_{j}$ at the end.
Then, proceed to {\it Step 2}.
\end{enumerate}    
\end{algorithm}

Let $(V,g)$ be a pair of a generalized SPIT and a positive integer.
Since $|V|$ is finite, these processes terminate after a finite number of steps.
We denote by $\Downseq{V}{g}$ and $\Upseq{V}{g}$ the final sequences obtained from \cref{subalgorithm: down} and \cref{subalgorithm: up}, respectively,
that is,
\begin{equation}\label{Eq: Dseq is decreasing and Useq is increasing}
\begin{aligned}
&\Downseq{V}{g} = \{a_j\}_{j = 1,2,\ldots, l}, \quad \text{where } a_j = g-j, \\
&\Upseq{V}{g} = \{b_j\}_{j = 1,2,\ldots, l'}, \quad \text{where } b_j = g+j-1.
\end{aligned}
\end{equation}
Here, $l$ and $l'$ are the lengths of $\Downseq{V}{g}$ and $\Upseq{V}{g}$, respectively.
Let 
\[
\Downseq{V}{g} \cdot V:= s_{a_l} s_{a_{l-1}} \cdots  s_{a_1} \cdot V
\quad \text{and} \quad 
\Upseq{V}{g} \cdot V:= s_{b_{l'}} s_{b_{l'-1}} \cdots  s_{b_1} \cdot V.
\]
We say that {\em $(V,g)$ satisfies the condition $(\mathbf{C1})$} if
\begin{enumerate}
\item both $g-1$ and $g$ are entries of $V$, with $g-1$ strictly below $g$, and

\item whenever $g+1 \in V$, $g+1$ is either weakly below $g-1$ or strictly above $g$.
\end{enumerate}
Similarly, we say that {\em $(V,g)$ satisfies the condition $(\mathbf{C2})$} if 
\begin{enumerate}
\item both $g$ and $g+1$ are entries of $V$, 
with $g+1$ strictly below $g$, and

\item whenever $g-1 \in V$, $g-1$ is either strictly below $g+1$ or weakly above $g$.
\end{enumerate}
It is evident that conditions $(\mathbf{C1})$ and $(\mathbf{C2})$ are mutually exclusive.

Given $T \in \SPIT(\alpha)$, we introduce a sequence $\{V^{i}(T)\}_{i=0,1,\ldots,n-\ourc}$ consisting of certain fillings.
Let
\[
V^0(T):=T.
\]
Given $1 \leq i \leq n-\ourc$, assume that $V^{j}(T)$'s have been defined for all $0 \leq j \leq i-1$. 
Let $g$ be the entry of the uppermost box in the rightmost column of $V^{i-1}(T)$.
We now define
\begin{equation}\label{Eq: ViT}
V^{i}(T):=
\begin{cases}
\left(\Downseq{V^{i-1}(T)}{g} \cdot V^{i-1}(T)\right)_{[1:n-i]} & \text{if $(V^{i-1}(T),g)$ satisfies $(\mathbf{C1})$,}  \\
\left(\Upseq{V^{i-1}(T)}{g} \cdot V^{i-1}(T)\right)_{[1:n-i]} & \text{if $(V^{i-1}(T),g)$ satisfies $(\mathbf{C2})$,} \\
V^{i-1}(T)_{[1:n-i]} & \text{otherwise.}
\end{cases}
\end{equation}

\begin{lemma}\label{Lem: ViT is a generalized SPIT}
Let $\alpha \in \PCp_n$ and $T \in \SPIT(\alpha)$.
For every $0 \leq i \leq n - \ourc$, $V^i(T)$ is a generalized SPIT.
\end{lemma}
\begin{proof}
We prove this assertion by induction on $i$.
The case where $i=0$, it is obvious since $V^0(T) = T$. 
Given $1 \leq i \leq n-\ourc$, assume that $V^{i-1}(T)$ is a generalized SPIT.
Let $(x,y)$ be the position of the uppermost box in the rightmost column of $V^{i-1}(T)$, and let $g:=V^{i-1}(T)_{x,y}$.
If $(V^{i-1}(T),g)$ satisfies $(\mathbf{C1})$, then $\Downseq{V^{i-1}(T)}{g} = \{g-j\}_{j=1,2,\ldots,l}$ for some $l \geq 1$.
For $1 \leq j < l$, one observes that $g-j-1$ is strictly below $g-j$ in $s_{g-j} \cdots s_{g-1} \cdot V^{i-1}(T)$.
Combining this with the induction hypothesis, we derive that $V^{i}(T) = \Downseq{V^{i-1}(T)}{g} \cdot V^{i-1}(T)$ is a generalized SPIT.
If $(V^{i-1}(T),g)$ satisfies $(\mathbf{C2})$, then $\Upseq{V^{i-1}(T)}{g} = \{g+j-1\}_{j=1,2,\ldots,l'}$ for some $l' \geq 1$.
For $1 \leq j < l$, one observes that $g+j$ is strictly below $g+j-1$ in $s_{g+j-1} \cdots s_{g} \cdot V^{i}(T)$.
Combining this with the induction hypothesis, we derive that $V^{i}(T) = \Upseq{V^{i-1}(T)}{g} \cdot V^{i-1}(T)$ is a generalized SPIT.
If neither $(\mathbf{C1})$ nor $(\mathbf{C2})$ holds, then $V^{i}(T) = (V^{i-1}(T))_{[1:|V^{i-1}(T)|-1]}$. 
Thus, by induction hypothesis, $V^{i}(T)$ is a generalized SPIT.
\end{proof}

Let $T \in \SPIT(\alpha)$.
Note that the shape of $V^{i}(T)$ is obtained from the shape of $V^{i-1}(T)$ by removing the uppermost box in the rightmost column.
For each $1 \leq i \leq n-\ourc$, let $y_i(T)$ be the unique entry in $V^{i-1}(T)$, but not in $V^{i}(T)$.
We now assume that the sequence $\{(V^i(T),y_i(T))\}_{i=1,\ldots,n-\ourc}$ has been obtained by the previously described process. 
Then we define the map
\[
\ourmapSS: \SPIT(\alpha) \rightarrow \SPIT(\alpha), \quad T \mapsto \ourmapSS(T),
\]
where $\ourmapSS(T)$ is defined to be the filling of $\tcd({\alpha})$ such that 
\begin{equation*}\label{Def: V_bullet(T)}
\rwRC{c}{\ourmapSS(T)}_i = 
\begin{cases}
\rwRC{c}{V^{n-\ourc}(T)}_i    & \text{if } 1 \leq i \leq \ourc,\\
y_{n+1-i}(T) & \text{if } \ourc+1 \leq i \leq n.
\end{cases}
\end{equation*}

For the reader's understanding, we give an example before proving that $\ourmapSS$ is well-defined.
\begin{example}\label{Ex: R to C case}
Let us consider the SPIT
\[
\ctab{T}{9 \\
5 & 6 & 7 & 12 \\ 
3 & 4 & 11 \\ 
1 & 2 & 8 & 10} \in \SPIT((4,3,4,1)).
\]
We explain how to obtain $V^{1}(T)$ from $V^{0}(T):=T$.
Start with $g:=12$.
Since $g-1 \in V^{0}(T)$ and $g-1$ is strictly below $g$, but $g+1 \notin V^{0}(T)$, $(V^0(T),g)$ satisfies $(\mathbf{C1})$.
Since $\Downseq{V^{0}(T)}{g} = \{a_j\}_{j=1,2}$, where $a_j=g-j$,
from~\eqref{Eq: ViT} we have 
\[
V^{1}(T):=\ctabT{9 \\ 5 & 6 & 7 \\ 3 & 4 & 12 \\ 1 & 2 & 8 & 11}   \quad  \left( 
\begin{array}{c}
\tableau{9 \\ 5 & 6 & 7 & 12 \\ 3 & 4 & 11 \\ 1 & 2 & 8 & 10}
\end{array}
\stackrel[]{s_{11}}{\longrightarrow}
\begin{array}{c}
\tableau{9 \\ 5 & 6 & 7 & 11 \\ 3 & 4 & 12 \\ 1 & 2 & 8 & 10}
\end{array}
\stackrel[]{s_{10}}{\longrightarrow}
\begin{array}{c}
\tableau{9 \\ 5 & 6 & 7 & 10 \\ 3 & 4 & 12 \\ 1 & 2 & 8 & 11}
\end{array} \right)
\]
For $1 \leq i \leq 4$, applying this process iteratively to the filling $V^{i}(T)$ yields the filling $V^{i+1}(T)$ as follows:
\begin{align*}
&V^{2}(T):=\ctabT{9 \\ 5 & 6 & 7 \\ 3 & 4 & 12 \\ 1 & 2 & 8 } 
\qquad 
V^{3}(T):=\ctabT{9 \\ 5 & 6 \\ 3 & 4 & 12 \\ 1 & 2 & 7 }  \quad  
\left( 
\begin{array}{c}
\tableau{9 \\ 5 & 6 & 7 \\ 3 & 4 & 12 \\ 1 & 2 & 8}
\end{array}
\stackrel[]{s_{7}}{\longrightarrow}
\begin{array}{c}
\tableau{9 \\ 5 & 6 & 8 \\ 3 & 4 & 12 \\ 1 & 2 & 7}
\end{array}\right)
\\
& V^{4}(T):=\ctabT{9 \\ 5 & 6 \\ 3 & 4 \\ 1 & 2 & 7 } \qquad  V^{5}(T):=\ctabT{9 \\ 5 & 6 \\ 3 & 4  \\ 1 & 2 } 
\end{align*}
From $V^i(T)$'s ($1 \leq i \leq 5$), we get 
\[
y_1=10, \ y_2=11, \ y_3=8, \ y_4=12, \ y_5=7.
\]
Finally, we have that
\[
\ctab{\ourmapSS(T)}{9 \\ 5 & 6 & 8 & 10 \\ 3 & 4 & 12 \\ 1 & 2 & 7 & 11}.
\]
Note that $\Des(\rwRC{r}{T}) = \{2,4,8,10,11\} = \Des(\rwRC{c}{\ourmapSS(T)})$.

In $\SPIT(\alpha)$, there exist four additional tableaux $T_i$ ($i=1,2,3,4$) with $\Des(\rwRC{c}{T_i})=\{2,4,8,10,11\}$.
The images of these tableaux under $\ourmapS$ are as follows:
\begin{align*}
& \ctab{T_1}{12 \\ 5 & 6 & 9 & 10 \\ 3 & 4 & 11 \\ 1 & 2 & 7 & 8} \stackrel{\ourmapS}{\longmapsto} \ctabT{12 \\ 5 & 6 & 9 & 11 \\ 3 & 4 & 10 \\ 1 & 2 & 7 & 8} 
& \ctab{T_2}{12 \\ 5 & 6 & 8 & 10 \\ 3 & 4 & 9 \\ 1 & 2 & 7 & 11} \stackrel{\ourmapS}{\longmapsto} \ctabT{12 \\ 5 & 6 & 7 & 11 \\ 3 & 4 & 9 \\ 1 & 2 & 8 & 10} \\
& \ctab{T_3}{12 \\ 5 & 6 & 8 & 10 \\ 3 & 4 & 7 \\ 1 & 2 & 9 & 11} \stackrel{\ourmapS}{\longmapsto} \ctabT{12 \\ 5 & 6 & 9 & 11 \\ 3 & 4 & 7 \\ 1 & 2 & 8 & 10} \ 
& \ctab{T_4}{12 \\ 5 & 6 & 7 & 8 \\ 3 & 4 & 11 \\ 1 & 2 & 9 & 10}   \stackrel{\ourmapS}{\longmapsto} \ctabT{12 \\ 5 & 6 & 7 & 9 \\ 3 & 4 & 11 \\ 1 & 2 & 8 & 10} 
\end{align*}
\end{example}

\begin{lemma}\label{lem: well-defined}
For every peak composition $\alpha$, 
the map $\ourmapSS$ is well-defined.
\end{lemma}
\begin{proof}
Let $T \in \SPIT(\alpha)$.
To verify that $\ourmapSS(T) \in \SPIT(\alpha)$, it suffices to show that
\begin{enumerate}[label = {\rm (\roman*)}]
\item if $y_p(T)$ is placed to the right of $y_q(T)$ in $\ourmapSS(T)$, then $y_p(T) > y_q(T)$, and

\item $\ourmapSS(T)^{(\leq 2)}$ is canonical.
\end{enumerate}

For (i), let $(c_p,d_p)$ $(1 \leq p \leq n)$ denote the position of the $p$th box in $\tcd(\alpha)$ from top to bottom and right to left.
According to the assumption, $p < q$, $d_p = d_q$ and $c_p = c_q+1$.
Now, suppose for contradiction that $y_p < y_q$. 
By \cref{Lem: ViT is a generalized SPIT}, $V^{p}(T)$ forms a generalized SPIT, 
implying that the entries in row $d_p$ increase from left to right. 
Since $V^{p}(T)_{c_q,d_q} < y_p$, there should exist $p < k \leq q$ such that 
\[
V^{k-1}(T)_{c_q,d_q} < y_p, 
\quad \text{but} \quad 
V^{k}(T)_{c_q,d_q} > y_p.
\]
However, such a situation does not arise because either the sequence $\Downseq{V^{k-1}(T)}{g}$ is a consecutively decreasing one or $\Upseq{V^{k-1}(T)}{g}$ is a consecutively increasing one, as shown in \eqref{Eq: Dseq is decreasing and Useq is increasing}.
Therefore, we deduce that $y_q < y_p$.

For (ii), due to the construction of $\ourmapSS(T)$, it is clear that $(\ourmapSS(T))^{(\leq 2)} = V^{n-\ourc}(T) = T^{(\leq 2)}$.
Since $T^{(\leq 2)}$ is canonical, so is $\ourmapSS(T)^{(\leq 2)}$.
\end{proof}

\begin{theorem}\label{Thm3-2: our map is bij.}
For every peak composition $\alpha$, the map $\ourmapS$ is a bijection.
\end{theorem}
\begin{proof}
For $T \in \SPIT(\alpha)$, we recall that 
$U_0(T) = T^{(\leq 2)}$ and $U_{i}(T)$ denotes the generalized SPIT obtained by applying \cref{Subalgorithm on SPIT} to $(T_{[1:\ourc+i]}, U_{i-1}(T))$ for $1 \leq i \leq n-\ourc$. Moreover, $\ourmapS(T) = U_{n-\ourc}(T)$. 
To confirm the assertion, it suffices to show that $\ourmapSS \circ \ourmapS = {\rm id}$. 
This equality can be established by proving the identities
\begin{align}\label{eq: reversibility}
U_{n-\ourc-i}(T) = V^{i}(\ourmapS(T)) \text{ for }i=0,1,\ldots,n-\ourc .
\end{align}
We prove \eqref{eq: reversibility} by induction on $i$.
First, consider $i=0$.
By definition, we have that
\[
U_{n-\ourc}(T) = \ourmapS(T)=V^{0}(\ourmapS(T)).
\]

Next, let $1 \leq i \leq n-\ourc$.
Assuming that $U_{n-\ourc-p}(T) = V^{p}(\ourmapS(T))$ for all $p \leq i-1$, we aim to show $U_{n-\ourc-i}(T) = V^{i}(\ourmapS(T))$.
To do this, let us consider the sequence $\{(U_{n-\ourc-i}(T))^h_j\}_{j=1,2,\ldots,l}$ obtained from the pair $(T_{[1:n-\ourc-i+1]},U_{n-\ourc-i}(T))$
by applying \cref{Subalgorithm on SPIT}.
Here, $h$ is the positive integer associated with this pair, equals $\rwRC{c}{T}_{n-\ourc-i+1}$, and 
$l$ is the positive integer such that $(\setA{j},\setB{j}) \neq (\emptyset,\emptyset)$ for $1 \leq j < l$, but $(\setA{l},\setB{l}) =  (\emptyset,\emptyset)$.
We now consider two cases: $l=1$ or $l > 1$.

First, we deal with the case where $l=1$.
By the construction of \cref{Subalgorithm on SPIT}, it holds that
\begin{equation}\label{Eq: Unci1 = Uncih1}
U_{n-\ourc-i+1}(T) = (U_{n-\ourc-i}(T))^h_1,
\end{equation}
that is, $U_{n-\ourc-i+1}(T)$ is obtained from $U_{n-\ourc-i}(T)$ by adding a box at position $(x,y)$ and filling that box with $h$.
Here, $(x,y)$ is the position of $h$ in $T_{[1:n-\ourc-i+1]}$.
By the induction hypothesis, $V^{i-1}(\ourmapS(T))=U_{n-\ourc-i+1}(T)$, so it follows from \eqref{Eq: Unci1 = Uncih1} that 
\[ 
V^{i-1}(\ourmapS(T)) = (U_{n-\ourc-i}(T))^h_1.
\]
Since the entry $g$ associated with $V^{i-1}(\ourmapS(T))$ equals $h$ and $\setA{1}=\setB{1}=\emptyset$,
we can observe that 
\begin{itemize}
\item if $h-1 \in (U_{n-\ourc-i}(T))^h_1$, then $h-1$ is weakly above $h$, and 

\item if $h+1 \in (U_{n-\ourc-i}(T))^h_1$, then $h+1$ is also weakly above $h$.
\end{itemize}
For the definition of $g$, see \eqref{Eq: ViT}.
This observation says that $(V^{i-1}(\ourmapS(T)),h)$ satisfies neither the condition $(\mathbf{C1})$ nor the condition $(\mathbf{C2})$.
By \eqref{Eq: ViT}, we have  
\[
V^{i}(\ourmapS(T)) = V^{i-1}(\ourmapS(T))_{[1:n-\ourc-i]} = U_{n-\ourc-i}(T).
\]

Second, we deal with the case where $l > 1$.
Consider the sequence $\{z_j\}_{j=1,\ldots,l-1}$ for the pair 
$(T_{[1:n-\ourc-i+1]}, U_{n-\ourc-i}(T))$.
By \cref{Lem1: case j=1} and \cref{Lem2: case j geq 2}, we see that either $z_j=h-j$ for all $j$'s or $z_j=h+j-1$ for all $j$'s.
In the former case, we have 
\begin{align*}
&V^{i-1}(\ourmapS(T)) \\
&= U_{n-\ourc-i+1}(T)  \hspace*{50mm} \text{(by induction hypothesis)}\\
&= s_{h-l+1} s_{h-l+2} \cdots s_{h-1} \cdot (U_{n-\ourc-i}(T))^h_1  \hspace*{10mm} \text{(by the construction of $U_{n-\ourc-i+1}(T) $)}.
\end{align*}
Following \cref{Subalgorithm on SPIT}, we will write $s_{h-l+1} s_{h-l+2} \cdots s_{h-1} \cdot (U_{n-\ourc-i}(T))^h_1$ as $U^h_l$.
Then, we consider the pair $(V^{i-1}(\ourmapS(T),h-l+1)$.
According to {\it Case 1} in \cref{Summary: observation h-2 h h-1}, we see that 
$h-l+2 \in U^h_l$, $h-l+2$ is strictly below $h-l+1$ in $U^h_l$, and one of the following conditions holds:
\begin{itemize}
\item $h-l \notin U^h_l$.
\item $h-l \in U^h_l$ and $h-l$ is strictly below $h-l+2$.
\item $h-l \in U^h_l$ and $h-l$ is weakly above $h-l+1$.
\end{itemize}
By these properties, the pair $(V^{i-1}(\ourmapS(T)),h-l+1)$ satisfies the condition $(\mathbf{C2})$, and therefore 
\begin{align*}
V^{i}(\ourmapS(T)) = \left(\Upseq{V^{i-1}(\ourmapS(T))}{h-l+1}\cdot V^{i-1}(\ourmapS(T))\right)_{[1:n-\ourc-i]}
\end{align*}
(see \eqref{Eq: ViT}).
On the other hand, the sequence $\Upseq{V^{i-1}(\ourmapS(T))}{h-l+1}$ is given by $\{b_j\}_{j=1,2,\ldots,l-1}$ with $b_j = h-l+j$.
This tells us that $z_j = b_{l-j}$ for all $j=1,2,\ldots,l-1$, and therefore
\begin{align*}
V^{i}(\ourmapS(T)) = ((U_{n-\ourc-i}(T))^h_1)_{[1:n-\ourc-i]} = U_{n-\ourc-i}(T).
\end{align*}
The latter case can be dealt with in the same manner as above.
\end{proof}

\begin{corollary}\label{Coro: new definition of QSQ}
Let $\alpha$ be a peak composition.
\begin{enumerate}[label = {\rm (\alph*)}]
\item As multisets, we have 
\[
\left\{\Des(\sigma) \mid \sigma \in \RWC(\alpha)\right\}
=\left\{\Des(\sigma) \mid \sigma \in \RWR(\alpha)\right\}.
\]

\item The quasisymmetric Schur $Q$-function is equal to the $K$-generating function of standard peak immaculate tableaux for the reading $\mathrm{w}_{\rm r}$, where $K$ is used to denote peak quasisymmetric functions.
More precisely,
\[
\QSQ{\alpha} = \sum_{T\in \SPIT(\alpha)}K_{\comp(\Peak(\Des(\rwRC{r}{T}))).}
\]
\end{enumerate}
\end{corollary}
\begin{proof}
(a) The assertion can be obtained by combining \cref{Thm3-1: bijection1 on SPITs} with \cref{Thm3-2: our map is bij.}.

(b) The assertion can be obtained by combining \cref{def:SchurQformula} with (a).
\end{proof}

\begin{remark}\label{Rmk: Xi_alpha (T)=T for alpha}
When $\alpha$ is a peak composition with at most one part greater than $2$, $\ourmapS$ is the identity map. 
In the case where $\alpha$ has no parts greater than $2$, this follows from $|\SPIT(\alpha)| = 1$. On the other hand, if $\alpha$ has only one part greater than $2$, \cref{lem: peak condition iff canonical} implies that for $T \in \SPIT(\alpha)$, $i \in \Des(\rwRC{r}{T})$ if and only if $i+1$ appears strictly upper-left of $i$ in $T$. Consequently, $U_i(T) = T_{[1:\ourc+i]}$ for all $1 \leq i \leq n-\ourc$, implying $\ourmapS(T) = T$.
\end{remark}

\subsection{A representation-theoretical interpretation of \texorpdfstring{\cref{Thm3-2: our map is bij.}}{Lg}}
\label{A representation-theoretical interpretation}

For each $\alpha \in \PCp_n$, we present two $0$-Hecke modules $\HnZmodC_{\alpha}$ and $\HnZmodR_{\alpha}$ such that  
\[
\ch(\HnZmodC_{\alpha}) = \sum_{\sigma \in \RWC(\alpha)} F_{\comp(\Des(\sigma))} \quad \text{and} \quad
\ch(\HnZmodR_{\alpha}) = \sum_{\sigma \in \RWR(\alpha)} F_{\comp(\Des(\sigma))},
\]
where $\ch$ is the quasisymmetric characteristic (see \eqref{quasi characteristic}).
The equidistribution demonstrated in \cref{A combinatorial bijection on SPIT preserving its descent} implies that these two modules share the same image under the quasisymmetric characteristic. 
However, it turns out that they are not isomorphic as $H_n(0)$-modules in general.

To begin with, we briefly review the representation theory of the $0$-Hecke algebra $H_n(0)$. 
For more details, see \cite{21CKNO,22CKNO2}.
The {\em $0$-Hecke algebra} $H_n(0)$ is an associative $\C$-algebra with $1$ generated by the elements $\opi_1,\opi_2,\ldots,\opi_{n-1}$ subject to the following relations:
\begin{align*}
\opi_i^2 &= -\opi_i \quad \text{for $1\le i \le n-1$},\\
\opi_i \opi_{i+1} \opi_i &= \opi_{i+1} \opi_i \opi_{i+1}  \quad \text{for $1\le i \le n-2$},\\
\opi_i \opi_j &= \opi_j \opi_i \quad \text{if $|i-j| \ge 2$}.
\end{align*}
Let $\pi_i:= \opi_i+1$. 
Then, the set $\{\pi_i \mid  i \in [n-1]\}$ is also a generating set of $H_n(0)$.
For any reduced expression $s_{i_1} s_{i_2} \cdots s_{i_p}$ for $\sigma \in \SG_n$, let $\opi_{\sigma} := \opi_{i_1} \opi_{i_2} \cdots \opi_{i_p}$ and $\pi_{\sigma} := \pi_{i_1} \pi_{i_2 } \cdots \pi_{i_p}$.
It is well known that these elements are independent of the choice of reduced expressions, and both $\{\opi_\sigma \mid \sigma \in \SG_n\}$ and $\{\pi_\sigma \mid \sigma \in \SG_n\}$ are $\mathbb C$-bases for $H_n(0)$.

According to \cite{79Norton}, there are $2^{n-1}$ pairwise inequivalent irreducible (left) $H_n(0)$-modules, which are naturally indexed by compositions of $n$.
For each $\alpha \in \Cp_n$, the corresponding irreducible $H_n(0)$-module $\mathbf{F}_{\alpha}$ is defined to be the $1$-dimensional $\mathbb{C}$-vector space spanned by a nonzero vector $v_{\alpha}$ with the $H_n(0)$-action 
\[
\opi_i \cdot v_\alpha = \begin{cases}
0 & i \notin \set(\alpha),\\
- v_\alpha & i \in \set(\alpha)
\end{cases}
\]
for each $1 \le i \le n-1$. 
Let $\modL$ be the category of $H_n(0)$-modules
and $\calR(\modL)$ be the $\Z$-span of the isomorphism classes of finite dimensional $H_n(0)$-modules. 
We denote by $[M]$ the isomorphism class corresponding to an $H_n(0)$-module $M$. 
The \emph{Grothendieck group} $G_0(\modL)$ is the quotient of $\calR(\modL)$ modulo the relations $[M] = [M'] + [M'']$ whenever there exists a short exact sequence $0 \rightarrow M' \rightarrow M \rightarrow M'' \rightarrow 0$. 
By abusing notation, we denote by $[M]$ the element of $G_0(\modL)$ corresponding to an $H_n(0)$-module $M$.
Then $G_0(\modL)$ is a free abelian group with a basis $\{[\bfF_{\alpha}] \mid \alpha \in \Cp_n\}$. Let
\[
\calG(\modLb) := \bigoplus_{n \ge 0} G_0(\modL).
\]
It was shown in \cite{96DKLT, 97KT} that the linear map
\begin{equation}\label{quasi characteristic}
\ch: \calG(\modLb) \rightarrow \Qsym, \quad [\bfF_{\alpha}] \mapsto F_{\alpha},
\end{equation}
called the \emph{quasisymmetric characteristic}, is a ring isomorphism.

In the remainder of this subsection, let us fix $\alpha \in \PCp_n$.
Before introducing the $H_n(0)$-modules  
$\HnZmodC_{\alpha}$ and $\HnZmodR_{\alpha}$, we need to study the interval structure of $\RWC(\alpha)$ and $\RWC(\alpha)$.
The \emph{(left) weak Bruhat order} $\preceq_L$ on $\SG_n$ is defined to be the partial order on $\SG_n$ whose covering relation $\preceq_L^c$ is given as follows: 
\begin{align*}
\sigma \preceq_L^c s_i \sigma \ \text{ if and only if } \ i \notin \Des(\sigma).
\end{align*}
Given $\sigma, \rho \in \SG_n$, if $\sigma \preceq_L \rho$, then the \emph{(left) weak Bruhat interval} from $\sigma$ to $\rho$ is defined by
\[
[\sigma,\rho]_L := \{\gamma \in \SG_n \mid \sigma \preceq_L \gamma \preceq_L \rho \}.
\]
The {\em right Bruhat order}, denoted as $\preceq_R$, on $\SG_n$, and the {\em right weak Bruhat interval} $[\sigma,\rho]_R$ can likewise be defined analogously.

For $\sigma = \sigma_1 \sigma_2 \cdots \sigma_n \in \SG_n$, let
\begin{align*}
\inv{\sigma} &:= \{(i,j)  \mid 1 \leq i < j \leq n, \text{ but } \sigma_i > \sigma_j \} \quad \text{and} \\ 
\coinv{\sigma} &:= \{(i,j)  \mid 1 \leq i < j \leq n, \text{ and } \sigma_i < \sigma_j \}.
\end{align*}
From \cite[Proposition 3.1.3]{06BB} or \cite[Proposition 3.1]{91BW} it follows that 
\begin{equation}\label{Eq: weak order iff inv}
\sigma \preceq_L \rho \ \text{ if and only if }  \ \inv{\sigma} \subseteq \inv{\rho} \text{ if and only if } \ \coinv{\sigma} \supseteq \coinv{\rho}.
\end{equation} 

With this notation, we first observe that $\RWR(\alpha)$ is always an interval under the weak Bruhat order.
To be more precise, let $\Tsourcealpha$ be the SPIT of shape $\alpha$ filled with entries $1, 2, \ldots, n$ from bottom to top and from left to right and 
$\Tsinkalpha$ the SPIT of shape $\alpha$ whose entries in column $1$ are the first $\ell(\alpha)$ odd numbers, whose entries in column $2$ are the first $\ell(\alpha)-1$ or $\ell(\alpha)$ even numbers (depending on whether or not the last part of $\alpha$ is equal to 1), and whose entries in subsequent rows from top to bottom are the remaining numbers, increasing consecutively up each row.
For example, 
\begin{equation*}\label{Ex: source and sink of SPIT}
\begin{tikzpicture}[baseline=0mm]
\node at (0,0) {$\ctab{T^{\rm source}_{(3,2,4,1)}}{10 \\ 6 & 7 & 8 & 9 \\ 4 & 5  \\ 1 & 2 & 3}$};
\node at (3,0) {and };
\node at (6,0) {$\ctab{T^{\rm sink}_{(3,2,4,1)}}{7 \\ 5 & 6 & 8 & 9 \\ 3 & 4  \\ 1 & 2 & 10}$};
\end{tikzpicture}.
\end{equation*}
Then the following lemma shows that $\RWR(\alpha)$ is an interval under the weak Bruhat order.

\begin{lemma}\label{Lem: row reading word is interval}
Let $\alpha$ be a peak composition of $n$.
Then we have
\[
\RWR(\alpha) = [\rwRC{r}{\Tsinkalpha},\rwRC{r}{\Tsourcealpha}]_L.
\]
\end{lemma}
\begin{proof}
First, let us show $\RWR(\alpha) \subseteq [\rwRC{r}{\Tsinkalpha},\rwRC{r}{\Tsourcealpha}]_L$.
Let $T \in \SPIT(\alpha)$.
If $(i,j) \in \inv{\rwRC{r}{\Tsinkalpha}}$, then  $\rwRC{r}{\Tsinkalpha}_i$ is strictly above $\rwRC{r}{\Tsinkalpha}_j$ 
and $\rwRC{r}{\Tsinkalpha}_j$ appears in the first or second column of $\Tsinkalpha$.
For such a pair, since $T$ is an SPIT,
it follows from \cref{lem: peak condition iff canonical} that $\rwRC{r}{T}_i > \rwRC{r}{T}_j$, 
which implies that $(i,j) \in \inv{\rwRC{r}{T}}$.
Thus, we have that $\inv{\rwRC{r}{\Tsinkalpha}} \subseteq \inv{\rwRC{r}{T}}$.
Similarly, 
if $(i,j) \in \coinv{\rwRC{r}{\Tsourcealpha}}$, then $\rwRC{r}{\Tsourcealpha}_i$ and $\rwRC{r}{\Tsourcealpha}_j$ appear in the same row of $\Tsourcealpha$.
For such a pair, 
it follows from the definition of SPITs that  $\rwRC{r}{T}_i < \rwRC{r}{T}_j$, 
which implies that $(i,j) \in \coinv{\rwRC{r}{T}}$.
Thus, we have that $\coinv{\rwRC{r}{T}} \supseteq \coinv{\rwRC{r}{\Tsourcealpha}}$.
Consequently, 
we deduce that $\RWR(\alpha) \subseteq [\rwRC{r}{\Tsinkalpha},\rwRC{r}{\Tsourcealpha}]_L$.

Next, let us show $\RWR(\alpha) \supseteq [\rwRC{r}{\Tsinkalpha},\rwRC{r}{\Tsourcealpha}]_L$.
Let $\sigma \in [\rwRC{r}{\Tsinkalpha},\rwRC{r}{\Tsourcealpha}]_L$ and $T_\sigma$ be the filling of shape $\alpha$ such that $\rwRC{r}{T_\sigma} = \sigma$.
To prove the assertion,
it suffices to show that $T_\sigma \in \SPIT(\alpha)$.
Suppose that $x$ is immediately left of $y$ in $T_\sigma$.
Let $i_x$ and $i_y$ denote the indices such that $\rwRC{r}{T_\sigma}_{i_x} = x$ and $\rwRC{r}{T_\sigma}_{i_y} = y$, respectively.
Since $(i_x,i_y) \in \coinv{\rwRC{c}{\Tsourcealpha}}$, 
it follows from \eqref{Eq: weak order iff inv} that $(i_x,i_y) \in \coinv{\rwRC{c}{T_\sigma}}$, that is, 
\[
x = \rwRC{c}{T_\sigma}_{i_x}  < \rwRC{c}{T_\sigma}_{i_y} = y.
\]
So, the entries in each row are increasing from left to right.
Moreover,
if follows 
from $\coinv{T_\sigma} \supseteq \coinv{\Tsourcealpha}$ that $(T_\sigma)_{2,i} > (T_\sigma)_{1,i+1}$ for all $1 \leq i < \ell(\alpha)$.
This shows that $T_\sigma^{\, (\leq 2)}$ is canonical, and therefore, by \cref{lem: peak condition iff canonical}, we conclude that $T_\sigma \in \SPIT(\alpha)$.
\end{proof}

On the other hand, $\RWC(\alpha)$ generally does not form an interval under the left weak Bruhat order. 
Let $\Tsup{\alpha}$ be the SPIT of shape $\alpha$ whose entries in column $1$ are the first $\ell(\alpha)$ odd numbers, whose entries in column $2$ are the first $\ell(\alpha)-1$ or $\ell(\alpha)$ even numbers (depending on whether or not the last part of $\alpha$ is equal to 1), and whose entries in subsequent columns from left to right are the remaining numbers, increasing consecutively up each column.
For example, 
\begin{equation}\label{Ex: Tsupalpha}
\begin{tikzpicture}[baseline=0mm]
\node at (4,0) {$\ctab{T^{\rm sup}_{(3,2,4,1)}}{7 \\ 5 & 6 & 9 & 10 \\ 3 & 4  \\ 1 & 2 & 8}$};
\end{tikzpicture} \ .
\end{equation}
In \cite[Lemma 4.10]{22Searles}, it was established that $\rwRC{c}{\Tsupalpha}$ is the unique minimal element in $\RWC(\alpha)$ under the left weak Bruhat order. However, it may not have a unique maximal element.
The subsequent lemma provides a characterization of $\alpha$ for which $\RWC(\alpha)$ forms an interval.

\begin{lemma}\label{Lem: rw_C(SPIT) is interval or not}
Let $\alpha$ be a peak composition of $n$.
The set $\RWC(\alpha)$ is an interval under the weak Bruhat order if and only if the number of parts of $\alpha$ greater than $2$ is at most $1$.
\end{lemma}
\begin{proof}
Let us prove the `if' part.
When $\alpha$ has no parts greater than $2$, 
the assertion is trivial since  
$\SPIT(\alpha) = \{\Tsupalpha\}$. 
When $\alpha$ has only one part greater than $2$, one can prove that 
\[
\RWC(\alpha) = [\rwRC{c}{\Tsupalpha},\rwRC{c}{\Tsourcealpha}]_L
\]
by following the same approach used to prove \cref{Lem: row reading word is interval}.

The converse is established through a proof by contradiction. 
To do this, we assume that $\alpha$ has at least $2$ parts greater than $2$.
Let $T^{\rm s}_\alpha$ be the SPIT of shape $\alpha$ defined in  the following way:
\begin{itemize}
\item If $\alpha_{\ell(\alpha)} = 2$, then the entries in column $1$ are the first $\ell(\alpha)$ odd numbers, the entries in column $2$ are the first $\ell(\alpha)$ even numbers, and the entries in subsequent rows are the remaining numbers, increasing consecutively from left to right in each row starting from the top.

\item If $\alpha_{\ell(\alpha)} = 1$, then for the peak composition $(\alpha_1,\alpha_2,\ldots,\alpha_{\ell(\alpha)-1})$, the entries from $1$ to $n-1$ are filled in the same manner as described above.
Then, the entry $n$ is placed in the box $(1,\ell(\alpha))$.
\end{itemize}
It can be observed that for every $i \notin  \Des(\rwRC{c}{T^s_\alpha})$, $s_i\cdot T^s_\alpha \notin \SPIT(\alpha)$.
That is, $\rwRC{c}{T^s_\alpha}$ is a maximal element in $\{\rwRC{c}{T} \mid T \in \SPIT(\alpha)\} \cap \SG_n$ under the weak Bruhat order.
Put $c:=c_1(\alpha)+c_2(\alpha)$ and $d:=c+c_3(\alpha)$.
By the assumption, $d-c \geq 2$.
Then it follows from the definition of $T^s_\alpha$ and $\Tsourcealpha$ that
\begin{align*}
&(c+1,d) \in \inv{\rwRC{c}{T^s_\alpha}}, \text{ but } (c+1,d) \notin \inv{\rwRC{c}{\Tsourcealpha}},\\
&(c,c+1) \notin \inv{\rwRC{c}{T^s_\alpha}}, \text{ but } (c,c+1) \in \inv{\rwRC{c}{\Tsourcealpha}}.
\end{align*} 
Given \eqref{Eq: weak order iff inv}, this observation implies that $\rwRC{c}{T^s_\alpha}$ and $\rwRC{c}{\Tsourcealpha}$ are incomparable under the weak Bruhat order. 
Consequently, $\RWC(\alpha)$ cannot be an interval under the weak Bruhat order, contradicting the assumption.
\end{proof}

Let us now introduce the $H_n(0)$-modules in our consideration.
Define $\HnZmodC_{\alpha}$ to be the $H_n(0)$-module with $\C\SPIT(\alpha)$ as its underlying space, equipped with the following $H_n(0)$-action: 
For $1 \leq i \leq n-1$ and $T \in \SPIT(\alpha)$,
\begin{align}\label{Hn0-action star}
\opi_i \star T:= 
\begin{cases}
-T & \text{if } i \in \Des(\rwRC{c}{T}),\\
0 & \text{if $i \notin \Des(\rwRC{c}{T})$ and $s_i \cdot T \notin \SPIT(\alpha)$},\\
s_i \cdot T & \text{if } i \notin \Des(\rwRC{c}{T}) \text{ and } s_i \cdot T \in \SPIT(\alpha).
\end{cases}
\end{align}
This module was originally introduced by Searles and is cyclically generated by \( T^{\sup}_\alpha \) (see \cite[Theorem 4.11]{22Searles}).  

If \( \alpha \) is a peak composition of \( n \) with at most one part greater than \( 2 \), then \( \HnZmodC_{\alpha} \) is isomorphic to a negative weak Bruhat interval module as defined in \cite{22JKLO}. 
Given a weak Bruhat interval \( I \) in $\SG_n$, the {\em negative interval module} associated with \( I \), denoted by \( \osfB(I) \), is the left \( H_n(0) \)-module with underlying space \( \C I \) and the \( H_n(0) \)-action given by  
\begin{align}\label{Hn0-action another star}
\opi_i \star \gamma  := 
\begin{cases}
-\gamma & \text{if } i \in \Des{}{(\gamma)}, \\
0 & \text{if } i \notin \Des{}{(\gamma)} \text{ and } s_i\gamma \notin I, \\
s_i \gamma & \text{if } i \notin \Des{}{(\gamma)} \text{ and } s_i\gamma \in I
\end{cases} 
\end{align}
for \( i \in [n-1] \) and \( \gamma \in I \) (see \cite[Definition 1]{22JKLO}).
To describe the isomorphism, note that \( \RWC(\alpha) \) forms a weak Bruhat interval by \cref{Lem: rw_C(SPIT) is interval or not}. Consider the bijection  
\[
\mathrm{w}_{\rm c} :\SPIT(\alpha) \to  \RWC(\alpha), \quad T \mapsto \rwRC{c}{T}.
\]
A comparison of \eqref{Hn0-action star} and \eqref{Hn0-action another star} shows that this bijection induces an isomorphism  
from \( \HnZmodC_{\alpha} \) to \( \osfB(\RWC(\alpha)) \).

Next, the $H_n(0)$-module $\HnZmodR_{\alpha}$ is defined as follows.
By \cref{Lem: row reading word is interval}, the set \( \RWR(\alpha) \) forms a weak Bruhat interval. Consider the bijection  
\[
\mathrm{w}_{\rm r} : \SPIT(\alpha) \to \RWR(\alpha), \quad T \mapsto \rwRC{r}{T}.
\]
Using this bijection, the \( H_n(0) \)-module structure on \( \C\SPIT(\alpha) \) is defined by transferring the \( H_n(0) \)-module structure from \( \osfB(\RWR(\alpha)) \). The resulting \( H_n(0) \)-module is denoted by \( \HnZmodR_{\alpha} \).
Therefore its $H_n(0)$-action is given by
\begin{align}\label{Hn0-action dot}
\opi_i \actionS T:= 
\begin{cases}
-T & \text{if } i \in \Des(\rwRC{r}{T}),\\
0 & \text{if $i \notin \Des(\rwRC{r}{T})$ and $s_i \cdot T \notin \SPIT(\alpha)$},\\
s_i \cdot T & \text{if } i \notin \Des(\rwRC{r}{T}) \text{ and } s_i \cdot T \in \SPIT(\alpha)
\end{cases}
\end{align}
for $1 \leq i \leq n-1$ and $T \in \SPIT(\alpha)$.

Before presenting the main result of this subsection, a tableau \( \Tspalpha \) and the notion of colored digraph isomorphisms are introduced.
Let $\Tspalpha$ be the SPIT of shape $\alpha$ such that 
\begin{itemize}
\item the subdiagram $\tcd((2^{j_2}))$ of $\tcd(\alpha)$ is filled with $1,2,\ldots,2j_2$ from left to right and from bottom to top, and 

\item the entries at $(3,j_1)$ and $(3,j_2)$ are $2j_2+1$ and $2j_2+2$ respectively, with the remaining boxes filled with increasing consecutive numbers, arranged from left to right and from bottom to top. 
\end{itemize}
Next, consider weak Bruhat intervals $I_1, I_2$ in $\SG_n$. 
A map $f: I_1 \to I_2$ is called a \emph{colored digraph isomorphism} if $f$ is bijective and satisfies that for all $\gamma, \gamma' \in I_1$ and $1 \leq i \leq n-1$, 
\begin{align*}
\gamma \overset{i}{\rightarrow} \gamma' 
\quad \text{if and only if} \quad
f(\gamma) \overset{i}{\rightarrow} f(\gamma'),
\end{align*}
where 
$\gamma \overset{i}{\rightarrow} \gamma'$ denotes that  
$\gamma \preceq_L \gamma'$ and $s_i \gamma = \gamma'$. 
If there exists a descent-preserving colored digraph isomorphism between two intervals $I_1$ and $I_2$, then $\osfB(I_1)$ is isomorphic to $\osfB(I_2)$.
For details, see \cite[Section 4]{22JKLO}.

\begin{theorem}\label{Thm: CRWC is interval if alpha one part2}
Let $\alpha$ be a peak composition of $n$.
\begin{enumerate}[label = {\rm (\alph*)}]
\item $\ch(\HnZmodC_{\alpha})= \ch(\HnZmodR_{\alpha})$.

\item
If $\alpha$ has at most one part greater than $2$, then $\HnZmodC_{\alpha}=\HnZmodR_{\alpha}$ as $H_n(0)$-modules.
Otherwise, $\HnZmodC_{\alpha}\not\cong\HnZmodR_{\alpha}$ as $H_n(0)$-modules.

\item 
If $\alpha$ has at most one part greater than $2$, then there exists a descent-preserving colored digraph between $\RWC(\alpha)$ and $\RWR(\alpha)$. 
Moreover, the intervals $[\gamma, \xi \gamma]_L$'s are equivalent up to descent-preserving colored digraph isomorphism for all $\gamma \in [\rwRC{c}{\Tsupalpha}, \rwRC{r}{\Tsupalpha}]_R$, where $\xi={\rwRC{c}{\Tsourcealpha}}\rwRC{c}{\Tsup{\alpha}}^{-1}$.
\end{enumerate}
\end{theorem}
\begin{proof}
(a) Following the same way as in \cite[Section 5]{15TW}, we can see that 
\[
\ch(\HnZmodC_{\alpha}) = \sum_{\sigma \in \RWC(\alpha)} F_{\comp(\Des(\sigma))} \quad \text{ and } \quad 
\ch(\HnZmodR_{\alpha}) = \sum_{\sigma \in \RWR(\alpha)} F_{\comp(\Des(\sigma))}.
\]
Therefore, the assertion follows from \cref{Coro: new definition of QSQ}.

(b) 
To begin with, we note that $\HnZmodC_{\alpha}$
and $\HnZmodR_{\alpha}$ share the same underlying space $\C\SPIT(\alpha)$.
Now, we assume that $\alpha$ has at most one part greater than $2$.
By \cref{Rmk: Xi_alpha (T)=T for alpha}, we see that 
$\Des(\rwRC{c}{T})=\Des(\rwRC{r}{T})$ for every $T \in \SPIT(\alpha)$.  
Considering $H_n(0)$-actions given in \eqref{Hn0-action star} and \eqref{Hn0-action dot}, we observe that 
$\HnZmodC_{\alpha}$
and $\HnZmodR_{\alpha}$ have exactly the same $H_n(0)$-action.
As a consequence, $\HnZmodC_{\alpha}$ 
and $\HnZmodR_{\alpha}$ are identical as $H_n(0)$-modules.

Next, let us consider the case where $\alpha$ has more than one part greater than $2$.
Recall that for an $H_n(0)$-module $M$, the notation $\upchi[M]$ denotes the $\upchi$-twist of $M$, where $\upchi$ is the anti-involution of $H_n(0)$ defined by
\begin{align*}
&\upchi: H_n(0) \rightarrow H_n(0), \quad \opi_i \mapsto \opi_i \quad \text{for $1 \le i \le n-1$}.
\end{align*}
Note that since \( \upchi(ab) = \upchi(b)\upchi(a) \), it follows that \( \upchi(\opi_\sigma) = \opi_{\sigma^{-1}} \) for \( \sigma \in \SG_n \).
The underlying space of $\upchi[M]$ is the dual space $M^*$ of $M$, and 
its $H_n(0)$-action is given by  
\[
(b \cdot^\chi \delta)(m) := \delta(\chi(b) \cdot m) \quad \text{ for $b \in H_n(0)$, $\delta \in M^*$, and $m\in M$}
\]
(refer to \cite[Section 3.4]{22JKLO} for more details).
In the following, we will show that $\upchi[\HnZmodC(\alpha)]$ is not isomorphic to $\upchi[\HnZmodR(\alpha)]$ as $H_n(0)$-modules.

Combining \cite[Theorem 4 (3)]{22JKLO} with \cref{Lem: row reading word is interval}, we derive that the $H_n(0)$-module $\upchi[\HnZmodR_{\alpha}]$ is generated by $(\Tsourcealpha)^*$. 
Here, $\{T^* \mid T \in \SPIT(\alpha)\}$ denotes the dual basis of 
the basis $\SPIT(\alpha)$ for $\C\SPIT(\alpha)$.
On the other hand, according to \cref{Lem: rw_C(SPIT) is interval or not}, the module $\upchi[\HnZmodC_{\alpha}]$ is generated by $T^*$'s, where $T$'s are the SPITs 
satisfying the condition
\[
\opi_i \star T = 0 \text{ or } T \text{ for all }1 \leq i \leq n-1.
\]
It can be easily seen that $\Tsourcealpha$ and $T^{\rm s}_\alpha$ satisfy this condition.
For the definition of $T^{\rm s}_\alpha$, see the proof of \cref{Lem: rw_C(SPIT) is interval or not}.
Moreover, from \cite[Lemma 3.11]{15BBSSZ} it follows that  
$\Tsourcealpha$ is the unique tableau in $\SIT(\alpha)$ satisfying the condition $\Des(\rwRC{r}{T}) \subseteq \set(\alpha)$.
Since  $\SPIT(\alpha) \subseteq \SIT(\alpha)$ and 
$\left\{\Des(\sigma) \mid \sigma \in \RWC(\alpha)\right\}=\left\{\Des(\sigma) \mid \sigma \in \RWR(\alpha)\right\}$ by 
\cref{Coro: new definition of QSQ} (a), 
any $H_n(0)$-module homomorphism $\Gamma_\alpha:\upchi[\HnZmodR] \rightarrow \upchi[\HnZmodC]$ should map $(\Tsourcealpha)^*$ to itself.

Suppose that $j_1<j_2$ are the smallest two subindices in $\{1\le i \le \ell(\alpha) \mid \alpha_i>2\}$.
Then, let $\Tspalpha$ be the SPIT of shape $\alpha$ such that 
\begin{itemize}
\item the subdiagram $\tcd((2^{j_2}))$ of $\tcd(\alpha)$ is filled with $1,2,\ldots,2j_2$ from left to right and from bottom to top, and 

\item the entries at $(3,j_1)$ and $(3,j_2)$ are $2j_2+1$ and $2j_2+2$ respectively, with the remaining boxes filled with increasing consecutive numbers, arranged from left to right and from bottom to top. 
\end{itemize}
Consider the permutation $\rho_\alpha := \rwRC{r}{\Tspalpha} \cdot (\rwRC{r}{\Tsourcealpha})^{-1}$.
We observe that
\begin{equation}\label{Eq: pi_rho actionS_upchi Tsourcealpha}
\pi_{\rho_\alpha} \actionS^\upchi  (\Tsourcealpha)^* = (\Tspalpha)^*
\ \text{ and } \
\pi_{2j_2+1} \pi_{\rho_\alpha} \actionS^\upchi (\Tsourcealpha)^* = s_{2j_2+1} (\Tspalpha)^*.
\end{equation}
Note that from the definitions of $\Tsourcealpha$ and $\Tspalpha$, we deduce that
\[
\rho_\alpha = 
s_{[k_{j_1}-\alpha_{j_1}+1+q,k_{j_1}-\alpha_{j_1}+2]} 
s_{[k_{j_1}-\alpha_{j_1}+3+q,k_{j_1}-\alpha_{j_1}+3]} 
\, \cdots \, 
s_{[k_{j_1}+q,k_{j_1}]},
\]
where $k_{j_1}:=\alpha_1+\cdots+\alpha_{j_1}$, $q:=2(j_2-j_1)$, and $s_{[a,b]}:= s_{a} s_{a-1} \cdots s_{b}$ ($a > b$).
This expression tells us that  
\begin{equation}\label{Eq: pi_rho star_upchi Tsourcealpha}
\pi_{\rho_\alpha} \star^\upchi \Gamma_\alpha((\Tsourcealpha)^*) = (\Tspalpha)^*, 
\ \text{ and } \
\pi_{2j_2+1} \pi_{\rho_\alpha} \star^\upchi \Gamma_\alpha((\Tsourcealpha)^*) = (\Tspalpha)^*.
\end{equation}
Here, the second equality follows from the fact that the entry $2 j_2+1$ is strictly below $2 j_2 + 2$ in column $3$ of $\Tspalpha$.
Comparing \eqref{Eq: pi_rho actionS_upchi Tsourcealpha} with \eqref{Eq: pi_rho star_upchi Tsourcealpha} yields that  
\[
\upchi[\HnZmodR_\alpha] \not\cong \upchi[\HnZmodC_\alpha] \quad \text{ as $H_n(0)$-modules.}
\]

(c) We observe the equality $\Tsup{\alpha}=\Tsinkalpha$, which follows from the fact that $\tcd(\alpha)$ has only one box in each column $3,4,\ldots,\alphamax$.
Therefore we have  
\[
\RWC(\alpha) = [\rwRC{c}{\Tsup{\alpha}},\rwRC{c}{\Tsourcealpha}]_L \quad \text{and} \quad \RWR(\alpha) = [\rwRC{r}{\Tsup{\alpha}},\rwRC{r}{\Tsourcealpha}]_L.
\]
Next, we observe that the readings 
$\mathrm{w}_{\rm c}: \SPIT(\alpha) \to \RWC(\alpha)$ and $\mathrm{w}_{\rm r}:  \SPIT(\alpha) \to \RWR(\alpha)$ induce the following $H_n(0)$-isomorphisms: 
\begin{align*}  
&\mathrm{w}_{\rm c}: \HnZmodC_{\alpha} \to \osfB(\RWC(\alpha)), \quad T \mapsto  \rwRC{c}{T}\\
&\mathrm{w}_{\rm r}: \HnZmodR_{\alpha} \to \osfB(\RWR(\alpha)), \quad T \mapsto  \rwRC{r}{T}.
\end{align*}
From the commutative diagram
\[ \begin{tikzcd}
\HnZmodC_{\alpha} \arrow{r}{\rm id} \arrow[swap]{d}{\mathrm{w}_{\rm c} \,(\cong)} & \HnZmodR_{\alpha} \arrow{d}{\mathrm{w}_{\rm r}\,(\cong)} \\
\osfB(\RWC(\alpha)) \arrow{r}{\mathrm{w}_{\rm r}\circ {\mathrm{w}_{\rm c}}^{-1}}& \osfB(\RWR(\alpha))  
\end{tikzcd}
\]
it follows that 
\[
\mathrm{w}_{\rm r}\circ {\mathrm{w}_{\rm c}}^{-1}(\opi_i \cdot \mathrm{w}_{\rm c}(T))=\mathrm{w}_{\rm r}(\opi_i \star T)=\mathrm{w}_{\rm r}(\opi_i \actionS T)
=\opi_i \cdot \mathrm{w}_{\rm r}(T).
\]
Note that $\rwRC{c}{T} \overset{i}{\rightarrow} s_i \rwRC{c}{T}$ 
if and only if $i \notin \Des(\rwRC{c}{T})$ and 
$s_i \rwRC{c}{T} \in \RWC(\alpha)$. 
In the proof of (b), we showed that 
$\Des(\rwRC{c}{T})=\Des(\rwRC{r}{T})$. 
Thus, if $i \notin \Des(\rwRC{c}{T})$ and 
$s_i \rwRC{c}{T} \in \RWC(\alpha)$, then 
\[
\mathrm{w}_{\rm r}\circ {\mathrm{w}_{\rm c}}^{-1}(s_i \rwRC{c}{T})=\mathrm{w}_{\rm r}(s_i \cdot T)=\mathrm{w}_{\rm r}(\opi_i \actionS T)
=\opi_i \cdot \mathrm{w}_{\rm r}(T) = s_i \rwRC{r}{T},
\]
as required.
As a consequence, we conclude that 
\[
\mathrm{w}_{\rm r}\circ {\mathrm{w}_{\rm c}}^{-1}: \RWC(\alpha)\to \RWR(\alpha)
\]
is a descent-preserving poset isomorphism.

On the other hand, \cite[Proposition 4.1 and Theorem 4.5]{23KLO} says that the equivalence class $C$ up to descent-preserving colored digraph isomorphism containing $\RWC(\alpha)$ and $\RWR(\alpha)$ is of the form
\[
C = \{[\gamma,  \xi \gamma]_L \mid \gamma \in [\sigma_0, \sigma_1]_R\}
\]
for some right weak Bruhat interval $[\sigma_0, \sigma_1]_R$.
One can easily check that $\rwRC{c}{\Tsupalpha} \preceq_R \rwRC{r}{\Tsupalpha}$ by using induction on $|\alpha|$.
Since 
\[
[\rwRC{c}{\Tsup{\alpha}}, \rwRC{r}{\Tsup{\alpha}}]_R  \subseteq 
[\sigma_0, \sigma_1]_R,
\]
it follows that $[\gamma, \xi \gamma]_L$'s for all $\gamma \in [\rwRC{c}{\Tsup{\alpha}},\rwRC{r}{\Tsup{\alpha}}]_R$ are contained in $C$, and which proves the second assertion. 
\end{proof}

\begin{example}
(a) Let us consider $\alpha = (2,3,2)$.
Let  
\[
\ctab{T_1:}{5 & 6 \\ 3 & 4 & 7 \\ 1 & 2}, \quad 
\ctab{T_2:}{5 & 7 \\ 3 & 4 & 6 \\ 1 & 2}, \quad\text{and}\quad 
\ctab{T_3:}{6 & 7 \\ 3 & 4 & 5 \\ 1 & 2}.
\]
Then $\SPIT(\alpha) = \{T_1,T_2,T_3\}$.
By a direct calculation, one observes that $\ourmapS(T) = T$ for $T \in \SPIT(\alpha)$ (see also \cref{Rmk: Xi_alpha (T)=T for alpha}). 
\cref{Fig: actions on SPIT_alpha} illustrates their $H_7(0)$-actions on the bases for $\osfB(\RWC(\alpha))$ and $\osfB(\RWR(\alpha))$.  
From this one knows that $\HnZmodC_\alpha = \HnZmodR_\alpha$ as $H_n(0)$-modules.
In \cref{Fig: actions on SPIT_alpha}, the symbol
\begin{tikzpicture}[baseline=-1mm]
\def \hp {3em}
\def \vp {4em}
\node at (\hp*0, \vp*0) {} edge [out=40,in=320, loop] ();
\node at (\hp*0 + 0.7*\hp, \vp*0) {\tiny $\opi_i$};
\end{tikzpicture}
represents that $\opi_i$ acts as $-\id$.

\vspace*{-4mm}
\begin{figure}[ht]
\begin{tikzpicture}
\def \hp {1em}
\def \vp {14mm}
\def \hstep {35mm}
\node at (0,0) (A11) {$\rwRC{c}{T_1}$};

\node at (0,-\vp) (A21) {$\rwRC{c}{T_2}$};

\node at (0,-\vp*2) (A31) {$\rwRC{c}{T_3}$};
\draw[->] (A11) -- (A21) node[right,midway] {\scriptsize$\opi_6$};
\draw[->] (A21) -- (A31) node[right,midway] {\scriptsize $\opi_5$};
\node[right,xshift=\hp*1] at (A11) {} edge [out=50,in=-50, loop] ();
\node[right,xshift=\hp*3] at (A11) {\scriptsize $\opi_2,\opi_4$};
\node[right,xshift=\hp*1] at (A21) {} edge [out=50,in=-50, loop] ();
\node[right,xshift=\hp*3] at (A21) {\scriptsize $\opi_2,\opi_4,\opi_6$};
\node[right,xshift=\hp*1] at (A31) {} edge [out=50,in=-50, loop] ();
\node[right,xshift=\hp*3] at (A31) {\scriptsize $\opi_2,\opi_5$};
\draw[->] (A11) -- ($(A11)+(-\hp*2,-\vp*0.6)$) node [left,pos=0.3] {\tiny $\opi_1,\opi_3,\opi_5$} node[xshift=-\hp*0.3,yshift=-\vp*0.1] {$0$};
\draw[->] (A21) -- ($(A21)+(-\hp*2,-\vp*0.6)$) node [left,pos=0.3] {\tiny $\opi_1,\opi_3$} node[xshift=-\hp*0.3,yshift=-\vp*0.1] {$0$};
\draw[->] (A31) -- ($(A31)+(-\hp*2,-\vp*0.6)$) node [left,pos=0.3] {\tiny $\opi_1,\opi_3,\opi_5,\opi_6$} node[xshift=-\hp*0.3,yshift=-\vp*0.1] {$0$};
\node at (0+\hstep*2,-\vp*2) (C31) {$\rwRC{r}{T_3}$};
\node at (0+\hstep*2,-\vp*1) (C21) {$\rwRC{r}{T_2}$};
\node at (0+\hstep*2,-\vp*0) (C11) {$\rwRC{r}{T_1}$};
\draw[->] (C21) -- (C31) node[right,midway] {\scriptsize $\pi_5$};
\draw[->] (C11) -- (C21) node[right,midway] {\scriptsize $\pi_6$};
\node[right,xshift=\hp*1] at (C31) {} edge [out=45,in=305, loop] ();
\node[right,xshift=\hp*3] at (C31) {\scriptsize $\opi_2,\opi_5$};
\node[right,xshift=\hp*1] at (C21) {} edge [out=45,in=305, loop] ();
\node[right,xshift=\hp*3] at (C21) {\scriptsize $\opi_2,\opi_4,\opi_6$};
\node[right,xshift=\hp*1] at (C11) {} edge [out=45,in=305, loop] ();
\node[right,xshift=\hp*3] at (C11) {\scriptsize $\opi_2,\opi_4$};
\draw[->] (C11) -- ($(C11)+(\hp*-2,-\vp*0.6)$) node [left,pos=0.3] {\tiny $\opi_1,\opi_3,\opi_5$} node[xshift=\hp*-0.5,yshift=\vp*-0.2] {$0$};
\draw[->] (C21) -- ($(C21)+(\hp*-2,-\vp*0.6)$) node [left,pos=0.2] {\tiny $\opi_1,\opi_3$} node[xshift=\hp*-0.5,yshift=\vp*-0] {$0$};
\draw[->] (C31) -- ($(C31)+(\hp*-2,-\vp*0.6)$) node [left,pos=0.3] {\tiny $\opi_1,\opi_3,\opi_5,\opi_6$} node[xshift=\hp*-0.5,yshift=\vp*-0.1] {$0$};
\end{tikzpicture}
\caption{$\osfB(\RWC(\alpha))$ and $\osfB(\RWR(\alpha))$}
\label{Fig: actions on SPIT_alpha}
\end{figure}

(b) Let $\alpha = (3,3)$.
Then
\[
\ctab{\Tsourcealpha}{4 & 5 & 6 \\ 1 & 2 & 3}, \quad 
\ctab{\Tspalpha}{3 & 4 & 6 \\ 1 & 2 & 5}, \quad
\ctab{\Tsupalpha}{3 & 4 & 5 \\ 1 & 2 & 6}.
\]
Put $\ctab{T_1:}{3 & 5 & 6 \\ 1 & 2 & 4}$, and then we see that 
$\SPIT(\alpha) = \{\Tsourcealpha,\Tspalpha,\Tsupalpha,T_1\}$.
Consider $\upchi[\HnZmodR]$ and $\upchi[\HnZmodC]$.
\cref{Fig: actionS and star_upchi at alpha 33} illustrates their $H_6(0)$-actions on the bases.
In \cref{Fig: actionS and star_upchi at alpha 33}, the symbol
\begin{tikzpicture}[baseline=-1mm]
\def \hp {3em}
\def \vp {4em}
\node at (\hp*0, \vp*0) {} edge [out=40,in=320, loop] ();
\node at (\hp*0 + 0.7*\hp, \vp*0) {\tiny $\pi_i$};
\end{tikzpicture}
represents that $\pi_i$ act as $\id$.

\vspace*{-4mm}
\begin{figure}[ht]
\begin{tikzpicture}
\def \hp {1em}
\def \vp {15mm}
\def \hstep {16mm}
\node at (0-\hstep*5,-\vp*1.6) (A31) {${\Tspalpha}^*$};
\node at (0-\hstep*5,-\vp*0.8) (A21) {$T_1^*$};
\node at (0-\hstep*5,-\vp*0) (A11) {${\Tsourcealpha}^*$};
\node at (0-\hstep*5,-\vp*2.4) (A41) {${\Tsupalpha}^*$};
\draw[->] (A21) -- (A31) node[right,midway] {\scriptsize$\pi_4$};
\draw[->] (A11) -- (A21) node[right,midway] {\scriptsize $\pi_3$};
\draw[->] (A31) -- (A41) node[right,midway] {\scriptsize $\pi_5$};
\node[right,xshift=\hp*0.4] at (A31) {} edge [out=45,in=305, loop] ();
\node[right,xshift=\hp*2.1,yshift=-\vp*0.05] at (A31) {\scriptsize $\pi_1,\pi_3,\pi_4$};
\node[right,xshift=\hp*0.4] at (A21) {} edge [out=45,in=305, loop] ();
\node[right,xshift=\hp*2.1,yshift=-\vp*0.05] at (A21) {\scriptsize $\pi_1,\pi_3,\pi_5$};
\node[right,xshift=\hp*1.1] at (A11) {} edge [out=45,in=305, loop] ();
\node[right,xshift=\hp*2.7,yshift=-\vp*0.05] at (A11) {\scriptsize $\pi_1,\pi_2,\pi_4,\pi_5$};
\node[right,xshift=\hp*0.7] at (A41) {} edge [out=45,in=305, loop] ();
\node[right,xshift=\hp*2.5,yshift=-\vp*0.05] at (A41) {\scriptsize $\pi_1,\pi_3,\pi_4,\pi_5$};
\draw[->] (A21) -- ($(A21)+(\hp*-2,-\vp*0.5)$) node [left,pos=0.1] {\tiny $\pi_2$} node[xshift=\hp*-0.2,yshift=\vp*-0.2] {$0$};
\draw[->] (A31) -- ($(A31)+(\hp*-2,-\vp*0.5)$) node [left,pos=0.1] {\tiny $\pi_2$} node[xshift=\hp*-0.3,yshift=\vp*-0.2] {$0$};
\draw[->] (A41) -- ($(A41)+(\hp*-2,-\vp*0.5)$) node [left,pos=0.3] {\tiny $\pi_2$} node[xshift=\hp*-0.2,yshift=\vp*-0.2] {$0$};
\node at (0,-\vp*2) (B31) {${\Tspalpha}^*$};
\node at (0-\hstep,-\vp) (B21) {$T^*_1$};
\node at (0-\hstep,-\vp*0) (B11) {${\Tsourcealpha}^*$};
\node at (0+\hstep,-\vp*1) (B41) {${\Tsupalpha}^*$};
\draw[->] (B21) -- (B31) node[right,midway] {\scriptsize$\pi_4$};
\draw[->] (B11) -- (B21) node[right,midway] {\scriptsize $\pi_3$};
\draw[->] (B41) -- (B31) node[right,midway] {\scriptsize $\pi_5$};
\node[right,xshift=\hp*0.6] at (B31) {} edge [out=45,in=305, loop] ();
\node[right,xshift=\hp*2.3,yshift=-\vp*0.05] at (B31) {\scriptsize $\pi_1,\pi_3,\pi_4,\pi_5$};
\node[right,xshift=\hp*0.4] at (B21) {} edge [out=45,in=305, loop] ();
\node[right,xshift=\hp*2.1,yshift=-\vp*0.05] at (B21) {\scriptsize $\pi_1,\pi_3,\pi_5$};
\node[right,xshift=\hp*1.1] at (B11) {} edge [out=45,in=305, loop] ();
\node[right,xshift=\hp*2.7,yshift=-\vp*0.05] at (B11) {\scriptsize $\pi_1,\pi_2,\pi_4,\pi_5$};
\node[right,xshift=\hp*0.6] at (B41) {} edge [out=45,in=305, loop] ();
\node[right,xshift=\hp*2.1,yshift=-\vp*0.05] at (B41) {\scriptsize $\pi_1,\pi_3,\pi_4$};
\draw[->] (B21) -- ($(B21)+(\hp*-2,-\vp*0.6)$) node [left,pos=0.5] {\tiny $\pi_2$} node[xshift=\hp*-0.5,yshift=\vp*-0.2] {$0$};
\draw[->] (B31) -- ($(B31)+(\hp*0,-\vp*0.6)$) node [left,pos=0.5] {\tiny $\pi_2$} node[xshift=\hp*-0,yshift=\vp*-0.2] {$0$};
\draw[->] (B41) -- ($(B41)+(\hp*2,-\vp*0.6)$) node [left,pos=0.5] {\tiny $\pi_2$} node[xshift=\hp*0.5,yshift=\vp*-0.1] {$0$};
\end{tikzpicture}
\caption{The $\actionS^\upchi$-action and the $\star^\upchi$-action on the dual basis of $\SPIT(\alpha)$}
\label{Fig: actionS and star_upchi at alpha 33}
\end{figure}

\noindent
By definition, $\rho_\alpha = s_{[4,3]}$.
From this one observes that 
\begin{align*}
& \pi_{\rho_\alpha} \actionS^\upchi {\Tsourcealpha}^* = {\Tspalpha}^* \quad \text{and} \quad  \pi_{5}\pi_{\rho_\alpha} \actionS^\upchi {\Tsourcealpha}^* = {\Tsupalpha}^*,\\
& \pi_{\rho_\alpha} \star^\upchi {\Tsourcealpha}^* = {\Tspalpha}^* \quad \text{and} \quad  \pi_{5}\pi_{\rho_\alpha} \star^\upchi {\Tsourcealpha}^* = {\Tspalpha}^*,
\end{align*}
which implies that $\upchi[\HnZmodR_\alpha] \not\cong \upchi[\HnZmodC_\alpha]$ as $H_6(0)$-modules.
\end{example}

\begin{remark}\label{remark on the equivalence classes on intervals}
The classification of left weak Bruhat intervals in $\SG_n$ up to descent-preserving isomorphism was explored in \cite[Section 4]{23KLO}. Let $\alpha$ denote a peak composition with at most one part greater than $2$. 
In the proof of \cref{Thm: CRWC is interval if alpha one part2},  
we see that the intervals $[\gamma, \xi \gamma]_L$ in \cref{Thm: CRWC is interval if alpha one part2} (c)
belong to the same equivalence class, say $C$, within this classification. 
However, it is still unclear whether these intervals represent the complete set of intervals within $C$, or if there is an additional interval yet to be identified.
It would be very nice to obtain a thorough description of this class, which can be achieved by the identification of the minimal and maximal elements in $\{\sigma \mid [\sigma, \rho]_L \in C\}$ under the right weak Bruhat order.  
\end{remark}

\section{Transition matrices for \texorpdfstring{$\{\QSQ{\alpha}\}$ and $\{\PYQS{\alpha}\}$}{Lg}}

\subsection{The transition matrix from \texorpdfstring{$\{\QSQ{\alpha}\}$ to $\{\PYQS{\alpha}\}$}{Lg}} 
\label{A generating function-based interpretation}

In this subsection, we present a detailed description of the transition matrix from $\{\QSQ{\alpha} \mid \alpha \in \PCp \}$ to $\{\PYQS{\alpha}\mid \alpha \in \PCp\}$.
A crucial tool in accomplishing this is the insertion algorithm developed by Allen, Hallam, and Mason \cite{18AHM}, which was obtained through a slight modification of the Schensted insertion method outlined in \cite[Section 6]{11HLMW}. This algorithm inserts a positive integer $k$ into a YCT $T$, resulting in another YCT denoted as $k \rightarrow T$.
Using this algorithm, they proved that the transition matrix from $\{\DIF{\alpha} \mid \alpha \in \Cp \}$ to $\{\YQS{\alpha}\mid \alpha \in \Cp\}$ consists of nonnegative integers (see \cite{18AHM}).
We show that this algorithm can also be applicable in describing the expansion of $\QSQ{\alpha}$ in the basis $\{\PYQS{\alpha}\mid \alpha \in \PCp\}$.

For $\alpha = (\alpha_1,\alpha_2,\ldots,\alpha_{\ell(\alpha)})$, let $\bar{\alpha} := (\alpha_1 +1 , \alpha_2+1, \ldots, \alpha_{\ell(\alpha)} +1)$. 
Given a YCT $T$ of shape $\alpha$, let $\bar{T}$ be the filling of $\tcd({\bar{\alpha}})$
in which the rightmost entry in each row is $\infty$ and the remaining boxes have the same entry as $T$.
For example, 
\begin{equation}\label{Ex:T132}
\begin{tikzpicture}[baseline=0mm]
\node at (0,0) {$\ctab{T}{6 & 8 \\ 2 & 4 & 7 \\ 1}$};
\node at (2.5,0) {and };
\node at (5,0) {$\ctab{\bar{T}}{6 & 8 & \infty \\ 2 & 4 & 7 & \infty \\ 1 & \infty }$};
\end{tikzpicture}. 
\end{equation}
Then the procedure $k \rightarrow T$ is defined in the following manner.

\begin{algorithm}{\rm (\cite[Procedure 3.1]{18AHM})}\label{alg: insertion}
Let \((c_1,d_1), (c_2,d_2), \ldots\) be the boxes of \(\bar{T}\), listed column by column from right to left, with the boxes in each column ordered from top to bottom.

\begin{enumerate}[label = {\it Step \arabic*.}]
\item
Set $k_0:=k$ and let $i$ be the smallest positive integer such that $\bar{T}_{c_i-1,d_{i-1}} \le k_0 < \bar{T}_{c_i,d_i}$.
If such an $i$ exists, then go to {\it Step 2}.
Otherwise, go to {\it Step 3}.

\item Consider the following two cases:
\begin{itemize}
\item[{\bf Case 1:}]
If $\bar{T}_{c_i,d_i} =\infty$, then place $k_0$ in box $(c_i,d_i)$ and terminate the procedure.

\item[{\bf Case 2:}] 
If $\bar{T}_{c_i,d_i} \neq \infty$, then set $k:=\bar{T}_{c_i,d_i}$,
place $k_0$ in box $(c_i,d_i)$, and repeat the procedure by inserting $k$ into the sequence of boxes $(c_{i+1},d_{i+1}), (c_{i+2},d_{i+2}), \ldots$.
In such a situation, we say that $\bar{T}_{c_i,d_i}$ is {\em bumped}.
\end{itemize}

\item 
Begin a new row (containing only $k_0$) in the highest position in the leftmost column such that all entries below $k_0$ in the leftmost column are smaller than $k_0$ and terminate the procedure.  
If the new row is not the top row of the diagram, shift all higher rows up by one. 
\end{enumerate}
\end{algorithm}

\begin{example}\label{insertExample}
For $T$ given in \eqref{Ex:T132}, consider the process of obtaining \( 5 \rightarrow T \) from \( T \).
First, enumerate the boxes of $\bar{T}$ as follows:
\[
(4,2), (3,3), (3,2), (2,3), (2,2), (2,1), (1,3), (1,2), (1,1)
\]
The first element bumped from $\bar{T}$ is the $7$ in column $3$.  
This $7$ is replaced by $5$ and $7$ is then inserted into the remaining sequence of boxes.  
The $7$ then bumps the $8$ in column $2$ and $8$ is inserted into the remaining boxes.  
The $8$ is placed to the right of the $1$ and the process terminates.  
\[
\begin{tikzpicture}
\def \hp {5mm}

\node at (0,0) {$\ctab{\bar{T}}{6 & 8 & \infty \\ 2 & 4 & {\color{blue} 7} & \infty \\ 1 & \infty}$};
\draw[<-,color=red] (\hp*1.3,0.15) to (\hp*1.3,1) node[above] {{\color{red} $5$}};

\draw [->,decorate,decoration={snake,amplitude=.4mm,segment length=2mm,post length=1mm}] (\hp*3,0) -- (\hp*4,0);

\node at (\hp*6,0) {$ \tableau{6 & {\color{blue} 8} & \infty \\ 2 & 4 & 5 & \infty \\ 1 & \infty }$};
\draw[<-,color=red] (\hp*5.6,0.6) to (\hp*5.6,1) node[above] {{\color{red} $7$}};

\draw [->,decorate,decoration={snake,amplitude=.4mm,segment length=2mm,post length=1mm}] (\hp*8,0) -- (\hp*9,0);

\node at (\hp*11,0) {$\tableau{6 & 7 & \infty \\ 2 & 4 & 5 & \infty \\ 1 & {\color{blue} \infty}}$};
\draw[<-,color=red] (\hp*10.6,-0.25) to (\hp*10.6,1) node[above] {{\color{red} $8$}};

\draw [->,decorate,decoration={snake,amplitude=.4mm,segment length=2mm,post length=1mm}] (\hp*13,0) -- (\hp*14,0);
\node at (\hp*16,0) {$\tableau{6 & 7 & \infty \\ 2 & 4 & 5 & \infty \\ 1 & 8}$};
\end{tikzpicture}
\]
Therefore, we have that $\ctab{5 \rightarrow T}{6 & 7  \\ 2 & 4 & 5 \\ 1 & 8 }$.
\end{example}

Now we introduce the procedure in \cite{18AHM} that maps a word to a pair of fillings of the same shape.
Begin with $(P_0,Q_0)=(\emptyset, \emptyset)$, where $\emptyset$ is the empty filling.  
Let $w_1$ be the first letter in the word $w = w_1 w_2 \cdots w_k$.  
Insert $w_1$ into $P_0$ using \cref{alg: insertion} and let $P_1$ be the resulting YCT.  
Record the location in $P_1$ where the new box was created by placing a ``1'' in $Q_0$ in the corresponding location and let $Q_1$ be the resulting filling.  
Next, assume the first $j-1$ letters of $w$ have been inserted.  
Let $w_j$ be the $j$th letter in $w$.  
Insert $w_j$ into $P_{j-1}$ and let $P_j$ be the resulting filling.  
Place the letter $j$ in the box of $Q_{j-1}$ corresponding to the new box in $P_j$ created from \cref{alg: insertion} and let $Q_j$ be the resulting filling.
We denote the final pair $(P_k, Q_k)$ as $(P(w), Q(w))$. 
We refer to $P(w)$ as the \emph{insertion tableau} of $w$ and $Q(w)$ as the \emph{recording tableau} of $w$.

Next, we introduce some properties required to develop our arguments, starting with a review of the relevant definitions and notation in \cite[Section 3]{18AHM}.
Let $\beta$ be a composition of $n$.
Let $Q$ be a filling of a diagram for $\beta$ with the entries in $[n]$, each appearing exactly once. 
A {\em row strip} of length $k$ of $Q$ is a maximal sequence $a_1, a_2, \ldots, a_k$ of $k$ consecutive integers such that for all $1 \leq i < k$ in the row strip, $a_i+1$ appears strictly right of $a_i$ in $Q$. 
The {\em row strip shape} of $Q$ is the composition $(\alpha_1, \alpha_2,\ldots, \alpha_{})$, where $\alpha_i$ is the length of the row strip sequence which starts with the number $\alpha_1+\alpha_2+\cdots+\alpha_{i-1}+1$.
Finally, a {\em dual immaculate recording tableau} (DIRT) $Q$ of shape $\beta$ is defined as a filling of $\tcd(\beta)$ with the entries in $[n]$ such that 
\begin{enumerate}[label = {\rm (\arabic*)}]
\item the entries in each row are weakly increasing from left to right,

\item the row strips start in the first (or leftmost) column,

\item the entries in the first column are strictly increasing from bottom to top, and 

\item if $i < j$ and $Q_{k,j} > Q_{k,i}$, then $Q_{k,j} > Q_{k+1,i}$.
(If the box $(i+1, g)$ is empty, then we assume that it contains the entry
infinity.) 
Pictorially, this condition states that if 
\[
\begin{tikzpicture}
\def \hp {2.8em}
\def \vp {1.9em}
\node at (\hp*0.5,\vp*0.5) {{\tiny $(k,j)$}};
\node at (\hp*1.5,\vp*-1.5) {{\tiny $(k+1,i)$}};
\node at (\hp*0.5,\vp*-0.4) {$\vdots$};
\node at (\hp*0.5,\vp*-1.5) {{\tiny $(k,i)$}};

\draw (0,0) rectangle (\hp,\vp);
\draw (0,\vp*-2) rectangle (\hp*1,\vp*-1);
\draw (\hp,\vp*-2) rectangle (\hp*2,\vp*-1);

\node[right] at (\hp*2.2,\vp*-0.5) {\small with $Q_{k,j} > Q_{k,i}$, then $Q_{k,j} > T_{k+1,i}$.};
\end{tikzpicture}
\]
\end{enumerate}
Let $\DIRT(\beta,\alpha)$ be the set of DIRTs of shape $\beta$ and with row strip shape $\alpha$.

With the above terminologies, the properties of the procedure in \cite{18AHM} to be introduced are as follows:
\begin{itemize}
\item $P(w)$ is an SYCT for every $w \in \SG_n$.

\item $Q(w)$ is a dual immaculate recording tableau (DIRT) with a row strip shape $\alpha^\rmr$ if and only if $w$ appears as $\rwRC{r}{T}$ for some $T \in \SIT(\alpha)$. 

\item $\Des(\rwRC{r}{T}) = \BDes(P(\rwRC{r}{T}))$ holds for all $T \in \SIT(\alpha)$.
\end{itemize}
Based on these properties, Allen, Hallam, and Mason established the bijection
\[
\begin{tikzcd}[arrows={-Stealth}]
\AHM: \SIT(\alpha) \longrightarrow \underset{\beta \in \Cp_n} {\bigcup}\SYCT(\beta) \times \DIRT(\beta,\alpha^\rev), 
\quad T \mapsto (P(\rwRC{r}{T}),Q(\rwRC{r}{T})),
\end{tikzcd}
\]
which yields the expansion: 
\begin{equation}\label{eq: dI to YQS}
\DIF{\alpha} = \sum_{\beta \in \Cp_n} c_{\alpha, \beta} \YQS{\beta}\,,
\end{equation}
where $c_{\alpha, \beta} = |\DIRT(\beta,\alpha^\rmr)|$.
Furthermore, they showed that $c_{\alpha, \alpha} = 1$ and $c_{\alpha, \beta}$ vanishes unless $\beta$ is strictly less than $\alpha$ in the dominance order. 
Consequently, the transition matrix from $\{\DIF{\alpha}\}$ to $\{\YQS{\alpha}\}$ is unitriangular and comprises nonnegative integer entries (see \cite[Theorem 1.1 and Lemma 4.2]{18AHM}).
Combining the equality $\Des(\rwRC{r}{T}) = \BDes(P(\rwRC{r}{T}))$ with \cite[Lemma 5.3]{22Searles} yields that 
\begin{equation}\label{Eq: des wr SIT = des wY SYCT}
\Des(\rwRC{r}{T}) = \Des(\readSYCT(P(\rwRC{r}{T}))) \quad \text{for }T \in \SIT(\alpha).  
\end{equation}

In the following, we investigate what can be obtained if the domain of $\AHM$ is restricted to $\SPIT(\alpha)$. 
For the proof of our result, we need the following lemma.

\begin{lemma}{\rm (\cite[Lemma 3.5]{18AHM})}\label{Lem: Property1 for AHM algorithm}
Assume that $x \leq y$ are inserted into a YCT $T$ to form $(y \rightarrow (x \rightarrow T))$. 
The new box created by the insertion of $y$ is strictly to the right of the new box created by the insertion of $x$. 
In particular, if a sequence $x_1 \leq x_2 \leq \cdots \leq x_m$ is inserted into a YCT in order from smallest to largest, then the
column indices of the resulting new boxes strictly increase.
\end{lemma}

\begin{proposition}\label{prop: bijection}
Let $\alpha$ be a peak composition of $n$. 
Then the image of $\SPIT(\alpha)$ under \emph{$\AHM$} is given by $\bigcup_{\beta} \SPYCT(\beta) \times \DIRT(\beta,\alpha^\rev)$, and therefore
the map  
\[
\begin{tikzcd}[arrows={-Stealth}]
\emph{\AHM}|_{\SPIT(\alpha)}: \SPIT(\alpha) \longrightarrow \underset{\beta \in \PCp_n} {\bigcup}\SPYCT(\beta) \times \DIRT(\beta,\alpha^\rev)
\end{tikzcd}\]
is a one-to-one correspondence satisfying that $\Des(\rwRC{r}{T}) = \Des(\readSYCT(P(\rwRC{r}{T})))$ for all $T \in \SPIT(\alpha)$.
\end{proposition}
\begin{proof}
First, we show that \( {\AHM}|_{\SPIT(\alpha)} \) is well-defined.
Let $T \in \SPIT(\alpha)$.
Since $P(\rwRC{r}{T})$ is an SYCT, we have only to show that $P(\rwRC{r}{T})^{(\leq 2)}$ is canonical.
By \cref{lem: peak condition iff canonical}, for $1 \leq i < \ell(\alpha)$, it holds that $T_{1,i} < T_{2,i} < u$ for all $u \in T_{(\geq i+1)}$.
Here, $T_{(\geq i+1)}$ is the subtableau consisting of the rows $T_{(i+1)}, T_{(i+2)}, \ldots,$ $T_{(\ell(\alpha))}$.
Due to this inequality, the entry $T_{1,i}$ appears at the box $(1,1)$ in $T_{1,i} \rightarrow P(\rwRC{r}{T_{(\geq i+1)}})$, 
and $T_{2,i}$ appears at the box $(2,1)$ in $T_{2,i} \rightarrow (T_{1,i} \rightarrow P(\rwRC{r}{T_{(\geq i+1)}}))$.
Therefore, 
\[
(T_{2,i} \rightarrow (T_{1,i} \rightarrow P(\rwRC{r}{T_{(\geq i+1)}})))_{2,j} < (T_{2,i} \rightarrow (T_{1,i} \rightarrow P(\rwRC{r}{T_{(\geq i+1)}})))_{1,j+1}
\]
for all $1 \leq j < \ell(\alpha)-i+1$.
Repeating this process, we can derive from \cref{Lem: Property1 for AHM algorithm} that 
\[
P(\rwRC{r}{T})_{2,j} < P(\rwRC{r}{T})_{1,j+1} \quad \text{for all $1 \leq j < \ell(\alpha)$}, 
\]
thus $P(\rwRC{r}{T})^{(\leq 2)}$ is canonical.

Second, we show that ${\AHM}|_{\SPIT(\alpha)}$ is surjective.  
For each $T \in \SIT(\alpha)$, we claim that if $P(\rwRC{r}{T}))$ is an SPYCT, then $T \in \SPIT(\alpha)$.
Assuming that the claim holds, the surjectivity can be established as follows.  
Let \( (U, V) \in \SPYCT(\beta) \times \DIRT(\beta, \alpha^\rev) \) for some \( \beta \). Since \( \AHM \) is a bijection, there exists a unique tableau \( T = {\AHM}^{-1} (U, V) \in \SIT(\alpha) \) such that \( P(\rwRC{r}{T}) = U \). By the assumption, it follows that \( T \in \SPIT(\alpha) \), as required.

We now prove the claim by contradiction. We assume that $T$ does not satisfy the peak-tableau condition.
As $T^{(\leq 2)}$ is not canonical, we can choose the smallest $1 \leq i < \ell(\alpha)$ such that
\begin{equation}\label{Eq: T1i < T1i1 < T2i}
T_{1,i} < T_{1,i+1} < T_{2,i}, \quad \text{but} \quad 
T_{1,k} < T_{2,k} < T_{1,k+1} \quad \text{for all } 1 \leq k \leq i-1.
\end{equation}
When $i=1$, $P(\rwRC{r}{T})_{1,1} =1$ and $P(\rwRC{r}{T})_{1,2} =2$ since $T_{1,1} = 1$ and $T_{1,2}=2$.
This violates the peak-tableau condition, which contradicts the assumption.
Consider the case where $i \geq 2$.
The latter inequality implies that $P(\rwRC{r}{T})_{1,k} = T_{1,k}$ and $P(\rwRC{r}{T})_{2,k} = T_{2,k}$ for all $1 \leq k \leq i-1$. 
Consider the process
\[
P(\rwRC{r}{T_{(\ell(\alpha))}}) 
\ 
\begin{tikzpicture}
\draw [->,decorate,decoration={snake,amplitude=.4mm,segment length=2mm,post length=1mm}] (0,0) -- (6mm,0);
\end{tikzpicture}  
\
P(\rwRC{r}{T_{(\geq \ell(\alpha)-1}}) 
\
\begin{tikzpicture}
\draw [->,decorate,decoration={snake,amplitude=.4mm,segment length=2mm,post length=1mm}] (0,0) -- (6mm,0);
\end{tikzpicture}  
\ \cdots \ 
\begin{tikzpicture}
\draw [->,decorate,decoration={snake,amplitude=.4mm,segment length=2mm,post length=1mm}] (0,0) -- (6mm,0);
\end{tikzpicture}  
\
P(\rwRC{r}{T_{(\geq 1)}}) = P(\rwRC{r}{T}).
\]
Note that the entry $P(\rwRC{r}{T})_{2,i}$ remains stabilized after the $(\ell(\alpha)-i+1)$st step, with the exception of being shifted one box upwards at each step, as illustrated below:
\[
P(\rwRC{r}{T_{(\geq i)}})_{2,1} = P(\rwRC{r}{T_{(\geq i-1)}})_{2,2} = \cdots = P(\rwRC{r}{T})_{2,i}.
\]
If the box $(2,i)$ is not contained in the composition diagram of $P(\rwRC{r}{T})$, then $P(\rwRC{r}{T})$ is not an SPYCT. 
Otherwise, $P(\rwRC{r}{T})_{2,i}$ belongs to the set of entries of $T_{(\geq i)}$.
On the other hand, the entries $T_{1,i}+1, T_{1,i}+2, \ldots, T_{1,i+1}-1$ belong to the set of entries of $T_{(\leq i-1)}$ due to \eqref{Eq: T1i < T1i1 < T2i}.
This indicates that $P(\rwRC{r}{T})_{2,i} > T_{1,i+1} = P(\rwRC{r}{T})_{1,i+1}$, and thus $P(\rwRC{r}{T})^{(\leq 2)}$ is not canonical. 
Therefore, by \cref{lem: SPYCT canonical iff}, $P(\rwRC{r}{T})$ is not an SPYCT. This verifies the claim.

Finally, we show that \( {\AHM}|_{\SPIT(\alpha)} \) is injective and descent-preserving, which follows straightforwardly from the facts that \( \AHM \) is injective and descent-preserving.
\end{proof}

Now, we are ready to state the main result of this subsection, which asserts that the transition matrix from $\{\QSQ{\alpha}\}$ to $\{\PYQS{\beta}\}$ is derived from the transition matrix from $\{\DIF{\alpha}\}$ to $\{\YQS{\alpha}\}$ by restricting the index set from $\Cp$ to $\PCp$.

\begin{theorem}\label{Thm: QSQ into PYQS}
For a peak composition $\alpha$ of $n$, we have 
\[
\QSQ{\alpha} = \sum_{\beta \in \PCp_n} c_{\alpha, \beta} \PYQS{\beta},
\]
where $c_{\alpha, \beta} = |\DIRT(\beta,\alpha^\rmr)|$.
\end{theorem}
\begin{proof}
Combining \cref{Coro: new definition of QSQ} with \cref{prop: bijection}, we deduce that 
\[
\begin{aligned}
\QSQ{\alpha} &= \sum_{T\in \SPIT(\alpha)}K_{\comp(\Peak(\Des(\rwRC{r}{T})))} \\
&= \sum_{\beta \in \PCp_n} c_{\alpha, \beta} \, \left(\sum_{P \in \SPYCT(\beta)} K_{\comp(\Peak(\Des(\readSYCT(P))))}\right) 
= \sum_{\beta \in \PCp_n} c_{\alpha, \beta} \PYQS{\beta} \, ,    
\end{aligned}
\]
as desired.
\end{proof}

Recall that the transition matrix from $\{\DIF{\alpha}\}$ to $\{\YQS{\alpha}\}$ is unitriangular in the dominance order and consists of nonnegative integer entries.
\cref{Thm: QSQ into PYQS} verifies that this property holds for the transition matrix from $\{\QSQ{\alpha}\}$ to $\{\PYQS{\beta}\}$.

This expansion was also obtained by Searles and Slattery-Holmes~\cite{24SS}, using a different method.

\subsection{The expansion of \texorpdfstring{$\PYQS{\alpha}$ in the basis $\{\QSQ{\alpha}\}$ in special cases}{Lg}}
\label{The expansion of PYQS in the basis QSQ in special cases}

Unlike the transition matrix from $\{\DIF{\alpha} \mid \alpha \in \Cp\}$ to $\{\YQS{\alpha} \mid \alpha \in \Cp\}$, very little is known about its inverse until now. Nevertheless, some important observation and conjectures were presented in \cite[Section 6]{18AHM}. 
Specifically, the authors noted that this transition matrix has integer entries. 
Furthermore, they proposed the following conjectures:
\begin{itemize}
\item (\cite[Conjecture 6.1]{18AHM})
For $\alpha\in \Cp_n$, $\YQS{\alpha}=\sum_{\beta} b_{\alpha,\beta}\DIF{\beta}$ with $b_{\alpha,\beta}\in \{-1,0,1\}$. Moreover, $\sum_{\beta}b_{\alpha,\beta}=1$ if $\alpha=(1^k, n-k)$ and $0$ otherwise. 

\item (\cite[Conjecture 6.2]{18AHM})
If $\lambda$ is a strict partition, then $\YQS{\lambda} = \sum_{\sigma \in \SG_{\ell(\lambda)}} (-1)^{\ell(\sigma)}\DIF{\sigma(\lambda)}$, where $\sigma(\lambda)=(\lambda_{\sigma(1)}, \ldots, \lambda_{\sigma(\ell(\lambda))})$.
\end{itemize}
Motivated by these conjectures, we study the expansion of $\PYQS{\alpha}$ in the basis $\{\QSQ{\alpha}\}$ in special cases.

Firstly, we deal with the case where $\alpha$ is a strict partition. 
Before introducing our result, recall the notion of standard shifted tableaux given in \cite[Section 6]{89Stem}. 
A {\em standard shifted tableau} of shape $\lambda$ is a filling of the shifted Young diagram of $\lambda$ with the numbers $\{1,2,\dots,n\}$, where $n$ is the number of boxes, such that
\begin{enumerate}[label = {\rm (\arabic*)}]
\item  the entries increase strictly from left to right in each row, and   
\item the entries increase strictly from bottom to top in each column.  
\end{enumerate}
For example, if $\lambda = (4,2)$, then one possible standard shifted tableau is  
\[
\begin{ytableau}  
\none & 3 & 6 \\  
1 & 2 & 4 & 5
\end{ytableau} \, .
\]
Let $\SShT(\lambda)$ be the set of all standard shifted tableaux of shape $\lambda$.

The following proposition is a peak analogue of \cite[Conjecture 6.2]{18AHM}.

\begin{proposition} \label{peak analogue:strict partition} 
Let $\lambda$ be a strict partition.
Then we have
\[
\PYQS{\lambda} = \sum_{\substack{\alpha \in \PCp \\ \lambda(\alpha) = \lambda}} (-1)^{\ell(\sigma_\alpha)}\QSQ{\alpha},
\]
where $\sigma_\alpha$ is a unique minimal permutation in $\SG_{\ell(\lambda)}$ such that $\sigma_\alpha \cdot \alpha = \lambda(\alpha)$.
\end{proposition}
\begin{proof}
For any $S\in \SShT(\lambda)$, we have $i\in \Des(\rwRC{r}{S})$ if and only if $i$ is positioned strictly below $i+1$ in $S$.
It was shown in \cite[Lemma 6.8]{22Searles} that there is a bijection  
$\mathsf{rect}:\SShT(\lambda) \to \SPYCT(\lambda)$ 
such that 
\[
\Des(\rwRC{r}{S})=\Des(\readSYCT(\mathsf{rect}(S))).
\]
On the other hand, it was shown in \cite[Section 2.4]{97Stem} that 
\[
Q_\lambda = \sum_{S \in \SShT(\lambda)} K_{\comp(\Peak(\Des(\rwRC{r}{S})))},
\]
where $Q_\lambda$ is the {\em Schur $Q$-function} associated with $\lambda$.
From this it follows that $ Q_\lambda=\PYQS{\lambda}$.
Now, the assertion can be derived from the equality
\[
Q_\lambda = \sum_{\substack{\alpha \in \PCp \\ \lambda(\alpha) = \lambda}} (-1)^{\ell(\sigma_\alpha)} \QSQ{\alpha}
\]
(see \cite[Theorem 4.12]{15JL}).
\end{proof}

Secondly, we study the case where $\alpha$ is a peak composition with $\alpha_i \leq 3$ for all $1\le i\le \ell(\alpha)$.
We show that a peak analogue of \cite[Conjecture 6.1]{18AHM} holds in this case. 
For such a peak composition $\alpha$, define
\[
I_\alpha := \{ i \mid  1 \leq i \leq \ell(\alpha)-1 \text{ and } \alpha_i = 3 \text{ and } \alpha_{i+1} \neq 1 \} \, .
\]
For any interval $[i,j]$ contained in $I_\alpha$, we define $\Delta_{[i,j]}(\alpha)$ by the composition such that  
\begin{align*}
\Delta_{[i,j]}(\alpha)_k=
\begin{cases}\alpha_k-1 &\text{ if } i\le k \le j\\
\alpha_k+(j-i+1) &\text{ if } k =j+1\\
\alpha_k &\text{ otherwise.} 
\end{cases}
\end{align*}
Given a subset $S$ of $I_\alpha$, 
write it as the disjoint 
union of intervals
\[
S = [i_i, j_1] \cup [i_2, j_2] \cup \cdots \cup [i_s, j_s].  
\]
For $1\le k < l\le s$, one observes that 
\[
\Delta_{[i_k, j_k]}(\Delta_{[i_l, j_l]}(\alpha)) = \Delta_{[i_l, j_l]}(\Delta_{[i_k, j_k]}(\alpha)),
\]
thus this operation does not depend on the order of applying $\Delta_{[i_k, j_k]}$'s.
Therefore, we can define $\Delta_S(\alpha)$ recursively as follows:
\begin{equation}\label{Eq: Delta_S is defined recursively}
\Delta_S(\alpha) := \Delta_{[i_s, j_s]}( \cdots \Delta_{[i_2, j_2]}((\Delta_{[i_1, j_1]}(\alpha)))\cdots ).  
\end{equation}

\begin{example}\label{Ex:Delta_S(alpha)}
Let $\alpha = (3,2,3,3,1)$. 
Then $I_\alpha = \{1,3\}$ and 
its subsets are $\emptyset$, $\{1 \}$, $\{3 \}$, and $\{1, 3\}$.
By the definition of $\Delta_S(\alpha)$, 
one has that   
\begin{align*}
&\Delta_\emptyset(\alpha) = (3,2,3,3,1), \ \Delta_{\{1 \}}(\alpha) = (2,3,3,3,1), \ \Delta_{\{3 \}}(\alpha) = (3,2,2,4,1), \text{ and} \\ 
&\Delta_{\{1,3 \}}(\alpha) = \Delta_{\{3 \}}\Delta_{\{1\}}(\alpha)=\Delta_{\{1 \}}\Delta_{\{3\}}(\alpha)=(2,3,2,4,1).
\end{align*}
\end{example}

For each $s \in I_\alpha$,
define 
\[
\mathcal{A}_\alpha(s) := \{T \in \SPIT(\alpha) \mid T_{\alpha_s, s} > T_{\alpha_{s+1}, s+1}\}.
\]
And, let 
\[
\mathcal{A}_\alpha := \{T \in \SPIT(\alpha) \mid T_{\alpha_s, s} > T_{\alpha_{s+1}, s+1} \text{ for some } s \in I_\alpha\}.
\]

\begin{lemma}\label{lem:initial case}
Assume that $\alpha$ is a peak composition with $\alpha_i \leq 3$ for all $1\le i\le \ell(\alpha)$.
\begin{enumerate}[label = {\rm (\alph*)}]
\item Let $S \subseteq I_\alpha$.
Then, as multisets,
\[
\{\Des(\rwRC{r}{T}) \mid T \in \bigcap_{s \in S}\mathcal{A}_\alpha(S) \} = \{\Des(\rwRC{r}{T}) \mid T \in \SPIT(\Delta_S(\alpha))\}.
\]

\item As multisets, 
\[
\{\Des(\rwRC{r}{T}) \mid T \in \SPIT(\alpha) \setminus \mathcal{A}_\alpha\} = \{\Des(\readSYCT(T)) \mid T \in \SPYCT(\alpha)\}.
\]
\end{enumerate}
\end{lemma}
\begin{proof}
(a)
Due to \eqref{Eq: Delta_S is defined recursively}, we only need to prove the assertion for $S = [i,j] \subset I_\alpha$ with $i \leq j$.
Note that 
\[
\bigcap_{s \in S} \mathcal{A}_\alpha(s) = \{ T \in  \SPIT(\alpha) \mid T_{3,i} > \cdots > T_{3,j} > T_{\alpha_{j+1},j+1}\}.
\]
Consider the map
\[
\Phi_S: \bigcap_{s \in S} \mathcal{A}_\alpha(s) 
\to 
\SPIT(\Delta_S(\alpha)), \quad 
T \mapsto \bfs_{i,j+1} ( \cdots (\bfs_{j,j+1}(T))),
\]
where $\bfs_{p,q}(T)$ ($1 \leq p<q \leq \ell(\alpha)$) is the filling obtained from $T$ by shifting the box located at $(\beta_i,i)$ along with its entry to the immediately right of the rightmost box of row $j$ in $T$.
From the definition of $\bfs_{p,q}$, it follows that $\Phi_S$ is well-defined and bijective.

Next, let us show that
\begin{equation}\label{Eq: DesrT = DesrYsT}
\Des(\rwRC{r}{T}) = \Des(\rwRC{r}{\Phi_S(T)})
\end{equation}
for all $T \in \bigcap_{s \in [i,j]} \mathcal{A}_\alpha(s)$. 
For each $i \leq s \leq j$, let 
\begin{align*}
T[s] &:=\bfs_{s,j+1} (\cdots (\bfs_{j,j+1}(T))), \text{ and }\\ 
x_s &:=T[s]_{3,s}.
\end{align*}
For convenience, we set $T[j+1]:=T$. 
Using this notation, the equality \eqref{Eq: DesrT = DesrYsT} is equivalent to 
\begin{equation}\label{Eq: Des_r(T[s]) = Des_r(T[s+1])}
\Des(\rwRC{r}{T[s]}) = \Des(\rwRC{r}{T[s+1]}) \quad \text{for all $s=1,2,\ldots,j$}.
\end{equation}
For $i \leq s \leq j+1$, $T[s]$ is an SPIT, and thus, by definition,
\begin{equation}\label{Eq: T[s] is a peak-tableau}
T[s]_{\alpha_{j+1},j+1} < x_s.
\end{equation}
Except for the entry $x_s$, $T[s]$ and $T[s+1]$ are the same.
So, to establish the equality \eqref{Eq: Des_r(T[s]) = Des_r(T[s+1])}, it suffices to show that
\begin{itemize}
\item [(i)] $x_s-1$ does not appear between row $s$ and row $j$, and 
\item [(ii)]
$x_s+1$ does not appear between row $s$ and row $j+1$.
\end{itemize}
If $x_s-1$ appears between row $s$ and row $j$, then $x_s-1$ should be in column $1$ or column $2$, so it follows from \cref{lem: peak condition iff canonical} that $T[s]_{\alpha_{j+1},j+1} > x_s$. 
It contradicts \eqref{Eq: T[s] is a peak-tableau}, thus (i) holds.
If $x_s+1$ appears between row $s$ and row $j+1$, one can obtain a contradiction in a similar way of (i), thus (ii) holds.

(b)
Consider a map
\[
\Psi_\alpha: \SPYCT(\alpha) \rightarrow \SPIT(\alpha)\setminus \mathcal{A}_\alpha, \quad T \mapsto \Psi_\alpha(T),
\]
where $\Psi_\alpha(T)$ is the filling of shape $\alpha$ obtained from $T$ by sorting the entries of $T^{(3)}$ in increasing order from bottom to top, while keeping $T^{(1)}$ and $T^{(2)}$ unchanged. 
We claim that $\Psi_\alpha$ is a well-defined bijection, and $\Des(\readSYCT(T)) = \Des(\rwRC{r}{\Psi_\alpha(T)})$ for $T \in \SPYCT(\alpha)$.

Firstly, to show that $\Psi_\alpha(T) \in \SPIT(\alpha) \setminus \mathcal{A}_\alpha$, consider $T \in \SPYCT(\alpha)$. 
By definition, $\Psi_\alpha(T)^{(1)} = T^{(1)}$ and $\sort(T)^{(2)} = T^{(2)}$. According to the Young triple condition, $\Psi_\alpha(T)_{2,i} < \Psi_\alpha(T)_{1,i+1}$ for all $1 \leq i < \ell(\alpha)$, ensuring $\Psi_\alpha(T) \in \SPIT(\alpha)$. 
If $j$ is an index such that $\alpha_{j+1} = 2$ and $\alpha_{j} = 3$, then $T_{2,j} > T_{3,r}$ for any $r < j$ where $\alpha_r=3$, due to the Young triple condition. It follows that $T_{3,r_2} < T_{3,r_1}$ for any $r_1 < j < r_2$ where $\alpha_{r_1}=\alpha_{r_2}=3$. Consequently, $\sort(T)_{2,j+1} > \sort(T)_{3,j}$, which confirms that $\Psi_\alpha(T) \in \SPIT(\alpha) \setminus \mathcal{A}_\alpha$.

Secondly, we show that $\Psi_\alpha$ is a bijection. 
For $T, S \in \SPIT(\alpha)\setminus \mathcal{A}_\alpha$, it holds that if $T^{(3)}$ and $S^{(3)}$ have the same set of entries, then $T = S$.
Therefore the injectiveness of $\Psi_\alpha$ can be derived by showing that the same phenomenon happens for $T, S \in \SPYCT(\alpha)$.
Assume for contradiction that $T^{(3)}$ and $S^{(3)}$ have the same set of entries, but $T \neq S$.
Choose the largest index $i$ such that $T_{3,i} \neq S_{3,i}$.
Without loss of generality, we may assume that $T_{3,i} < S_{3,i}$. 
Then there exists an index $q < i$ such that $S_{3,q} = T_{3,i}$. 
Given that $T_{2,i} < T_{3,i}$, it follows that $S_{2,i} = T_{2,i} < S_{3,q} < S_{3,i}$, yielding a violation of the Young triple condition.
This contradicts the assumption that $S \in \SPYCT(\alpha)$.
We next derive the surjectiveness of $\Psi_\alpha$.
For $T \in \SPIT(\alpha) \setminus \mathcal{A}_\alpha$, let $S_T$ be the filling of shape $\alpha$ subject to the following conditions:
\begin{itemize}
\item The first column and second column of $S_T$ are those of $T$, respectively.

\item The third column of $S_T$ is constructed as follows:
Let $(x_1,x_2,\ldots,x_l)$ be the sequence of the entries in the third column of $T$ from bottom to top.
For $1 \leq p \leq l$, let $y_p$ be the uppermost entry in the second column of $S_T$ such that that $y_p<x_p$ and the box to the immediate right of $y_p$ is empty.
Then, place $x_p$ in the box to the immediate right of $y_p$.
\end{itemize}
The latter condition shows that $S_T$ satisfies the Young triple condition, therefore $S_T \in \SPYCT(\alpha)$.

Finally, we show that $\Des(\readSYCT(T)) = \Des(\rwRC{r}{\Psi_\alpha(T)})$ for $T \in \SPYCT(\alpha)$.
Observe that $i \in \Des(\readSYCT(T))$ if and only if $i+1,i \in T^{(3)}$ and $i+1$ is strictly below $i$.
Now, the desired equality follows from the definition of the map $\mathtt{sort}$.
\end{proof}

\begin{proposition}\label{Thm: peak Young Quasi Schur to quasi Schur Q}
Assume that $\alpha$ is a peak composition with $\alpha_i \leq 3$ for all $1\le i\le \ell(\alpha)$.
Then we have 
\begin{equation}\label{Eq: PYQS to QSQ}
\PYQS{\alpha} = \sum_{S \, \subseteq \, I_\alpha} (-1)^{|S|} \,  \QSQ{\Delta_S(\alpha)}.    
\end{equation}
In particular, the sum of all coefficients is $1$ if $I_\alpha =\emptyset$ and $0$ otherwise.
\end{proposition}
\begin{proof}
Note that 
\begin{equation}\label{Eq: QS = YQ + KA}
\begin{split}
&\sum_{S \in \SPIT(\alpha)} K_{\comp(\Peak(\Des(\rwRC{r}{S})))} \\ 
&=\sum_{T \in \SPIT(\alpha) \setminus \mathcal{A}_\alpha} K_{\comp(\Peak(\Des(\rwRC{r}{T})))} + \sum_{U \in \mathcal{A}_\alpha} K_{\comp(\Peak(\Des(\rwRC{r}{U})))}.
\end{split}
\end{equation}
From \cref{lem:initial case}(b) it follows that 
\[
\sum_{T \in \SPIT(\alpha) \setminus \mathcal{A}_\alpha} K_{\comp(\Peak(\Des(\rwRC{r}{T})))} = \PYQS{\alpha}.
\]
On the other hand, using the inclusion-exclusion principle, we can derive that 
\[
\mathbbm{1}_{\mathcal{A}_\alpha}(T) =
\sum_{\emptyset \neq S \subseteq I_\alpha} (-1)^{|S| +1} \mathbbm{1}_{\bigcap_{s \in S}\mathcal{A}_\alpha(s)}(T) \quad \text{ for $T \in \mathcal{A}_\alpha$.}
\]
Here, $\mathbbm{1}_X$ is the indicator function, 
that is, $\mathbbm{1}_X(T) = 1$ if $T \in X$ and $0$ otherwise. 
Combining this equality with \cref{lem:initial case}(a) yields that 
\[
\sum_{U \in \mathcal{A}_\alpha} K_{\comp(\Peak(\Des(\rwRC{r}{T})))} = \sum_{\emptyset \neq S \subseteq I_\alpha} (-1)^{|S|+1} \QSQ{\Delta_S(\alpha)}.
\]
Therefore, the assertion follows from that the left-hand side of \eqref{Eq: QS = YQ + KA} is $\QSQ{\alpha}$ (see  \cref{Coro: new definition of QSQ}(b)).
\end{proof}

We remark that the summation of the right-hand side of \eqref{Eq: PYQS to QSQ} is cancellation-free.

\begin{example}\label{example of raising operator}
Let $\alpha = (3,2,3,3,1)$. 
Combining \cref{Ex:Delta_S(alpha)} with \cref{Thm: peak Young Quasi Schur to quasi Schur Q} shows that
\[
\PYQS{(3,2,3,3,1)} = \QSQ{(3,2,3,3,1)} - \QSQ{(2,3,3,3,1)} - \QSQ{(3,2,2,4,1)} + \QSQ{(2,3,2,4,1)}.
\]
\end{example}

\begin{remark}
One may expect that \cref{Thm: peak Young Quasi Schur to quasi Schur Q} extends to general peak compositions. However, this is not the case. Let  \( \alpha = (3,3,4) \). Then \( I_\alpha = \{1,2\} \), whose subsets are \( \emptyset, \{1\}, \{2\}, \) and \( \{1,2\} \). By the definition of \( \Delta_S(\alpha) \),  
\[
\Delta_\emptyset(\alpha) = (3,3,4), \quad \Delta_{\{1\}}(\alpha) = (2,4,4), \quad \Delta_{\{2\}}(\alpha) = (3,2,5), \quad \text{and} \quad \Delta_{\{1,2\}}(\alpha) = (2,3,5).
\]  
Thus,  
\begin{equation*}
\sum_{S \subseteq I_\alpha} (-1)^{|S|} \QSQ{\Delta_S(\alpha)} = \QSQ{(3,3,4)} - \QSQ{(2,4,4)} - \QSQ{(3,2,5)} + \QSQ{(2,3,5)}.    
\end{equation*}  
On the other hand, the expansion of \( \PYQS{(3,3,4)} \) in terms of \( \QSQ{\beta} \) is given by  
\[
\PYQS{(3,3,4)} = \QSQ{(3,3,4)} - \QSQ{(2,4,4)} - \QSQ{(3,2,5)} + \QSQ{(2,2,6)}.
\]  
Since these expressions do not coincide, \cref{Thm: peak Young Quasi Schur to quasi Schur Q} does not hold for all peak compositions. A suitable generalization of \( I_\alpha \) that accounts for this discrepancy has yet to be determined.  
\end{remark}

Thirdly, we study the case where $\tcd(\alpha)$ is of hook-like shape within the context of peak compositions.
\begin{proposition}\label{Prop: peak hook shape1}
Assume that $\alpha$ is a peak composition with at most one part greater than $2$.
Then we have
\[
\PYQS{\alpha}  =
\begin{cases}
 \QSQ{\alpha} & \text{ if } \alpha=(2,2,\ldots, 2, k) \text{ or } \alpha=(2,2,\ldots, 2, k, 1) \text{ for } k \geq 2,\\
 \QSQ{\alpha} - \QSQ{s_i \cdot \alpha} & \text{ otherwise.}
\end{cases}
\]
Here, $i$ denotes the index such that $\alpha_i = k$.
\end{proposition}
\begin{proof}
Let $\alpha$ be a peak composition of shape $(2^a, k, 2^b)$ or $(2^a, k, 2^b, 1)$ for some $a, b \geq 0$ and $k \geq 2$.
First, assume that $b>0$ and $\alpha_i = k >2$.
Then there are two cases $T_{i+2, 2} > T_{i,3}$ and $T_{i+2, 2} < T_{i,3}$ where $T \in \SPIT(\alpha)$.
In the former case, since $P(\rwRC{r}{T}) = T$,  $T \in \SPYCT(\alpha)$.
By \eqref{Eq: des wr SIT = des wY SYCT} it follows that $\Des(\rwRC{r}{T}) = \Des(\readSYCT(T))$.
In the latter case, let $\hat{T}$ be the tableau obtained from $T$ by shifting the subfilling 
\begin{tikzpicture}[baseline = 2mm]
\def \hp {7mm}
\def \vp {6mm}
\draw (0,0) rectangle (\hp,\vp) node[xshift=-\hp*0.5,yshift=-\vp*0.5] {\footnotesize $T_{3,i}$};
\draw (\hp,0) rectangle (\hp*2,\vp) node[xshift=-\hp*0.5,yshift=-\vp*0.5] {\footnotesize $T_{4,i}$};
\draw (\hp*2,0) rectangle (\hp*3,\vp) node[xshift=-\hp*0.5,yshift=-\vp*0.5] {\footnotesize $\cdots$};
\draw (\hp*3,0) rectangle (\hp*4.1,\vp) node[xshift=-\hp*0.55,yshift=-\vp*0.5] {\footnotesize $T_{\alpha_i,i}$};
\end{tikzpicture}
upwards by one cell.
Then the map
\[ 
\{T \in \SPIT(\alpha) \mid T_{2,i+1} < T_{3,i}\} \rightarrow \SPIT(s_i \cdot \alpha), \quad T \mapsto \hat{T} 
\] 
is a bijection satisfying that $\Des(\rwRC{r}{T}) = \Des(\rwRC{r}{\hat{T}})$.
Combining both cases, we see that 
\[
\SPIT(\alpha) =  \SPYCT(\alpha) \cup \SPIT(s_i \cdot \alpha)
\]
and 
\begin{align*}
&\{\Des(\rwRC{r}{T}) \mid T \in \SPIT(\alpha)\} \\
&= \{\Des(\readSYCT(T)) \mid T \in \SPYCT(\alpha) \} \cup \{\Des(\rwRC{r}{T}) \mid T \in \SPIT(s_i \cdot \alpha) \}.
\end{align*}
Therefore we have $\QSQ{\alpha} = \PYQS{\alpha} + \QSQ{s_i \cdot \alpha}$.

Next, assume that $b = 0$ and $k \leq 2$.
Let $T \in \SPIT(\alpha)$. Since $P(\rwRC{r}{T})= T$, it follows that $T \in \SPYCT(\alpha)$, and therefore $\SPIT(\alpha) = \SPYCT(\alpha)$.
Now, the assertion follows from \eqref{Eq: des wr SIT = des wY SYCT}.
\end{proof}

Finally, we characterize peak compositions $\alpha$ such that $\PYQS{\alpha} = \QSQ{\alpha}$.

\begin{proposition}\label{characterization:PYQS=QSD}
Let $\alpha$ be a peak composition of $n$.
Then $\PYQS{\alpha} = \QSQ{\alpha}$ if and only if $\alpha = (2,2,\ldots, 2, k)$ or $(2,2,\ldots, 2, k, 1)$ for some $k\geq 2$.
\end{proposition}
\begin{proof}
The `if' part is already shown in \cref{Prop: peak hook shape1}.
We here show the `only if' part.
For the sake of contradiction, let us suppose that $\alpha$ is neither of the forms $(2^{\ell(\alpha)-1},k)$ nor $(2^{\ell(\alpha)-2)},k,1)$ for some $k \geq 2$.
For a such composition $\alpha$, we consider the SPIT $\Tsupalpha$ (see \eqref{Ex: Tsupalpha}). 
By the map $\AHM$ one can see that the shape of $P(\rwRC{r}{\Tsupalpha})$ is
\[
\beta:=
\begin{cases}
(2^{\ell(\alpha)-1}, n-2\ell(\alpha)+2)  & \text{if } \alpha_{\ell(\alpha)} = 2, \\ (2^{\ell(\alpha)-2}, n-2\ell(\alpha)+3,1) & \text{if } \alpha_{\ell(\alpha)} = 1.    
\end{cases}
\]
Therefore, considering \cref{Thm: QSQ into PYQS},
$c_{\alpha, \beta} \geq 1$ for $\beta \neq \alpha$ except when $\alpha$ equals either $(2,2,\ldots, 2, k)$ or $(2,2,\ldots, 2, k, 1)$ for some $k\geq 2$.
This contradicts the assumption that $\PYQS{\alpha} = \QSQ{\alpha}$, confirming the assertion.
\end{proof}

\begin{remark}
It was pointed out in~\cite[Section 4]{19Oguz} that $\QSQ{\alpha}$ may not necessarily exhibit positive expansion in the bases $\{\QS{\alpha} \mid \alpha \in \Cp\}$, $\{\YQS{\alpha} \mid \alpha \in \Cp\}$, and $\{\DIF{\alpha} \mid \alpha \in \Cp\}$. 
One can observe a similar phenomenon for peak Young quasisymmetric Schur functions.
For example, 
$\PYQS{(2,2,1)} = -\QS{(1^2, 2, 1)} + \QS{(1, 2^2)} + \QS{(1, 3, 1)} + \QS{(2, 1, 2)}$,  
$\PYQS{(2,2,1)}  =
- \YQS{(1, 2, 1^2)} + \YQS{(1, 3, 1)} + \YQS{(2, 1, 2)} + \YQS{(2^2,1)}$, and 
$\PYQS{(3,1)} =
- \DIF{(1, 3)} + \DIF{(2, 1^2)} + \DIF{(2^2)} + \DIF{(3, 1)}$.
\end{remark}

\section{Combinatorial properties of \texorpdfstring{$\SPIT(\alpha)$  and $\SPYCT(\alpha)$}{Lg}}
\label{Sec: combinatorial properties of two sets of tableaux}

\subsection{A hook length formula for \texorpdfstring{$\SPIT(\alpha)$}{Lg}}
\label{A hook length formula}

Recall that we have a surprising enumeration formula for standard immaculate tableaux that is analogous to the hook length formula for standard Young tableaux (see \cite[Proposition 3.13]{14BBSSZ}.
The aim of this subsection is to provide a peak analogue of this formula.

Let $\alpha$ be a peak composition and $k$ be a positive integer. 
In this subsection, we present a closed formula for the cardinality of $\SPIT(\alpha)$.
Given a box $(i,j) \in \tcd(\alpha)$,
we define
\[
\tilde{h}_\alpha((i,j)) :=
\begin{cases}
\alpha_j + \alpha_{j+1} + \cdots + \alpha_{\ell(\alpha)} & \text{if $i=1$,}\\
\alpha_j + \alpha_{j+1} + \cdots + \alpha_{\ell(\alpha)}-1 & \text{if $i=2$,} \\
\alpha_j - i +1 & \text{otherwise.}
\end{cases}
\]
When $i \leq 2$, $\tilde{h}_\alpha((i,j))$ counts the number of boxes above or weakly to the right of $(i,j)$ in $\tcd(\alpha)$. For $i \geq 3$, it counts the number of boxes weakly to the right of $(i,j)$.
With this notation, we can derive the following hook length formula for $\SPIT(\alpha)$.

\begin{proposition}\label{hook length formula for SPIT}
Let $\alpha$ be a peak composition of $n$.
Then we have 
\[
|\SPIT(\alpha)| = \frac{n!}{\prod_{c \in \tcd(\alpha)} \tilde{h}_\alpha(c)}\, .
\]
\end{proposition}
\begin{proof}
We will prove the assertion by applying mathematical induction on $\ell(\alpha)$.
If $\ell(\alpha) = 1$, the assertion holds, as both sides equal $1$.
Now, assume that the assertion holds for all peak compositions of $n$ of length $l$.
Let $\alpha$ be a peak composition of $n$ with $\ell(\alpha)=l+1$. 
For each $T \in \SPIT(\alpha)$, recall that $T_{(1)}$ denotes the first row of $T$, and let $T \setminus T_{(1)}$ denote the tableau obtained from $T$ by removing $T_{(1)}$. 
Note that the first two entries of $T_{(1)}$ are $1$ and $2$.
Consider the map 
\begin{align*}
\phi: \SPIT(\alpha) & \rightarrow \left\{U \subseteq [3,n] \mid  |U|=\alpha_1-2 \right\} \times \SPIT((\alpha_2, \ldots, \alpha_{\ell(\alpha)}))\\
T \hspace{5mm} & \mapsto  \hspace{12mm} (\{T_{1,j} \mid 3\le j \le \alpha_1\}, \stan(T \setminus T_{(1)})), 
\end{align*}
where $\stan(T \setminus T_{(1)})$ is the standardization of $T \setminus T_{(1)}$.
It can be easily seen that this map is a bijection, and therefore
\[
|\SPIT(\alpha)| = \binom{n-2}{\alpha_1-2} \cdot |\SPIT((\alpha_2, \ldots, \alpha_{\ell(\alpha)}))|.
\]
Now, the assertion follows from the induction hypothesis.
\end{proof}

\begin{example}
(a) Let $\alpha=(3,2,4,2)$. Then the filling $T$ of $\tcd(\alpha)$ with 
$T_{i,j}= \tilde{h}_\alpha((i,j))$ is given by
\[
\begin{ytableau}
2 & 1 \\
6 & 5 & 2 & 1 \\
8 & 7 \\
11 & 10 & 1
\end{ytableau} \ .
\]
By \cref{hook length formula for SPIT}, we see that 
$|\SPIT(\alpha)| = 54$.

(b) Let $\alpha=(3^k)$ for some positive integer $k$. 
By \cref{hook length formula for SPIT}, we see that 
$$|\SPIT(3^k)|=\prod_{1\le i \le k}(3i-2).$$
The sequence $\{a_n\}_{n \ge 0}$ given by $a_n=|\SPIT(3^{n+1})|$ 
coincides with the sequence \oeis{A007559} in \cite{OEIS}.
\end{example}

\subsection{Word models for \texorpdfstring{$\SYCT$s and $\SPYCT$s}{Lg}}
\label{SPYCT for rectangular shaped peak compositions}

Unlike standard peak immaculate tableaux, a closed-form enumeration formula for standard peak Young composition tableaux is currently unknown. 
In this subsection, we demonstrate that SYCTs and SPYCTs can be each bijectively mapped to words
satisfying suitable conditions.

To define the set of words in bijection with $\SYCT(\alpha)$ we collect some definitions.
A {\em lattice word} is a word composed of positive integers, in which every prefix contains at least as many positive integers $i$ as integers $i + 1$.
For a composition $\alpha$, let $c_i(\alpha)$ be the number of boxes in column $i$ of $\tcd(\alpha)$ and $\alpha_{\text{max}}$ the largest part of $\alpha$.

Given a lattice word \( w = w_1 w_2 \ldots w_l \) and a positive integer \( i \geq 1 \), let  
\[
w_{p_{i,1}} w_{p_{i,2}} \ldots w_{p_{i,r_i}}
\]  
denote the subword consisting of all occurrences of \( i \) in \( w \), where \( p_{i,1}, p_{i,2}, \ldots, p_{i,r_i} \) are the positions in \( w \) where \( i \) appears, and \( r_i \) is the total number of occurrences of \( i \) in \( w \).
And, define a word \( \text{rs}(w) = a_1 a_2 \ldots a_l \) recursively as follows:
\begin{itemize}
\item For \( i = 1 \) and \( 1 \leq k \leq r_1 \), set \( a_{p_{1,k}} := k \). 

\item For \( i \ge 2 \) and \( 1 \leq k \leq r_i \), assume that the values of \( a_{p_{i-1,j}} \) for \( 1 \leq j \leq r_{i-1} \) and \( a_{p_{i,j}} \) for \( 1 \leq j \leq k-1 \) have already been determined.  
Then, set \( a_{p_{i,k}} := \max A_k \), where  
\[
A_k = \{ a_{p_{i-1, j}} \mid 1 \leq j \leq r_{i-1} \text{ and } p_{i-1, j} < p_{i,k} \} \setminus \{ a_{p_{i,j}} \mid 1 \leq j \leq k-1 \}.
\]
\end{itemize}
Indeed, the entries are determined in the order: first, the entries indexed by \( p_{1,1}, p_{1,2}, \dots, p_{1,r_1} \), followed by those indexed by \( p_{2,1}, p_{2,2}, \dots, p_{2,r_2} \), and so on, continuing in this manner for each subsequent group indexed by \( p_{i,1}, p_{i,2}, \dots, p_{i,r_i} \).
For example, consider $w = 1 \ 1 \ 2 \ 3 \ 2 \ 1$.
The following figure illustrates how to obtain ${\rm rs}(w)$:
\[
\begin{array}{rcccccc}
& 1 {}_{\color{red} p_{1,1}} & . & . & . & .  & . \\
\begin{tikzpicture}
\draw [->,decorate,decoration={snake,amplitude=.4mm,segment length=2mm,post length=1mm}] (\hp*0,0) -- (\hp*2,0);    
\end{tikzpicture} & 1 {}_{\color{red} p_{1,1}} &  2 {}_{\color{red} p_{1,2}} & . & . &  . & .  \\
\begin{tikzpicture}
\draw [->,decorate,decoration={snake,amplitude=.4mm,segment length=2mm,post length=1mm}] (\hp*0,0) -- (\hp*2,0);    
\end{tikzpicture} &1 {}_{\color{red} p_{1,1}} & 2 {}_{\color{red} p_{1,2}} & . & . & . & 3 {}_{\color{red} p_{1,3}}   \\  
\begin{tikzpicture}
\draw [->,decorate,decoration={snake,amplitude=.4mm,segment length=2mm,post length=1mm}] (\hp*0,0) -- (\hp*2,0);    
\end{tikzpicture} & 1 {}_{\color{red} p_{1,1}} & 2 {}_{\color{red} p_{1,2}} & 2 {}_{\color{red} p_{2,1}} & .  & . & 3 {}_{\color{red} p_{1,3}}   \\     
\begin{tikzpicture}
\draw [->,decorate,decoration={snake,amplitude=.4mm,segment length=2mm,post length=1mm}] (\hp*0,0) -- (\hp*2,0);    
\end{tikzpicture} &1 {}_{\color{red} p_{1,1}} & 2 {}_{\color{red} p_{1,2}} & 2 {}_{\color{red} p_{2,1}} & . & 1 {}_{\color{red} p_{2,2}} & 3 {}_{\color{red} p_{1,3}}   \\     
{\rm rs}(w) = & 1 {}_{\color{red} p_{1,1}} & 2 {}_{\color{red} p_{1,2}} & 2 {}_{\color{red} p_{2,1}} & 2 {}_{\color{red} p_{3,1}} & 1 {}_{\color{red} p_{2,2}} & 3 {}_{\color{red} p_{1,3}}    
\end{array}
\]

With this preparation, we define \( {\setWord{\SYCT(\alpha)}} \) as the set of words \( w \) satisfying the following conditions:  
\begin{itemize}
\item[\rm (s1)] \( w \) is a lattice word,  
\item[\rm (s2)] For each \( 1 \leq i \leq \alphamax \), the letter \( i \) appears exactly \( c_i(\alpha) \) times in \( w \), and  
\item[\rm (s3)] For each \( 1 \leq j \leq \ell(\alpha) \), the letter \( j \) appears exactly \( \alpha_j \) times in \( {\rm rs}(w) \).  
\end{itemize}  
Then, consider the map 
\[
f_\alpha: \SYCT(\alpha) \to \setWord{\SYCT(\alpha)}, \quad T \mapsto f_\alpha(T),
\]  
where \( f_\alpha(T) \) is the word whose \( i \)th entry is the index of the column in \( T \) containing \( i \).

\begin{proposition} \label{words for SYCT}
For every composition $\alpha$ of $n$, 
the map $f_\alpha$ is bijective.
\end{proposition}
\begin{proof}
First, we show that \( f_\alpha \) is well-defined. 
For each \( T \in \SYCT(\alpha) \), there is a uniquely determined sequence of composition diagrams given by  
\begin{equation}\label{sequence of composition diagrams}
D_0 := \emptyset \subset D_1 \subset D_2 \subset \cdots \subset D_n = \tcd(\alpha),
\end{equation}  
where \( D_i \) is the shape of \( T[\leq i] \) for \( 1 \leq i \leq n \).  
This sequence ensures that \( f_\alpha(T) \) satisfies {\rm (s1)} and {\rm (s2)}.  
Furthermore, by the Young triple condition for SYCTs, for each \( 1 < i \leq n \), the new box must be placed in the uppermost available position in the column containing \( i \) in \( T \).  
This guarantees that \( f_\alpha(T) \) satisfies {\rm (s3)}.  

Next, we show that \( f_\alpha \) is injective.  
Suppose \( S, T \in \SYCT(\alpha) \) and \( S \neq T \).  
Then there exists an entry \( i > 1 \) such that for all \( 1 \leq j < i \), the column indices of \( j \) in \( S \) and \( T \) are identical, but the column index of \( i \) differs between \( S \) and \( T \).  
This implies that \( f_\alpha(S) \neq f_\alpha(T) \), proving injectivity.  

Finally, we show that \( f_\alpha \) is surjective.  
Let \( w = w_1 w_2 \ldots w_n \in \setWord{\SYCT(\alpha)} \).  
We construct a sequence of composition diagrams  
\[
D_0 = \emptyset \subset D_1 \subset \cdots \subset D_n,
\]
where  
\begin{itemize}
\item \( D_1 \) consists of a single box at \( (1,1) \), and  

\item for \( 1 < i \leq n \), \( D_i \) is obtained from \( D_{i-1} \) by adding a box at \( (w_i, {\rm rs}(w)_i) \).
\end{itemize}
By the definitions of \( w \) and \( {\rm rs}(w) \), we know that the final diagram \( D_n \) is equal to \( \tcd(\alpha) \).  
Moreover, this sequence satisfies the following conditions:  
\begin{itemize}
\item[\rm (i)] Boxes in each row are placed from left to right. 

\item[\rm (ii)] Boxes in the leftmost column are stacked from bottom to top.  

\item[\rm (iii)] If \( i < j \) and the box \( (k, j) \) is placed before \( (k+1, i) \), then \( (k+1, j) \in \tcd(\alpha) \) and is also placed before \( (k+1, i) \).  
\end{itemize}
Let \( T_w \) be the tableau by recording the order in which the boxes are placed.   
Since the conditions above ensure that $T_w$ satisfies the definition of an SYCT, we conclude that $f_\alpha(T_w) = w$.
Thus, $f_\alpha$ is surjective, completing the proof.
\end{proof}

\begin{example} \label{examples on WSYCT}
(a) 
The following illustrates the correspondence between $\SYCT((2,3))$ and $\setWord{\SYCT((2,3))}$:
\begin{gather*}
\ctab{T_1}{3 & 4 & 5 \\ 1 & 2} \longleftrightarrow 
\begin{cases}
f_\alpha(T_1) = 1 \ 2 \ 1 \ 2 \ 3 \\ 
{\rm rs}(f_\alpha(T_1)) = 1 \ 1 \ 2 \ 2 \ 2
\end{cases} 
\
\ctab{T_2}{2 & 3 & 4 \\ 1 & 5} \longleftrightarrow 
\begin{cases}
f_\alpha(T_2) = 1 \ 1 \ 2 \ 3 \ 2 \\  
{\rm rs}(f_\alpha(T_2)) = 1 \ 2 \ 2 \ 2 \ 1
\end{cases}
\\ 
\ctab{T_3}{2 & 3 & 5 \\ 1 & 4} \longleftrightarrow
\begin{cases}
f_\alpha(T_3) = 1 \ 1 \ 2 \ 2 \ 3 \\
{\rm rs}(f_\alpha(T_3)) = 1 \ 2 \ 2 \ 1 \ 2
\end{cases}
\hspace*{80mm}
\end{gather*}

(b) 
Let \( \alpha = (m^k) \) with \( m \geq 2 \) and \( k \geq 1 \). 
In the proof of \cref{words for SYCT}, the pair \( (w_i, {\rm rs}(w)_i) \) represents the position of the box \( D_i \) for each \( i \). 
Since \( \tcd(\alpha) \) has a rectangular shape, each number \( i \) appears exactly \( k \) times for \( 1 \leq i \leq m \). 
Consequently, the condition (s3) in the definition of \( \setWord{\SYCT(\alpha)} \) is automatically satisfied and can be omitted. 
It follows that \( \setWord{\SYCT(\alpha)} \) is the set of lattice words \( w = w_1 w_2 \dots w_{mk} \) in which each \( i \) appears exactly \( k \) times for \( 1 \leq i \leq m \).
In the case where $m=2$, it is well known that 
$\setWord{\SYCT(\alpha)}$ is a combinatorial model for the $k$th Catalan number $C_k:=\frac{1}{k+1}\binom{2k}{k}$ (for instance, see \cite{15Stan}). 
Therefore, from \cref{words for SYCT}, we deduce that $|\SYCT((2^k))|=\frac{1}{k+1}\binom{2k}{k}$.
\end{example}

Now, we introduce a peak analogue of \cref{words for SYCT}.  
A lattice word \( w \) is called a {\em peak lattice word} if it satisfies the condition  
\[
\#(\text{occurrences of } 1 \text{ in } w) \geq 2 \cdot \#(\text{occurrences of } 2 \text{ in } w).
\]  
Given a peak lattice word \( w = w_1 w_2 \cdots w_n \), define 
\[
{\rm ps}(w):= v_1 v_2 \ldots v_n,
\]
where 
\[
v_i=\begin{cases} 
1 & \text{ if $w_i=1$ and the count of 1's preceding $w_i$ is even},\\
2 & \text{ if $w_i=1$ and the count of 1's preceding $w_i$ is odd},\\
w_i+1 & \text{ if } w_i \ne 1
\end{cases}
\]
For instance, $1112312344$ is mapped to $1213423455$.

With this preparation, we define \( {\setWord{\SPYCT(\alpha)}} \) as the set of words \( w \) satisfying the following conditions:  
\begin{itemize}
\item[$(s1')$] \( w \) is a peak lattice word,  
\item[$(s2')$] $(c_1(\alpha)+c_2(\alpha))$ occurrences of $1$ and $c_{i+1}(\alpha)$ occurrences of $i$ for $2\le i \le \alphamax-1$, and  
\item[$(s3')$] For each \( 1 \leq j \leq \ell(\alpha) \), the letter \( j \) appears exactly \( \alpha_j \) times in \( {\rm rs}({\rm ps}(w)) \).  
\end{itemize}
Then, consider the map
\[
\tilde{f}_\alpha: \SPYCT(\alpha) \rightarrow \setWord{\SPYCT(\alpha)}, \quad T \mapsto \tilde{f}_\alpha(T),
\]
where $\tilde{f_\alpha}$ is the word whose $i$th entry is $1$ if $i$ is in the first or second column; 
otherwise, it is the index of the column of $T$ to which $i$ belongs, subtracted by $1$.
Define 
\[\bar \iota: \setWord{\SPYCT(\alpha)} \to \setWord{\SYCT(\alpha)}, \quad w \mapsto {\rm ps}(w).
\]
Clearly $\bar \iota$ is injective.
Using a similar approach as the proof of \cref{words for SYCT}, we can establish the following proposition. 

\begin{proposition}\label{prop: word model for spyct}
For every peak composition $\alpha$,  
the map $\tilde{f}_\alpha$ is bijective. 
Furthermore, the diagram 
\[ 
\begin{tikzcd}
\SPYCT(\alpha) \arrow{r}{\iota} \arrow[swap]{d}{\tilde{f}_\alpha} & \SYCT(\alpha)  \arrow{d}{f_\alpha} \\
\setWord{\SPYCT(\alpha)} \arrow{r}{\bar \iota}& \setWord{\SYCT(\alpha)}
\end{tikzcd}
\]
is commutative, where $\iota:\SPYCT(\alpha) \to \SYCT(\alpha)$ is the natural injection.
\end{proposition}

\begin{example}\label{Example: 3^k and 4^k}
Let $\alpha=(m^k)$ with $m\ge 3$ and $k\ge 1$.
In the same manner as in \cref{examples on WSYCT} (b), it follows that \( \setWord{\SPYCT(\alpha)} \) consists of words \( w = w_1 w_2 \dots w_{mk} \) in which \( 1 \) appears exactly \( 2k \) times and each \( i \) appears exactly \( k \) times for \( 2 \leq i \leq m-1 \), satisfying the following conditions for every prefix \( w_1 w_2 \dots w_k \) (\( k = 1,2, \dots, mk \)): 
\begin{align*}
&\sharp (\text{$1$ occurrences}) \ge 2\cdot \sharp (\text{$2$ occurrences}), \text{ and }\\
&\sharp (\text{$i$ occurrences}) \ge \sharp (\text{$i+1$ occurrences})\text{ for all } 2\le i \le m-2.
\end{align*}
\begin{itemize}
\item In the case where $m=3$, one can easily see that $\setWord{\SPYCT((3^k))}$ is in bijection with the set of the lattice paths originating from $(0,0)$ and terminating at $(2k,k)$ using steps from the set $\{(1,0),(0,1)\}$ that do not cross above the diagonal. 
For instance, when $k=2$, the word $111212$ is mapped to $\DYCKPsix{2}{0}{3}{0}{3}{1}{4}{1}$ \,.
It is well known that the latter set is a combinatorial model for the number $C_k^{(3)}:=\frac{1}{2k+1}\binom{3k}{k}$ (for instance, see \cite[Exercise 6.33(c)]{99Stanley} or \cite{19CD}).
The numbers $C_k^{(3)}$ ($k\ge 1$) are called \emph{the Fuss--Catalan numbers} or \emph{Raney numbers} (see \cite[A001764]{OEIS}). 
Therefore $|{\SPYCT((3^k))}|=\frac{1}{2k+1}\binom{3k}{k}$.

\item
In the case where $m=4$, under the assignment $1 \mapsto (1,0)$, $2 \mapsto (1,1)$, and $3 \mapsto (1,-1)$, $\setWord{\SPYCT((4^k))}$ is in bijection with the set ${\sf Sch}((4^k))$ of the Schr\"{o}der paths originating from $(0,0)$ and terminating at $(4k,0)$ using $(1,0)$ steps $2k$ times, $(1,1)$ steps $k$ times, and $(1,-1)$ steps $k$ times such that at each lattice point $P$ along the path, $a_P\ge 2b_P$ and $b_P\ge c_P$. 
Here, $a_P, b_P, c_P$ denote the numbers of $(1,0)$-steps, $(1,1)$-steps, and $(0,1)$-steps, respectively, taken to reach $P$ from $(0,0)$. 
The first six $|{\sf Sch}((4^k))|$'s are given by 
\[
1,\, 9, \, 153, \, 3579, \, 101630, \, 3288871.
\]
It should be mentioned that this sequence is already listed in \cite[A361190]{OEIS}, which enumerates the number of $4k$-step closed paths commencing and concluding at $(0,0)$ without crossing above the diagonal $y=x$ or below the $x$-axis, utilizing steps from the set $\{(1,1), (1,-1), (-1,0)\}$.
Establishing a bijection between ${\sf Sch}((4^k))$ and the set of paths mentioned above is quite straightforward. 
To elaborate, our Schr\"oder paths can be mapped onto the aforementioned closed paths by associating the $i$th step of the Schr\"oder path with the $(4k-i+1)$st step of the closed path, following the rule: 
$(1,0) \longleftrightarrow (-1,0)$, $(1,1) \longleftrightarrow (1,-1)$, and $(1,-1) \longleftrightarrow (1,1)$.
For instance, the following figure illustrates the connections among the SPYCT $T$, its corresponding lattice word, the associated Schr\"oder path, and the corresponding closed path.

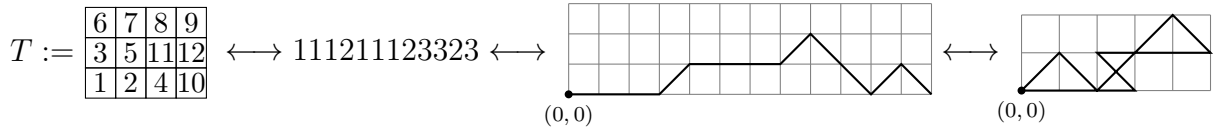
\[
\begin{tikzpicture}[baseline=-3mm]
\def \hp {4mm}
\def \vp {4mm}
\node[left] at (0,\vp*-0.4) {$T:=$};
\draw[draw=black] (0,0) rectangle (\hp,\vp) node at (0.5*\hp, 0.5*\vp) {\small 6};
\draw[draw=black] (\hp,0) rectangle (\hp*2,\vp) node at (1.5*\hp, 0.5*\vp) {\small 7};
\draw[draw=black] (\hp*2,0) rectangle (\hp*3,\vp) node at (2.5*\hp, 0.5*\vp) {\small 8};
\draw[draw=black] (\hp*3,0) rectangle (\hp*4,\vp) node at (3.5*\hp, 0.5*\vp) {\small 9};
\draw[draw=black] (0,\vp*-1) rectangle (\hp,\vp*0) node at (0.5*\hp, -0.5*\vp) {\small 3};
\draw[draw=black] (\hp,\vp*-1) rectangle (\hp*2,\vp*0) node at (1.5*\hp, -0.5*\vp) {\small 5};
\draw[draw=black] (\hp*2,\vp*-1) rectangle (\hp*3,\vp*0) node at (2.5*\hp, -0.5*\vp) {\small 11};
\draw[draw=black] (\hp*3,\vp*-1) rectangle (\hp*4,\vp*0) node at (3.5*\hp, -0.5*\vp) {\small 12};
\draw[draw=black] (0,\vp*-2) rectangle (\hp,\vp*-1) node at (0.5*\hp, -1.5*\vp) {\small 1};
\draw[draw=black] (\hp,\vp*-2) rectangle (\hp*2,\vp*-1) node at (1.5*\hp, -1.5*\vp) {\small 2};
\draw[draw=black] (\hp*2,\vp*-2) rectangle (\hp*3,\vp*-1) node at (2.5*\hp, -1.5*\vp) {\small 4};
\draw[draw=black] (\hp*3,\vp*-2) rectangle (\hp*4,\vp*-1) node at (3.5*\hp, -1.5*\vp) {\small 10};
\end{tikzpicture}
\longleftrightarrow 
1 1 1 2 1 1 1 2 3 3 2 3 
\longleftrightarrow 
\hspace*{-3mm}
\begin{tikzpicture}[baseline=4.5mm]
\def \hp {4mm}
\def \vp {4mm}
\foreach \c in {0,1,...,12}{
    \draw[-,black!50] (\hp*\c,\vp*0) -- (\hp*\c,\vp*3);
}
\foreach \c in {0,1,2,3}{
    \draw[-,black!50] (\hp*0,\vp*\c) -- (\hp*12,\vp*\c);
}
\draw[-,thick] (\hp*0,0) -- (\hp*3,\vp*0) -- (\hp*4,\vp*1) -- (\hp*7,\vp*1) -- (\hp*8,\vp*2) -- (\hp*9,\vp*1) -- (\hp*10,\vp*0) -- (\hp*11,\vp*1) -- (\hp*12,\vp*0);
\node[circle,fill=black,inner sep=0pt,minimum size=3pt] at (\hp*0,\vp*0) (Z0) {};
\node[below] at (Z0) {\tiny $(0,0)$};
\end{tikzpicture}
\longleftrightarrow
\hspace*{-3mm}
\begin{tikzpicture}[baseline=4mm]
\def \hp {5mm}
\def \vp {5mm}
\foreach \c in {0,1,2,3,4,5}{
    \draw[-,black!50] (\hp*\c,\vp*0) -- (\hp*\c,\vp*2);
}
\foreach \c in {0,1,2}{
    \draw[-,black!50] (\hp*5,\vp*\c) -- (\hp*0,\vp*\c);
}
\draw[-,thick] (\hp*0,0) -- (\hp*1,\vp*1) -- (\hp*2,\vp*0) -- (\hp*3,\vp*1) -- (\hp*4,\vp*2) -- (\hp*5,\vp*1) -- (\hp*4,\vp*1) -- (\hp*3,\vp*1) -- (\hp*2,\vp*1) -- (\hp*3,\vp*0) -- (\hp*2,\vp*0) -- (\hp*1,\vp*0) -- (\hp*0,\vp*0);
\node[circle,fill=black,inner sep=0pt,minimum size=3pt] at (\hp*0,\vp*0) (Z0) {};
\node[below] at (Z0) {\tiny $(0,0)$};
\end{tikzpicture}
\]

\end{itemize}
\end{example}

\subsection{Other bijections}\label{Other bijections}

Let $\lambda$ be a partition, and let $\SYT(\lambda)$ denote the set of standard Young tableaux of shape $\lambda$.
Consider the map,
\[
f_\lambda: \bigcup_{\substack{\alpha \in \Cp \\ \lambda(\alpha)=\lambda}} \SYCT(\alpha) \rightarrow \SYT(\lambda), \quad
T \mapsto \SPITtoYT{T},
\]
where $\SPITtoYT{T}$ is the standard Young tableau of shape $\lambda(\alpha)$ by placing the entries in the $i$th column of $T$ into the $i$th column of $\SPITtoYT{T}$ in the increasing order from bottom to top.
For example, if $T =$ {\footnotesize $\begin{ytableau} 5 & 6 & 7 \\ 3 & 4 & 8 \\ 1 & 2 \end{ytableau}$}, then  
$\SPITtoYT{T} =$ {\footnotesize $\begin{ytableau} 5 & 6 \\ 3 & 4 & 8 \\ 1 & 2 & 7 \end{ytableau}$}.
By making slight modifications to the arguments presented in \cite[Section 3 and Section 4]{08Mason} or \cite[Section 4]{13LMvW}, one can easily see that $f_\lambda$ is a bijection.
\footnote{In \cite[Section 3 and Section 4]{08Mason}, the authors consider standard reverse composition tableaux and standard reverse Young tableaux instead of standard Young composition tableaux and standard Young tableaux.}
Now, assume that $\lambda$ is a partition with at most one part equal to $1$, equivalently, a partition that is also a peak composition.
We refer to such a partition as a {\em peak partition}.
For simplicity, let 
\begin{align*}
\widehat{\SYCT(\lambda)}:={\bigcup_{\substack{\alpha \in \Cp \\ \lambda(\alpha)=\lambda}}  \SYCT(\alpha)}
\quad \text{and} \quad 
\widehat{\SPYCT(\lambda)}:={\bigcup_{\substack{\alpha \in \PCp \\ \lambda(\alpha)=\lambda}} \SPYCT(\alpha)},
\end{align*}
and consider the restriction of $f_\lambda$
\[
f_\lambda|: \widehat{\SPYCT(\lambda)} \rightarrow \SYT(\lambda), \quad
T \mapsto \SPITtoYT{T}.
\]
A standard Young tableau $T$ of shape $\lambda$ is referred to as a {\em standard peak Young tableau} (SPYT) if for every $1 \leq k \leq n$, the subdiagram of $\tcd(\lambda)$ consisting of boxes with entries at most $k$ in $T$ forms the diagram of a peak composition.
We denote by $\SPYT(\lambda)$ the set of standard peak Young tableaux of shape $\lambda$.

\begin{proposition}\label{Prop: f_lambda restriction and hook length formula}
Let $\lambda$ be a peak partition. 
\begin{enumerate}[label = {\rm (\alph*)}]
\item 
The image of $f_\lambda|$ equals $\SPYT(\lambda)$.

\item
Assume that each part of $\lambda$ is less than or equal to $3$.
Then we have 
\[
|\SPYT(\lambda)| = \frac{|\lambda|!}{\prod_{c \in \tcd(\lambda)}\hat{h}_\lambda(c)},
\]    
\end{enumerate}
where 
\[
\hat{h}_\lambda((i,j)) :=
\begin{cases}
\lambda'_1-j+\lambda'_2-j+\lambda_j & \text{if $i = 1$,}\\
\lambda'_1-j+\lambda'_2-j+\lambda_j-1 & \text{if $i = 2$,} \\
\lambda'_3-j+1 & \text{if $i = 3$.}
\end{cases}
\]
\end{proposition}
\begin{proof}
(a) Let $T \in \SYCT(\alpha)$ for some $\alpha$ with $\lambda(\alpha)=\lambda$. 
Note that for a column where the entries increase from bottom to top, $f_\lambda$ does not alter the relative positions of the entries.
Therefore, according to \cref{lem: SPYCT canonical iff}, if $T \in \SPYCT(\alpha)$, then $T$ and $\sort(T)$ share the first and second columns. 
Once more, according to \cref{lem: SPYCT canonical iff}, we observe that $\sort(T) \in \SPYT(\lambda)$.

Next, choose any tableau $U \in \SPYT(\lambda)$.
Since $\SPYT(\lambda) \subseteq \SYT(\lambda)$, $U$ is an SYT. 
Therefore there is a unique $T_U \in \SYCT(\beta)$ for some $\beta \in \Cp$ with $\lambda(\beta) = \lambda$ such that $f_\lambda(T_U)=U$.
Specifically, the tableau $T_U$ is constructed as follows:
\begin{itemize}
\item $(T_U)_{1,j} := U_{1,j}$ for $1 \leq j \leq \ell(\lambda)$.

\item For $i = 2, 3, \ldots, \alphamax$, let $(x_1,x_2,\ldots,x_q)$ be the sequence of the entries in column $i$ of $U$ from bottom to top.
For $p = 1,2,\ldots, q$, find the uppermost entry $y$ in column $j-1$ in $T_U$ such that 
$y < x_p$ and the box to the right of $y$ is empty.
Then, place $x_p$ in the box to the right of $T_U^{-1}(y)$.
\end{itemize}
(See \cite[Section 3 and Section 4]{08Mason} or \cite[Section 4]{13LMvW}).
Since $U \in \SPYT(\lambda)$, $U^{(\leq 2)}$ is canonical, meaning
\[
U_{1,1} < U_{2,1} < U_{1,2} < U_{2,2} < U_{1,3} < \cdots.
\]
Therefore, according to the construction of $T_U$, $U_{2,j}$ should be placed in the box $(2,j)$ in $T_U$. Consequently,
$(T_U)_{2,j} = U_{2,j}$ for all $1 \leq j \leq \ell(\lambda)$, indicating that $T_U$ is an SPYCT.

(b) We prove the first assertion by using induction on $|\lambda|$.
In the cases where $|\lambda|=1, 2$, or $3$, the assertion is straightforward. 

Let ${\rm PP}_{\le n}$ denote the set of peak partitions of size $ \leq n$ such that each part is at most $3$.
We now assume that the assertion is true for all $\lambda \in {\rm PP}_{n-1}$.
Suppose that $\lambda \in {\rm PP}_{\le n}\setminus{\rm PP}_{\le n-1}$.
If every part of $\lambda$ is at most two parts, that is, $\lambda = (2^a,1)$ or $(2^a)$, then one observes that $|\SPYT(\lambda)| = 1$.
The only tableau is obtained by filling integers $1,2,\ldots,n$ in this order from left to right and from bottom to top.
By a direct calculation, one can see that the hook length type formula gives the value $1$, as required.

Next, consider the cases where $\lambda$ has at least one part equal to $3$, that is, 
\[
\lambda=(3^{a}\, , 2^{b}) \text{ or }(3^{a} \, , 2^{b}\, , 1)
\text{ for some $a \geq 1$, and $b \geq 0$}.
\]
Let $T \in \SPYT(\lambda)$. To begin with, we observe that 
the subtableau $T[\leq n-1]$ of $T$, obtained by excluding the box filled with $n$, can be viewed as an SPYT within $\SPYT(\mu)$ for some $\mu \in {\rm PP}_{n-1}$.
From this observation, it follows that 
\[
|\SPYT(\lambda)|=\sum_{\text{some $\mu$'s in ${\rm PP}_{n-1}$}}|\SPYT(\mu)|.
\]
Here, we will deal with the most intricate case: 
$\lambda = (3^{a} \, , 2^{b}\, , 1)$ for $b \geq 1$.
The remaining cases can be verified in the same manner.
Since the entries in each column of $T$ are increasing from bottom to top,
the shape of $T[\leq n-1]$ aligns with one of the peak compositions
\[
\mu(1): = (3^{a-1} \, , 2^{b+1}\, , 1), \quad  
\mu(2) := (3^{a} \, , 2^{b-1}\, ,1, 1), \quad 
\mu(3) := (3^{a} \, , 2^{b}).
\]
However, $\mu(2)$ is not a peak composition, leading to the equality
\[
|\SPYT(\lambda)| = |\SPYT(\mu(1))| + |\SPYT(\mu(3))| \, .
\]
In the case where $T[\leq n-1] \in \SYCT(\mu(1))$, discrepancies in hook lengths between $\tcd(\lambda)$ and $\tcd(\mu(1))$ are confined to the $a$th row and third column, thus we see that 
\[
\prod_{c \in \tcd(\lambda)}\hat{h}_\lambda(c) 
= \frac{(n-3(a-1))a}{n-3a+1} \prod_{c \in \tcd(\mu(1))}\hat{h}_{\mu(1)}(c) \, .
\]
On the other hand, in the case where $T[n-1] \in \SYCT(\mu(3))$, discrepancies in hook lengths between $\tcd(\lambda)$ and $\tcd(\mu(3))$ are confined to the first and second column, thus we see that 
\[
\prod_{c \in \tcd(\lambda)}\hat{h}_\lambda(c) 
= \frac{(n-a+1)(n-3a)}{n-3a+1} \prod_{c \in \tcd(\mu(3))}\hat{h}_{\mu(3)}(c) \, .
\]
Merging these equations under the induction hypothesis shows that
\begin{align*}
|\SPYT(\lambda)|  &= \frac{(n-1)!}{\prod_{c \in \tcd(\mu(1))}\hat{h}_{\mu(1)}(c)} \, + \, \frac{(n-1)!}{\prod_{c \in \tcd(\mu(3))}\hat{h}_{\mu(3)}(c)}   \\
&=  \frac{(n-1)!}{\prod_{c \in \tcd(\lambda)}\hat{h}_\lambda(c)} \left( \frac{(n-3(a-1))a}{n-3a+1} \, + \, \frac{(n-a+1)(n-3a)}{n-3a+1}\right) \\
&= \frac{(n-1)!}{\prod_{c \in \tcd(\lambda)}\hat{h}_\lambda(c)} \cdot \frac{n(n-3a+1)}{n-3a+1} \, ,
\end{align*}
as required.
\end{proof}

\begin{remark}
It is worth noting that \cref{Prop: f_lambda restriction and hook length formula}(b) may not apply to all peak partitions.
For example, when $\lambda = (5,3)$, we have $\frac{|\lambda|!}{\prod_{c \in \tcd(\lambda)}\hat{h}_\lambda(c)} = 20$, while $|\SPYT(\lambda)| = 19$.
\end{remark}

From now on, we assume that $\lambda$ is a peak partition with $\ell(\lambda) \geq 2$ and no parts equal to $1$.
The rest of this subsection is devoted to showing that $\widehat{\text{SPYCT}(\lambda)}$ is in bijection with the set of standard set-valued Young tableaux satisfying suitable conditions.
For any partition $\lambda=(\lambda_1,\lambda_2,\ldots, \lambda_k)$ and a composition $\rho=(\rho_1, \rho_2, \ldots, \rho_k)$, a \emph{standard set-valued Young tableau} of shape $\lambda$ and density $\rho$ is a filling $\tcd(\lambda)$ with distinct entries in $[\,\sum_{i=1}^k\lambda_i\rho_i\,]$ such that 
\begin{enumerate}[label = {\rm (\arabic*)}]
\item every box in the $j$th row contains precisely $\rho_j$ distinct integers, 

\item for each $(i,j), (i,j+1)\in \tcd(\lambda)$,
every integer in $(i,j)$ is smaller than every integer in $(i,j+1)$, and 

\item for $(i,j), (i+1,j)\in \tcd(\lambda)$,
every integer in $(i,j)$ is smaller than every integer in $(i+1,j)$.
\end{enumerate}
We denote by $\SVT(\lambda, \rho)$ the set of all standard set-valued Young tableaux of shape $\lambda$ and density $\rho$.

We present the operations required to establish our bijection.
Given a standard Young tableau $T$ of any shape, 
we define $\mergef_{\rm c}(T)$ (resp., $\mergef_\rmr(T)$) as the tableau obtained from $T$ by merging the first and second columns (resp., rows) in such a way that the box $(1,i)$ (resp., the box $(i,1)$) of the new filling is filled with the entries $T_{1,i}$ and $T_{2,i}$ (resp., the entries $T_{i,1}$ and $T_{i,2}$) if they exist. 
Note that if $T$ is of shape $\alpha = (\alpha_1,\alpha_2,\ldots,\alpha_{\ell(\alpha)})$, then $\mergef_{\rm c}(T)$ is of shape $\beta$ with $\beta_i = \alpha_i+\delta_{\alpha_i,1}-1$ for $1 \leq i \leq \ell(\alpha)$ and $\mergef_\rmr(T)$ is of shape $(\max\{\alpha_1,\alpha_2\},\alpha_2,\ldots,\alpha_{\ell(\alpha)})$, where $\delta_{i,j}$ is the {\em Kronecker delta function}.
Next, given a standard Young tableau or standard set-valued Young tableau $T$ of any shape, we define $\reflectf(T)$ as the tableau obtained from $T$ by reflecting it along the line $y=x$.
When the tableau shape is fixed, the first two operations yield functions $\mergef_{\rm c}, \mergef_{\rm r}$ from the set of standard Young tableaux of the given shape to the set of standard set-valued Young tableaux of a certain shape and density.
Similarly, the third operation yields a function $\reflectf$ from the set of standard Young tableaux of the given shape (resp., the set of standard set-valued Young tableaux of the given shape and density) to the set of standard Young tableaux of a certain shape (resp., the set of standard set-valued Young tableaux of a certain shape and density). 
For simplicity, we omit specifying the domains and ranges of these functions.

Let $\lambda$ be a peak partition whose last part is greater than $1$, and let $\lambda'_{-}$ be the partition $(\lambda'_1,\lambda'_3,\lambda'_4, \ldots,\lambda'_{\lambda_1})$, where $\lambda'=(\lambda'_1,\lambda'_2,\ldots, \lambda'_{\lambda_1})$ is the conjugate of $\lambda$.
For each positive integer $s$, we define $\rho_{s}$ as the composition $(2,1,1,\ldots,1)$ of length $s$.
Now, we consider the map 
\[
g_\lambda: \SPYT(\lambda) \longrightarrow \SVT(\lambda'_{-}, \rho_{\lambda_1-1}), \quad T \mapsto \reflectf \circ \mergef_{\rm c} (T).
\]
For instance, we have 
\[
\begin{tikzpicture}[baseline=-1em]
\def \hp {1.5em}
\def \vp {1.1em}
\def \hstep {10em}
\node at (\hstep*0.9,0) (A1) {$\ctab{T}{
8 & 9 & \text{\footnotesize $10$} &  \\
3 & 4 & 7 & \text{\footnotesize $12$} & \\
1 & 2 & 5 & 6 & \text{\footnotesize $11$} }$};

\node at (\hstep*2,0) (A2) {
$\begin{tikzpicture}[baseline=-2em]
\def \hpp {1.3em}
\def \vpp {1.1em}
\foreach \r in {0,1}{
    \draw (\hp*\r,\vp*0) rectangle (\hp*\r+\hp,\vp*0+\vp);
}
\foreach \r in {0,1,2}{
    \draw (\hp*\r,\vp*-1) rectangle (\hp*\r+\hp,\vp*-1+\vp);
}
\foreach \r in {0,1,2,3}{
    \draw (\hp*\r,\vp*-2) rectangle (\hp*\r+\hp,\vp*-2+\vp);
}
\node at (\hp*0.5,\vp*0.5) {{\small $8 \ 9$}};
\node at (\hp*1.5,\vp*0.5) {{\small $10$}};
\node at (\hp*0.5,\vp*-0.5) {{\small $3 \ 4$}};
\node at (\hp*1.5,\vp*-0.5) {{\small $7$}};
\node at (\hp*2.5,\vp*-0.5) {{\small $12$}};
\node at (\hp*0.5,\vp*-1.5) {{\small $1 \ 2$}};
\node at (\hp*1.5,\vp*-1.5) {{\small $5$}};
\node at (\hp*2.5,\vp*-1.5) {{\small $6$}};
\node at (\hp*3.5,\vp*-1.5) {{\small $11$}};
\end{tikzpicture}$};

\node at (\hstep*3.2,0) (A3) {$g_{(5,4,3)}(T)=
\begin{tikzpicture}[baseline=-1mm]
\def \hpp {1.3em}
\def \vpp {1.1em}
\foreach \c in {0,-1,-2}{
    \foreach \r in {0,1}{
        \draw (\hp*\r,\vp*\c) rectangle (\hp*\r+\hp,\vp*\c+\vp);
    }
}
\draw (\hp*0,\vp*1) rectangle (\hp*1,\vp*2);
\draw (\hp*2,\vp*-1) rectangle (\hp*3,\vp*0);
\draw (\hp*2,\vp*-2) rectangle (\hp*3,\vp*-1);
\node at (\hp*0.5,\vp*1.5) {{\small $11$}};
\node at (\hp*0.5,\vp*0.5) {{\small $6$}};
\node at (\hp*1.5,\vp*0.5) {{\small $12$}};
\node at (\hp*0.5,\vp*-0.5) {{\small $5$}};
\node at (\hp*1.5,\vp*-0.5) {{\small $7$}};
\node at (\hp*2.5,\vp*-0.5) {{\small $10$}};
\node at (\hp*0.5,\vp*-1.5) {{\small $1$ $2$}};
\node at (\hp*1.5,\vp*-1.5) {{\small $3$ $4$}};
\node at (\hp*2.5,\vp*-1.5) {{\small $8$ $9$}};
\end{tikzpicture}$};

\draw[->] (A2) -- (A3) node[above,midway] {\tiny $\reflectf$};
\draw[->] (A1) -- (A2) node[above,midway] {\tiny $\mergef_{\rm c}$};
\end{tikzpicture}\, .
\]

\begin{proposition}\label{Prop: bijections among SPYCT SPYT and SVT }
Let $\lambda$ be a peak partition with its last part greater than $1$.
The map $g_\lambda$ is a bijection.
Furthermore, the diagram
\[ \begin{tikzcd}
\widehat{\SYCT(\lambda)} \arrow{r}{f_\lambda}  & \SYT(\lambda)  \arrow{r}{\reflectf}  & \SYT(\lambda') \arrow{d}{\mergef_\rmr} \\
\widehat{\SPYCT(\lambda)} \arrow{u}{\iota} \arrow{r}{f_\lambda|} & \SPYT(\lambda) \arrow{u}{\iota} \arrow{r}{g_\lambda} &  \SVT(\lambda'_-,\rho_{\lambda_1-1}) 
\end{tikzcd}
\]
is commutative, where $\iota$'s are the natural injections.
\end{proposition}
\begin{proof}
We prove the first assertion by directly constructing the inverse of $g_\lambda$. Let $U \in \SVT(\lambda'_{-}, \rho_{\lambda_1-1})$.
Define $\splitf_\rmr(U)$ as the tableau obtained by splitting the first column into two columns. During this process, the smaller entry from the $(1,i)$ box of $U$ is placed in the $(1,i)$ box of $\splitf_\rmr(U)$, while the larger entry is positioned in the $(2,i)$ box of $\splitf_\rmr(U)$.
Now, consider the map 
\[
h_\lambda:\SVT(\lambda'_{-}, \rho_{\lambda_1-1}) \rightarrow \SPYT(\lambda), \quad 
V \mapsto \splitf_{\rm c} \circ \reflectf (V).
\]
Given that $\reflectf$ is an involution and $\splitf_{\rm c} \circ \mergef_{\rm c}(T) = T$ for $T \in \SPYT(\lambda)$, it follows that $h_\lambda$ is the inverse of $g_\lambda$. 
And, the commutativity of the diagram follows from the equality $\reflectf \circ \mergef_{\rm c} = \mergef_\rmr \circ \reflectf$.
\end{proof}

Combining \cref{SPYCT for rectangular shaped peak compositions} with \cref{Other bijections}, we observe that the sets
\begin{align*}
\SPYCT((m^k)),
\quad \setWord{\SPYCT((m^k))},
\quad \SPYT((m^k)),
\quad \SVT((k^{m-1}),\rho_{m-1}), \quad (m,\, k \geq 2),  
\end{align*}
are mutually bijective.
The cardinality of these sets can be regarded as a novel generalization of the Fuss--Catalan number $C_k^{(3)}$ because it equals $C_k^{(3)}$ when $m=3$.


\vspace*{10mm}

\begin{thebibliography}{10}

\bibitem{18AHM}
E.~E. Allen, J.~Hallam, and S.~K. Mason.
\newblock {D}ual immaculate quasisymmetric functions expand positively into
  {Y}oung quasisymmetric {S}chur functions.
\newblock {\em J. Combin. Theory Ser. A}, 157:70--108, 2018.

\bibitem{14BBSSZ}
C.~Berg, N.~Bergeron, F.~Saliola, L.~Serrano, and M.~Zabrocki.
\newblock A lift of the {S}chur and {H}all--{L}ittlewood bases to
  non-commutative symmetric functions.
\newblock {\em Canad. J. Math.}, 66(3):525--565, 2014.

\bibitem{15BBSSZ}
C.~Berg, N.~Bergeron, F.~Saliola, L.~Serrano, and M.~Zabrocki.
\newblock Indecomposable modules for the dual immaculate basis of
  quasi-symmetric functions.
\newblock {\em Proc. Amer. Math. Soc.}, 143(3):991--1000, 2015.

\bibitem{17BBSSZ}
C.~Berg, N.~Bergeron, F.~Saliola, L.~Serrano, and M.~Zabrocki.
\newblock Multiplicative structures of the immaculate basis of non-commutative
  symmetric functions.
\newblock {\em J. Combin. Theory Ser. A}, 152:10--44, 2017.

\bibitem{04BHT}
N.~Bergeron, F.~Hivert, and J.-Y. Thibon.
\newblock The peak algebra and the {H}ecke--{C}lifford algebras at q= 0.
\newblock {\em J. Combin. Theory Ser. A}, 107(1):1--19, 2004.

\bibitem{16BSZ}
N.~Bergeron, J.~S{\'a}nchez-Ortega, and M.~Zabrocki.
\newblock The pieri rule for dual immaculate quasi-symmetric functions.
\newblock {\em Ann. Comb.}, 20:283--300, 2016.

\bibitem{11BLW}
C.~Bessenrodt, K.~Luoto, and S.~van Willigenburg.
\newblock Skew quasisymmetric {S}chur functions and noncommutative {S}chur
  functions.
\newblock {\em Adv. Math.}, 226(5):4492--4532, 2011.

\bibitem{16BTvW}
C.~Bessenrodt, V.~Tewari, and S.~van Willigenburg.
\newblock {L}ittlewood--{R}ichardson rules for symmetric skew quasisymmetric
  {S}chur functions.
\newblock {\em J. Combin. Theory Ser. A}, 137:179--206, 2016.

\bibitem{06BB}
A.~Bjorner and F.~Brenti.
\newblock {\em Combinatorics of {C}oxeter groups}, volume 231.
\newblock Springer Science \& Business Media, 2006.

\bibitem{91BW}
A.~Bj{\"o}rner and M.~L. Wachs.
\newblock Permutation statistics and linear extensions of posets.
\newblock {\em J. Combin. Theory Ser. A}, 58(1):85--114, 1991.

\bibitem{19CD}
C.~Ceballos and R.~S. Gonz\'{a}lez~D'Le\'{o}n.
\newblock Signature {C}atalan combinatorics.
\newblock {\em J. Comb.}, 10(4):725--773, 2019.

\bibitem{21CKNO}
S.-I. Choi, Y.-H. Kim, S.-Y. Nam, and Y.-T. Oh.
\newblock Modules of the $0$-{H}ecke algebra arising from standard permuted
  composition tableaux.
\newblock {\em J. Combin. Theory Ser. A}, 179:105389, 2021.

\bibitem{22CKNO2}
S.-I. Choi, Y.-H. Kim, S.-Y. Nam, and Y.-T. Oh.
\newblock Homological properties of 0-{H}ecke modules for dual immaculate
  quasisymmetric functions.
\newblock {\em Forum Math. Sigma}, 10:e91, 2022.

\bibitem{22CKNO1}
S.-I. Choi, Y.-H. Kim, S.-Y. Nam, and Y.-T. Oh.
\newblock The projective cover of tableau-cyclic indecomposable
  ${H}_n(0)$-modules.
\newblock {\em Trans. Amer. Math. Soc.}, 375(11):7747--7782, 2022.

\bibitem{24CKO}
S.-I. Choi, Y.-H. Kim, and Y.-T. Oh.
\newblock Poset modules of the 0-{H}ecke algebras and related quasisymmetric
  power sum expansions.
\newblock {\em Eur. J. Comb.}, 120:103965, 2024.

\bibitem{96DKLT}
G.~Duchamp, D.~Krob, B.~Leclerc, and J.-Y. Thibon.
\newblock Fonctions quasi-sym\'{e}triques, fonctions sym\'{e}triques non
  commutatives et alg\`ebres de {H}ecke \`a {$q=0$}.
\newblock {\em C. R. Acad. Sci. Paris S\'{e}r. I Math.}, 322(2):107--112, 1996.

\bibitem{11HLMW}
J.~Haglund, K.~Luoto, S.~Mason, and S.~van Willigenburg.
\newblock Quasisymmetric {S}chur functions.
\newblock {\em J. Combin. Theory Ser. A}, 118(2):463--490, 2011.

\bibitem{11HLMW2}
J.~Haglund, K.~Luoto, S.~Mason, and S.~van Willigenburg.
\newblock Refinements of the {L}ittlewood--{R}ichardson rule.
\newblock {\em Trans. Amer. Math. Soc.}, 363(3):1665--1686, 2011.

\bibitem{15JL}
N.~Jing and Y.~Li.
\newblock A lift of {S}chur's {$Q$}-functions to the peak algebra.
\newblock {\em J. Combin. Theory Ser. A}, 135:268--290, 2015.

\bibitem{22JKLO}
W.-S. Jung, Y.-H. Kim, S.-Y. Lee, and Y.-T. Oh.
\newblock Weak {B}ruhat interval modules of the 0-{H}ecke algebra.
\newblock {\em Math. Z.}, 301(4):3755--3786, 2022.

\bibitem{19Oguz}
E.~Kantarc{\i}~O{\u{g}}uz.
\newblock A note on {J}ing and {L}i’s type {B} quasisymmetric {S}chur
  functions.
\newblock {\em Ann. Comb.}, 23(1):159--170, 2019.

\bibitem{23KLO}
Y.-H. Kim, S.-Y. Lee, and Y.-T. Oh.
\newblock Regular {S}chur labeled skew shape posets and their $0$-{H}ecke
  modules.
\newblock {\em Forum Math. Sigma}, 12:e110, 2024.

\bibitem{97KT}
D.~Krob and J.-Y. Thibon.
\newblock Noncommutative symmetric functions. {IV}. {Q}uantum linear groups and
  {H}ecke algebras at {$q=0$}.
\newblock {\em J. Algebraic Combin.}, 6(4):339--376, 1997.

\bibitem{13LMvW}
K.~Luoto, S.~Mykytiuk, and S.~van Willigenburg.
\newblock {\em An introduction to quasisymmetric {S}chur functions}.
\newblock SpringerBriefs in Mathematics. Springer, New York, 2013.

\bibitem{08Mason}
S.~Mason.
\newblock A decomposition of {S}chur functions and an analogue of the
  {R}obinson--{S}chensted--{K}nuth algorithm.
\newblock {\em Sémin. Lothar. Comb.}, 75(B57e), 2006.

\bibitem{09Mason}
S.~Mason.
\newblock An explicit construction of type {A} {D}emazure atoms.
\newblock {\em J. Algebraic Combin.}, 29(3):295--313, 2009.

\bibitem{14MW}
P.~R. McNamara and R.~E. Ward.
\newblock Equality of {$P$}-partition generating functions.
\newblock {\em Ann. Comb.}, 18(3):489--514, 2014.

\bibitem{79Norton}
P.~N. Norton.
\newblock {$0$}-{H}ecke algebras.
\newblock {\em J. Austral. Math. Soc. Ser. A}, 27(3):337--357, 1979.

\bibitem{OEIS}
{OEIS Foundation Inc.}
\newblock The {O}n-{L}ine {E}ncyclopedia of {I}nteger {S}equences, 2024.
\newblock Published electronically at \url{http://oeis.org}.

\bibitem{22Searles}
D.~Searles.
\newblock Diagram supermodules for $0$-{H}ecke--{C}lifford algebras.
\newblock {\em Math. Z.}, 310(3):43, 2025.

\bibitem{24SS}
D.~Searles and M.~Slattery-Holmes.
\newblock {E}xpanding quasisymmetric {S}chur $q$-functions into peak {Y}oung
  quasisymmetric schur functions.
\newblock {\em arXiv preprint}, arXiv:2406.12751, 2024.

\bibitem{72Stanley}
R.~Stanley.
\newblock {\em Ordered structures and partitions}.
\newblock American Mathematical Society, Providence, R.I., 1972.
\newblock Memoirs of the American Mathematical Society, No. 119.

\bibitem{99Stanley}
R.~Stanley.
\newblock {\em Enumerative combinatorics. {V}ol. 2}, volume~62 of {\em
  Cambridge Studies in Advanced Mathematics}.
\newblock Cambridge University Press, Cambridge, 1999.

\bibitem{15Stan}
R.~P. Stanley.
\newblock {\em Catalan numbers}.
\newblock Cambridge University Press, New York, 2015.

\bibitem{97Stem}
J.~Stembridge.
\newblock Enriched $p$-partitions.
\newblock {\em Trans. Amer. Math. Soc.}, 349(2):763--788, 1997.

\bibitem{89Stem}
J.~R. Stembridge.
\newblock Shifted tableaux and the projective representations of symmetric
  groups.
\newblock {\em Adv. Math.}, 74(1):87--134, 1989.

\bibitem{15TW}
V.~Tewari and S.~van Willigenburg.
\newblock Modules of the {$0$}-{H}ecke algebra and quasisymmetric {S}chur
  functions.
\newblock {\em Adv. Math.}, 285:1025--1065, 2015.

\bibitem{19TW}
V.~Tewari and S.~van Willigenburg.
\newblock Permuted composition tableaux, {$0$}-{H}ecke algebra and labeled
  binary trees.
\newblock {\em J. Combin. Theory Ser. A}, 161:420--452, 2019.

\end{thebibliography}
\end{document}